\documentclass[pdflatex,sn-basic]{sn-jnl}

\usepackage{graphicx}
\usepackage{amsmath,amssymb,bm}
\usepackage{booktabs,longtable,tabularx}
\usepackage{placeins}

\begin{document}
\title[Robust weighting for GNSS ambiguity resolution]{Robust weighting for GNSS integer ambiguity resolution under heavy-tailed errors and localized biases}
\author[1]{\fnm{Ruifeng} \sur{Zhang}}\email{zhangrf24@m.fudan.edu.cn}
\author*[1,2]{\fnm{Yunxin} \sur{Zhang}}\email{xyz@fudan.edu.cn}
\affil[1]{\orgdiv{School of Mathematical Sciences}, \orgname{Fudan University}, \orgaddress{\city{Shanghai}, \postcode{200433}, \country{China}}}
\affil[2]{\orgdiv{Shanghai Key Laboratory for Contemporary Applied Mathematics, Laboratory of Mathematics for Nonlinear Science}, \orgname{Fudan University}, \orgaddress{\city{Shanghai}, \postcode{200433}, \country{China}}}
\abstract{Precise baseline estimation in carrier-phase positioning with global navigation satellite systems (GNSS) depends on resolving the integer ambiguity vector. Conventional integer least squares derives its float solution and search metric from nominal weighted least squares (WLS). Heavy-tailed errors and localized biases can distort the float estimate under fixed nominal weights, so we compare weighting rules under a common integer search. Paired simulations compare WLS, four iteratively reweighted least-squares (IRLS) rules and regularized Gaussian-mixture weighting across error laws, satellite geometries and bias durations. Adaptive weighting improves full-vector recovery over WLS at every tested duration under Gaussian-mixture or Cauchy errors. WLS recovers more vectors in some longer-bias cases with Gaussian or Student-$t_3$ background noise. In 493 ten-epoch windows from one measured short baseline on two dates, we add sparse Gaussian-mixture errors in original-order Cholesky-whitened coordinates. Huber and Cauchy IRLS match the internal reference vector in 96.3\% and 96.1\% of windows, compared with 9.7\% for WLS and 92.1\% for iterative data snooping followed by integer least squares. Before screening, the 95th-percentile baseline error relative to the header-coordinate benchmark is 6.3 cm for Huber and 198.0 cm for WLS. Further tests show that whitening can alter IRLS recovery and acceptance. Under frozen cross-date thresholds, methods can rank differently by recovery and by the fraction of windows yielding accepted baselines within 10 cm of the benchmark. For short-baseline GNSS positioning, weighting rules should therefore be evaluated with the intended whitening and screening procedures to determine how often they yield accurate fixed baselines.}
\AddToHookNext{cmd/Titlefont/before}{\vskip-55pt\relax}
\maketitle

\section{Introduction}\label{sec:intro}
The carrier-phase GNSS float solution treats the ambiguities as real-valued. Fixing them to integers allows the baseline to be re-estimated conditional on that vector \citep{groves2013,leick2015,teunissen2017}. A single wrong integer component can degrade the fitted baseline.

In integer least squares (ILS), observation weights determine both the float ambiguity estimate and the quadratic metric used to rank integer candidates. LAMBDA and MLAMBDA reduce and search this quadratic problem efficiently \citep{teunissen1995,chang2005}. ILS maximizes the success probability in its usual class when the float ambiguity is Gaussian and correctly weighted \citep{teunissen1999,hassibi1998}. Heavy tails, multipath and localized biases can invalidate that stochastic description even when the search solves its stated problem exactly \citep{hsu2018,zhu2018}. Changing observation weights can reduce the influence of outlying residuals and can also reduce information in the fit. We use full-vector recovery, or recovery, for agreement of all estimated ambiguities with the simulated truth or the internal long-arc reference constructed in Section~\ref{sec:realdata}. We examine how recovery changes and how these changes carry through to baseline error and acceptance.

Robust equivalent weights and phase models have long been used in geodesy \citep{yang1994,yang2002,wieser2000,knight2009}. Recent GNSS strategies combine robust residuals with satellite selection and partial ambiguity resolution \citep{par2025}, or combine heavy-tailed distributions, candidate estimators and fixing forms \citep{zhan2026}. Bias-constrained ILS incorporates bounded real biases \citep{khodabandeh2024}, while the ambiguity-resolved detector addresses validation of a mixed-integer model \citep{teunissen2024}. These studies motivate comparisons across estimation and decision strategies. We hold the integer-search routine and the full-vector recovery criterion fixed to compare how observation processing changes the solution on paired inputs.

Related theory addresses other estimator classes or assumptions. Best integer equivariant estimation under elliptical distributions uses a mean-squared-error criterion \citep{teunissen2020}. Generalized ILS success-rate bounds average a Gaussian bound over a scalar mixing scale \citep{teunissen2026}. The adaptive estimators studied here determine weights from residuals in the working observation basis. The distribution of the resulting float ambiguities affects the validity of a calibrated ratio threshold \citep{teunissenverhagen2009,verhagen2013}.

Nominal weighted least squares (WLS) is the baseline. The four iteratively reweighted least-squares (IRLS) configurations use Huber, Cauchy, IGG~III and Tukey weights with a common initialization, residual-scale and update protocol (Section~\ref{sec:irls}). The auxiliary Gaussian mixture model (GMM) fits regularized full component covariance matrices (Section~\ref{sec:gmm}). The main weighting tables show GMM results for $K=2$. Tables~S1--S4 in the Supplementary Information give simulation results for $K=3$ and $K=4$. Both families supply the working metric to the same integer search.

We compare estimators on identical inputs within each design cell. The controlled random and between-satellite simulations vary the error law, geometry and local-bias duration to identify where adaptive weighting improves full-vector recovery. The semi-synthetic study adds controlled errors to measured short-baseline observations and includes iterative data snooping followed by integer least squares (IDS--ILS) as a sequential residual-screening comparator. The internal long-arc reference supplies the evaluation integers for those windows.

Using paired inputs and a common integer search, we identify when alternative observation-weighting rules improve full-vector recovery and assess how those gains translate into baseline error before screening and accepted-solution quality under frozen cross-date thresholds. Full-vector recovery is the primary endpoint. Baseline error, the quality and frequency of accepted solutions under frozen thresholds, and computation time measure further consequences of the estimator choice. The whitening comparison uses the same contaminated raw observations to assess sensitivity to the representation supplied to the entire estimator. Section~\ref{sec:model} defines the estimators and endpoints, Section~\ref{sec:simulations} reports the simulations, Section~\ref{sec:realdata} reports recovery and positioning on measured observations, and Section~\ref{sec:conditions} examines acceptance and applicability.
\section{Model and estimation methods}\label{sec:model}
All comparisons estimate a shared real parameter vector and the full integer ambiguity vector from several observation epochs. This section first defines the weighted estimation problem and its integer search, then explains how each method constructs the weights or selects observations. Acceptance is a separate decision applied after the candidates have been obtained.

\subsection{Observation model and integer search}\label{sec:common-model}
One problem contains $N$ measurement epochs, indexed by $j=1,\ldots,N$. Each epoch supplies $m$ observation components in a consistent order. The real parameter vector $\bm x\in\mathbb R^p$ and integer ambiguity vector $\bm z\in\mathbb Z^{n_z}$ are shared across epochs within the problem. The model is
\begin{equation}
\bm y_j=\bm A_j\bm x+\bm B_j\bm z+\bm e_j,
\qquad \bm x\in\mathbb R^p,\quad \bm z\in\mathbb Z^{n_z}.
\label{eq:model}
\end{equation}
Here $\bm y_j$ and $\bm e_j$ have $m$ entries, $\bm A_j$ is $m\times p$ and $\bm B_j$ is $m\times n_z$. The common model leaves the error law unspecified. Geometry can change with $j$. Before screening, the observation components, their order and the ambiguity dimension are fixed within a problem. Their numbers may differ across problems.

A weighting method supplies a positive-definite working weight matrix $\bm W_j$ for each epoch. With these matrices held fixed, define the residual and weighted least-squares objective by
\begin{equation}
\bm r_j(\bm x,\bm z)=\bm y_j-\bm A_j\bm x-\bm B_j\bm z,
\qquad
J(\bm x,\bm z)=\sum_{j=1}^{N}\bm r_j^{\mathsf T}\bm W_j\bm r_j.
\label{eq:joint}
\end{equation}
Each solve holds the working weights fixed and minimizes $J$ over real $\bm x$ and integer $\bm z$. For diagonal $\bm W_j$, each diagonal entry multiplies the squared residual of its component. WLS uses nominal weights. IRLS and GMM compute a new working metric for each update. The corresponding solves incorporate the added equations in~\eqref{eq:mean-floor} or the rows retained by IDS--ILS.

To describe this solver, stack the observations into $\bm y_{\rm st}=(\bm y_1^{\mathsf T},\ldots,\bm y_N^{\mathsf T})^{\mathsf T}$ and let $\bm F$ have block rows $[\bm A_j\ \bm B_j]$. Its columns correspond to $\bm x$ followed by $\bm z$. Let $\bm W$ be block diagonal with blocks $\bm W_1,\ldots,\bm W_N$. The stacked quantities $\bm y_{\rm st},\bm F,\bm W$ have dimensions $Nm\times1$, $Nm\times(p+n_z)$ and $Nm\times Nm$, respectively. The normal matrix and right-hand side are
\begin{equation}
\bm N=\bm F^{\mathsf T}\bm W\bm F
=\begin{bmatrix}\bm N_{xx}&\bm N_{xz}\\\bm N_{zx}&\bm N_{zz}\end{bmatrix},
\qquad
\bm h=\bm F^{\mathsf T}\bm W\bm y_{\rm st}
=\begin{bmatrix}\bm h_x\\\bm h_z\end{bmatrix}.
\label{eq:normal}
\end{equation}
Bold $\bm N$ denotes the $(p+n_z)\times(p+n_z)$ normal matrix, whereas $N$ counts epochs. The subscripts $x,z$ index the real and integer blocks. The diagonal blocks have orders $p$ and $n_z$, and $\bm N_{xz}=\bm N_{zx}^{\mathsf T}$. The vectors $\bm h_x,\bm h_z$ have lengths $p,n_z$.

First let $\bm z$ vary over $\mathbb R^{n_z}$. Write $(\hat{\bm x}_{\rm f},\hat{\bm z})$ for the joint float solution. Its two block normal equations are
\begin{equation}
\bm N_{xx}\hat{\bm x}_{\rm f}+\bm N_{xz}\hat{\bm z}=\bm h_x,
\qquad
\bm N_{zx}\hat{\bm x}_{\rm f}+\bm N_{zz}\hat{\bm z}=\bm h_z.
\label{eq:block-normal}
\end{equation}
For a full-column-rank design and positive-definite weights, solving the first equation for $\hat{\bm x}_{\rm f}$ and substituting it into the second gives
\begin{equation}
\bm Q_z^{-1}=\bm N_{zz}-\bm N_{zx}\bm N_{xx}^{-1}\bm N_{xz},
\qquad
\hat{\bm z}=\bm Q_z(\bm h_z-\bm N_{zx}\bm N_{xx}^{-1}\bm h_x).
\label{eq:schur}
\end{equation}
The Schur complement $\bm Q_z^{-1}$ is the ambiguity information matrix for this fit. Its inverse $\bm Q_z$ is the $n_z\times n_z$ working covariance. For a specified integer vector $\bm z$, the real-parameter estimate conditional on that vector follows from the first block equation:
\begin{equation}
\hat{\bm x}(\bm z)=\bm N_{xx}^{-1}(\bm h_x-\bm N_{xz}\bm z).
\label{eq:fixed-real}
\end{equation}
Substituting this fit into $J$ leaves the candidate-independent float objective $J(\hat{\bm x}_{\rm f},\hat{\bm z})$ plus the ambiguity distance. The distance ranks candidates for the acceptance ratio in~\eqref{eq:ratio}. The chosen integer vector is a minimizer of this distance:
\begin{equation}
\check{\bm z}\in\operatorname*{arg\,min}_{\bm z\in\mathbb Z^{n_z}}
(\bm z-\hat{\bm z})^{\mathsf T}\bm Q_z^{-1}(\bm z-\hat{\bm z}).
\label{eq:ils-obj}
\end{equation}
When the fixed weight matrix equals the inverse covariance of the stacked observation errors, $\bm Q_z$ is the sampling covariance of $\hat{\bm z}$. With adaptive weights, $\bm Q_z^{-1}$ specifies the working search metric and the sampling covariance requires separate analysis. For the augmented IRLS fit, the normal equations and float objective include the added equations.

All methods use MLAMBDA-related candidate generation \citep{chang2005} followed by a QR-based Schnorr--Euchner sphere search in the original basis. Candidate distances are recomputed in the original working metric. A completed sphere search solves~\eqref{eq:ils-obj} to numerical precision. The sphere search has a 20,000-node budget per call, with candidate generation outside this count. An exhausted search is recorded as incomplete. The corresponding planned problem counts as a recovery failure and is ineligible for acceptance.

The eight weighting configurations are WLS, four IRLS rules and GMM at $K=2,3,4$. The principal weighting tables show WLS, the four IRLS rules and GMM at $K=2$. Experiments that also include IDS--ILS compare nine configurations.

\medskip
\noindent\textbf{Main notation.}\par\smallskip
{\small\renewcommand{\arraystretch}{0.98}
\noindent\begin{tabularx}{\linewidth}{@{}l>{\raggedright\arraybackslash}X@{}}
\toprule
Symbol & Meaning \\
\midrule
$N,m,p,n_z$ & Epochs, components per epoch, real parameters and integer ambiguities \\
$\bm y_j,\bm e_j,\bm A_j,\bm B_j$ & Observation and error vectors, real and integer designs \\
$\bm x,\bm z,\bm z_0,\bm z_{\rm ref}$ & Real and integer unknowns, simulated truth and internal reference \\
$\hat{\bm x}_{\rm f},\hat{\bm z},\check{\bm z}$ & Real-parameter float estimate, float ambiguity and chosen integer vector \\
$\bm W_j,\bm W,\bm N,\bm h,\bm Q_z$ & Epoch and stacked weights, normal matrix, right-hand side and working ambiguity covariance \\
$\bm\Sigma_j,\bm L_j,\hat\sigma^{(t)},\eta_{ij}^{(t)}$ & Nominal covariance, lower Cholesky factor, robust scale and standardized residual magnitude \\
$a_i,b_i,\tilde a_i,\delta_i$ & Total weight, weighted observation sum, bounded total weight and added weight \\
$K,\pi_k,\gamma_{jk},\nu_k,\bm\Sigma_k$ & Component count, mixing weight, responsibility, effective size and component covariance \\
$J_i,R,\tau$ & Candidate distance, guarded distance ratio and acceptance threshold \\
$n,n_a,n_w,n_{10},n_{a,10}$ & Planned problems, accepted problems, wrong accepted vectors, solutions within 10\,cm and accepted solutions within 10\,cm \\
$A,q,u,F_{10},V_{10}$ & Acceptance, conditional and unconditional wrong-fix frequencies, and the two all-planned position fractions \\
$\xi_{ij}$ & Binary contamination state of satellite $i$ at epoch $j$ \\
$\bm c_j,c_{ij}$ & Binary injection mask and its components \\
$\bm T_j^{\rm rev},\bm T_j^{\rm sym}$ & Transforms to reversed-order Cholesky and symmetric whitening \\
$\bm H,\hat{\bm v}$ & Hat matrix and residuals of the screening fit \\
\bottomrule
\end{tabularx}\par}

\subsection{WLS and IRLS weighting}\label{sec:irls}
WLS uses $\bm W_j=\bm\Sigma_j^{-1}$ for the nominal epoch covariance $\bm\Sigma_j$. In the random benchmark, $\bm\Sigma_j=\sigma^2\bm I_m$ for every generating law (Section~\ref{sec:random}). In the GNSS comparisons, nominal whitening transforms the observations and both designs so that the working covariance becomes the identity. WLS then uses identity weights. Sections~\ref{sec:gnss} and~\ref{sec:realdata} specify the transforms.

IRLS repeatedly assigns weights to the current residuals and solves the resulting weighted problem \citep{huber2009}. The initial pair $(\bm x^{(0)},\bm z^{(0)})$ comes from a unit-weight ILS fit to a vector containing the median of each observation component across epochs. The corresponding designs are the epoch means $N^{-1}\sum_j\bm A_j$ and $N^{-1}\sum_j\bm B_j$. Subsequent fits use the original epoch-specific designs. GMM uses the same initialization.

At outer iteration $t$, write the current estimate as $(\bm x^{(t)},\bm z^{(t)})$. For $i=1,\ldots,m$ and $j=1,\ldots,N$, the residual in component $i$ of epoch $j$ is
\begin{equation}
r_{ij}^{(t)}=y_{ij}-\bm a_{ij}^{\mathsf T}\bm x^{(t)}-\bm b_{ij}^{\mathsf T}\bm z^{(t)},
\label{eq:irls-residual}
\end{equation}
where $\bm a_{ij}^{\mathsf T}$ and $\bm b_{ij}^{\mathsf T}$ are row $i$ of $\bm A_j$ and $\bm B_j$, respectively, and $y_{ij}$ is component $i$ of $\bm y_j$. Pooling all $mN$ residual entries, define the median absolute deviation and scale by
\begin{equation}
\mathrm{MAD}^{(t)}=\operatorname{median}_{i,j}\left|r_{ij}^{(t)}-\operatorname{median}_{i',j'}r_{i'j'}^{(t)}\right|,
\qquad \hat\sigma^{(t)}=\max(1.4826\,\mathrm{MAD}^{(t)},10^{-6}).
\label{eq:irls-scale}
\end{equation}
The primed indices range over the same residual entries. The inner median centres the residuals for the scale calculation. The normal-consistency factor $1/\Phi^{-1}(0.75)$ is approximated by 1.4826, where $\Phi$ is the standard normal distribution function. The lower bound $10^{-6}$ keeps the scale positive. The weight argument is the original residual magnitude $\eta_{ij}^{(t)}=|r_{ij}^{(t)}|/\hat\sigma^{(t)}$. Each rule maps a non-negative standardized residual magnitude $\eta$ to a nominal weight in $[0,1]$.

\noindent\textbf{Huber IRLS} keeps small residuals at unit weight and reduces the weight of large residuals:
\begin{equation}
w_{\rm H}(\eta)=\begin{cases}1,&\eta\le1.345,\\1.345/\eta,&\eta>1.345.\end{cases}
\label{eq:huber}
\end{equation}
\noindent\textbf{Cauchy IRLS} reduces the weight continuously without setting it to zero:
\begin{equation}
w_{\rm C}(\eta)=\frac{1}{1+(\eta/2.3849)^2}.
\label{eq:cauchy}
\end{equation}
\noindent\textbf{IGG III IRLS} reduces intermediate residuals and gives zero nominal weight beyond its upper cutoff \citep{yang1994,yang2002}:
\begin{equation}
w_{\rm I}(\eta)=\begin{cases}
1,&\eta\le1.5,\\
\dfrac{1.5(3-\eta)^2}{\eta(3-1.5)^2},&1.5<\eta<3,\\[0.4em]
0,&\eta\ge3.
\end{cases}
\label{eq:igg}
\end{equation}
\noindent\textbf{Tukey IRLS} also has a finite cutoff:
\begin{equation}
w_{\rm T}(\eta)=\begin{cases}(1-(\eta/4.685)^2)^2,&\eta\le4.685,\\0,&\eta>4.685.\end{cases}
\label{eq:tukey}
\end{equation}
The constants 1.345, 2.3849 and 4.685 use conventional Gaussian-scale tuning for Huber, Cauchy and Tukey. IGG~III uses the cutoffs 1.5 and 3. The weight function is selected separately from the error-generating law. Supplementary Tables~S29--S30 give sensitivity results under the earlier MAD scaling and row-centre conventions.

For the weights actually used in iteration $t+1$, set $w_{ij}^{(t)}=\max(w(\eta_{ij}^{(t)}),10^{-12})$. The diagonal entries of $\bm W_j^{(t)}$ are these weights. A second lower bound applies to the total weight of each component across epochs. In the fixed-weight step, $a_i$ denotes this total, $b_i$ the weighted observation sum and $\tilde a_i$ the bounded total. The index $\ell$ runs over the $m$ components. Each applied weight is positive, so $a_i>0$. Suppressing the iteration superscript,
\begin{equation}
a_i=\sum_{j=1}^N w_{ij},\qquad b_i=\sum_{j=1}^N w_{ij}y_{ij},\qquad
\tilde a_i=\max(a_i,10^{-3}\max_\ell a_\ell,10^{-8}).
\label{eq:floor}
\end{equation}
The rule applies to all four IRLS configurations. An equation is added when $\tilde a_i>a_i$. Let $\bm\theta=(\bm x^{\mathsf T},\bm z^{\mathsf T})^{\mathsf T}$ and let $\bm d_{ij}^{\mathsf T}=[\bm a_{ij}^{\mathsf T}\ \bm b_{ij}^{\mathsf T}]$. The weighted means are $\bar y_i=b_i/a_i$ and $\bar{\bm d}_i=a_i^{-1}\sum_jw_{ij}\bm d_{ij}$, and the added weight is $\delta_i=\tilde a_i-a_i$. With $J(\bm\theta)$ denoting the fixed-weight objective in~\eqref{eq:joint}, the next fit minimizes
\begin{equation}
J_{\rm mean}(\bm\theta)=J(\bm\theta)+\sum_{i=1}^m\delta_i(\bar y_i-\bar{\bm d}_i^{\mathsf T}\bm\theta)^2.
\label{eq:mean-floor}
\end{equation}
When $\delta_i>0$, the fit gains one equation with response $\bar y_i$, design row $\bar{\bm d}_i^{\mathsf T}$ and weight $\delta_i$. The retained epoch equations and the added equation have the same weighted observation and design centres as the original weighted epoch equations. Both centres use the original total $a_i$. We call this mean-preserving augmentation. Supplementary Sections 7.5 and 7.6 give the fixed-geometry compression and the augmentation identity, respectively.

After each fit, IRLS and GMM check whether the Euclidean change in the real-parameter estimate is below $10^{-6}$ and the integer vector is unchanged. They stop when both conditions hold or after 20 updates. If another update is required, its weights are formed from the new residuals. Fits reaching the iteration budget remain in the results. These criteria measure parameter stabilization within a finite budget. Supplementary Section 7.7 reports stopping diagnostics.

\subsection{Regularized Gaussian mixture weighting}\label{sec:gmm}
\begin{samepage}
Let $\bm r_j=\bm r_j(\bm x,\bm z)$ be the current $m$-vector of residuals at epoch $j$. The model has $K$ zero-mean Gaussian components. The non-negative mixing weights $\pi_k$ and component covariances $\bm\Sigma_k\in\mathbb R^{m\times m}$ are shared across epochs. The full covariance matrices describe the joint residual scatter within each component. The responsibility is
\begin{equation}
\gamma_{jk}=\frac{\pi_k\varphi_m(\bm r_j;\bm0,\bm\Sigma_k)}{\sum_{\ell=1}^K\pi_\ell\varphi_m(\bm r_j;\bm0,\bm\Sigma_\ell)},
\label{eq:responsibility}
\end{equation}
Here $\varphi_m$ is the $m$-variate Gaussian density and $\gamma_{jk}$ is the conditional responsibility of component $k$ for residual $\bm r_j$ \citep{dempster1977,mclachlan2000}. When the denominator is positive, $\sum_k\gamma_{jk}=1$.
\end{samepage}

Let $\nu_k=\sum_{j=1}^N\gamma_{jk}$ be the effective size of component $k$. For $k\in\mathcal K_{\rm act}=\{k:\nu_k\ge10^{-8}\}$, set $\pi_k=\nu_k/N$ and compute the zero-centred scatter $\bm S_k=\nu_k^{-1}\sum_j\gamma_{jk}\bm r_j\bm r_j^{\mathsf T}$. Initialize the mixing weights at $1/K$. Set inactive component weights to zero and omit those components from the precision sum. After deactivation, the mixing weights can sum to less than one. Equation~\eqref{eq:responsibility} divides each weighted density by the sum of the weighted densities. Rescaling the active mixing weights to sum to one therefore gives the same responsibilities when this sum of weighted densities is positive. The covariance update adds a positive ridge:
\begin{equation}
\bm\Sigma_k^{\rm new}=\bm S_k+10^{-4}\max\!\left(\frac{\operatorname{tr}(\bm S_k)}m,10^{-12}\right)\bm I_m,
\qquad
\bm W_j=\sum_{k\in\mathcal K_{\rm act}}\gamma_{jk}(\bm\Sigma_k^{\rm new})^{-1}.
\label{eq:gmm-ridge}
\end{equation}
The ridge coefficient $10^{-4}$ multiplies the mean scatter variance after applying a positive lower bound to that variance scale. The precision sum can have off-diagonal entries. We first compute responsibilities from the current residuals and covariances, then update the covariances and use the same responsibilities with the updated precisions in the next fit. The initial covariance of component $k$ is $(2k-1)\bm I_m$ in squared working-observation units. Supplementary Section 7.2 gives numerical fallback rules, including those for zero total precision or failed matrix inversion.

The procedure alternates regularized scatter updates and integer refits under the stopping rule in Section~\ref{sec:irls}. With $N$ residual vectors, the unregularized scatter rank is at most $\min(N,m)$. At $N=5$ and $m=16$, at least eleven directions obtain their variance entirely from the ridge. The results compare $K=2,3,4$ within this implementation. Tables~S9--S10 describe the sensitivity to the ridge and initialization on selected short records.

\subsection{Iterative data snooping}\label{sec:ids}
IDS--ILS is the sequential residual-screening comparator. Let $\bm F$ and $\bm y_{\rm st}$ denote the currently retained nominally whitened design and observation rows. With all $p+n_z$ parameters treated as real in the screening fit, define the fitted-value projector and residual as follows:
\begin{equation}
\bm H=\bm F(\bm F^{\mathsf T}\bm F)^{-1}\bm F^{\mathsf T},
\qquad \hat{\bm v}=(\bm I-\bm H)\bm y_{\rm st}.
\label{eq:ids-projection}
\end{equation}
For $\ell_{\rm rows}$ retained rows, $\bm H$ is $\ell_{\rm rows}\times\ell_{\rm rows}$ and $\hat{\bm v}$ has $\ell_{\rm rows}$ entries. For retained row $i$, $H_{ii}$ is its leverage and $1-H_{ii}$ is its residual variance under the nominal unit-covariance model. Evaluate $|\hat v_i|/\sqrt{1-H_{ii}}$ for rows with $1-H_{ii}>10^{-10}$. Remove the row with the largest statistic if it exceeds 3.29 and refit the continuous model \citep{baarda1968,knight2009}. Once the testable statistics are all at or below 3.29, apply the common ILS search to the retained rows. The procedure returns no integer fit if it has no testable residual or if the current fit or a proposed deletion violates full column rank or leaves fewer than one residual degree of freedom. Such problems count as recovery failures and are ineligible for acceptance. The nominal two-sided probability at 3.29 is approximately 0.001 for a single standardized Gaussian residual. The method screens whitened rows and retains the original integer dimension.

\subsection{Candidate acceptance and general endpoints}\label{sec:acceptance-protocol}
Let $\bm z_0$ be the simulated true vector and $\bm z_{\rm ref}$ the internal reference for a measured window. A problem counts as recovered when a completed estimate equals the applicable vector in every component. The denominator includes every planned problem, with incomplete searches counted as failures. Acceptance uses the final candidates and a calibrated threshold.

Freeze the metric of the retained fit, including any augmented equations or retained rows. Rank distinct candidates $\bm z_1,\bm z_2$ by the quadratic term in~\eqref{eq:ils-obj}, with $\bm z_1=\check{\bm z}$. Let $J_i$ denote that term for candidate $i=1,2$. Define
\begin{equation}
J_i=(\bm z_i-\hat{\bm z})^{\mathsf T}\bm Q_z^{-1}(\bm z_i-\hat{\bm z}),
\qquad R=\frac{J_2}{\max(J_1,10^{-15})}.
\label{eq:ratio}
\end{equation}
The ratio uses the second-to-first ambiguity-distance convention of \citet{teunissenverhagen2009}. The candidates are ranked by~\eqref{eq:ils-obj}, so $J_1\le J_2$. When $J_1\ge10^{-15}$, the ratio is $J_2/J_1\ge1$. For smaller $J_1$, use the guarded ratio in~\eqref{eq:ratio}. A valid ratio requires a completed search with both candidates and a finite $R$. Accept the first candidate when the ratio is valid and $R\ge\tau$.

\begin{samepage}
Among $n$ planned evaluation problems, $n_a$ have accepted solutions and $n_w$ have accepted vectors with at least one component differing from the applicable reference. The acceptance fraction, conditional wrong-fix frequency and unconditional wrong-fix frequency are
\begin{equation}
A=\frac{n_a}{n},\qquad q=\frac{n_w}{n_a}\quad(n_a>0),\qquad u=\frac{n_w}{n}.
\label{eq:acceptance-rates}
\end{equation}
When $n_a>0$, $u=Aq$. At zero acceptance, $A=u=0$ and $q$ is undefined. Table columns marked \% display the corresponding proportions multiplied by 100. The supplementary tables also report $u$. Section~\ref{sec:realdata} defines the position endpoints.
\end{samepage}

Thresholds selected on calibration inputs are frozen for evaluation (Section~\ref{sec:conditions}). Supplementary Sections 7.2 and 7.5 describe the working metric and compression in fixed-geometry cases.

\section{Controlled simulations}\label{sec:simulations}
Each design cell fixes a combination of simulation factors. A comparison is paired when all configurations use the same generated problem and its known integer truth. The random benchmark tests the general mixed-integer model. The between-satellite experiment uses a three-dimensional simulated baseline and varies geometry and localized bias.

\subsection{Random mixed-integer benchmark}\label{sec:random}\label{sec:noise}
Each random problem has $N=40$ epochs, $m=13$ observations per epoch, $p=4$ real parameters and $n_z=6$ ambiguities. The designs are fixed across epochs within a problem and their entries are independent standard normal draws. The components of the true real vector $\bm x_0$ are independent $\mathcal N(0,10^2)$, and the true integer vector $\bm z_0$ is uniform on $\{-10,\ldots,9\}^{n_z}$. For each error law, all eight weighting configurations analyse the same 1,000 generated problems. The dimensions and nominal scale were chosen during exploratory work. The reported trials use this fixed design.

Let $\varepsilon$ denote a scalar observation error and let G, M, T and C index Gaussian, Gaussian-mixture, Student-$t_3$ and Cauchy laws. These four data-generating conditions are crossed with the weighting configurations; a method name such as Cauchy IRLS identifies a weight function, whereas Cauchy noise identifies the distribution of the simulated errors. For a nominal observation scale $\sigma$, the four laws are
\begin{align}
\varepsilon_{\mathrm G}&\sim \mathcal N(0,\sigma^2), &
\varepsilon_{\mathrm M}&\sim0.8\mathcal N(0,\sigma^2)+0.2\mathcal N(0,(10\sigma)^2),\nonumber\\
\varepsilon_{\mathrm T}&=a_T T_3,\quad T_3\sim t_3, &
\varepsilon_{\mathrm C}&\sim\operatorname{Cauchy}(0,a_C),
\label{eq:noise-laws}
\end{align}
with
\begin{equation}
a_T=\frac{\sigma\Phi^{-1}(0.75)}{F_{t_3}^{-1}(0.75)},
\qquad a_C=\sigma\Phi^{-1}(0.75).
\label{eq:noise-scales}
\end{equation}
The Gaussian-mixture law draws from its narrow component with probability 0.8 and its wide component with probability 0.2. Observation-level draws are independent. Let $\Phi$ and $F_{t_3}$ be the standard Gaussian and Student-$t_3$ distribution functions. For these continuous symmetric laws, the 0.75 quantile equals the median absolute deviation from zero. At the nominal scale $\sigma=2$, \eqref{eq:noise-scales} matches the Student-$t_3$ and Cauchy MADs to the Gaussian MAD. The Gaussian-mixture variance is $20.8\sigma^2$. The Cauchy law has no finite mean or variance.

Under Gaussian errors, WLS has the highest recovery at 98.4\% (Table~\ref{tab:core}). Under Gaussian-mixture errors, WLS recovers 4.1\%, the four IRLS configurations recover from 89.8\% to 94.3\%, and GMM at $K=2$ recovers 75.6\%. Under Cauchy errors the corresponding results are 0.2\%, 97.8\% to 98.2\%, and 49.0\%. Table~S1 gives all component counts. Tables~S29--S30 report sensitivities under historical settings, including use of the unscaled MAD and an earlier convention for computing the weighted row centres (Supplementary Section 7.6).
\begin{table}[htbp]
\caption{Full-vector recovery (\%) in the random mixed-integer benchmark. Each error law has 1,000 paired problems with 40 epochs per problem.}\label{tab:core}
\centering\small
\renewcommand{\arraystretch}{1.12}
\setlength{\tabcolsep}{2pt}
\begin{tabular*}{\textwidth}{@{\extracolsep{\fill}}lrrrr@{}}\toprule
Method & Gaussian & Student-$t_3$ & Gaussian mixture & Cauchy \\
\midrule
WLS & 98.4 & 85.9 & 4.1 & 0.2 \\
Huber IRLS & 97.7 & 97.8 & 89.8 & 97.8 \\
Cauchy IRLS & 97.2 & 97.4 & 90.5 & 98.1 \\
IGG~III IRLS & 96.3 & 97.3 & 94.3 & 98.2 \\
Tukey IRLS & 97.5 & 97.7 & 94.0 & 98.1 \\
GMM, $K=2$ & 95.6 & 96.9 & 75.6 & 49.0 \\
\bottomrule\end{tabular*}
\par\smallskip\parbox{\textwidth}{\footnotesize Columns identify the data-generating error laws. Recovery requires agreement with simulated truth in every integer component. The denominator includes all planned problems, with incomplete searches counted as failures. WLS denotes weighted least squares, IRLS iteratively reweighted least squares, and GMM a Gaussian mixture model with $K$ fitted components. IRLS uses the common scaled-MAD protocol. Table~S1 gives the results for $K=3,4$.}
\end{table}

\subsection{Between-satellite model and localized bias}\label{sec:gnss}
We next generate an idealized linear between-satellite model. Its real parameter $\bm x\in\mathbb R^3$ is a simulated baseline expressed in the same local Cartesian frame as the satellite directions. The measured-observation experiment in Section~\ref{sec:realdata} starts from two-receiver observations and forms receiver and satellite differences.

Each problem has nine satellites. Satellite $s$ has azimuth $\alpha_s$, elevation $e_s$ and unit direction $\bm u_s=(\cos e_s\cos\alpha_s,\cos e_s\sin\alpha_s,\sin e_s)^{\mathsf T}$. The $9\times3$ matrix $\bm U$ stacks these directions. The highest-elevation satellite is the reference. Each row of $\bm D$ places $+1$ in a target-satellite column and $-1$ in the reference column. The eight target columns form an identity block, so $\bm D$ has row rank eight. Geometry and the reference remain fixed across epochs within a problem. At epoch $j$, $\bm y_j^\rho$ and $\bm y_j^\phi$ are the eight-dimensional code and carrier-phase observations in metres, generated by
\begin{equation}
\bm y_j^\rho=\bm D\bm U\bm x+\bm D\bm\varepsilon_j^\rho,
\qquad
\bm y_j^\phi=\bm D\bm U\bm x+\lambda\bm z+\bm D\bm\varepsilon_j^\phi.
\label{eq:sd}
\end{equation}
The superscripts $\rho$ and $\phi$ denote code and carrier phase. The nine-dimensional $\bm\varepsilon_j^\rho$ and $\bm\varepsilon_j^\phi$ contain the errors before differencing. The entries of $\bm z\in\mathbb Z^8$ are in cycles and the GPS L1 wavelength $\lambda=0.1903\,\mathrm m$ converts cycles to range. The true baseline components are independent $\mathcal N(0,8^2)$ in metres and the true integer components are uniform on $\{-20,\ldots,20\}$. Stacking the code and phase blocks gives
\[
\bm y_j=\begin{bmatrix}\bm y_j^\rho\\\bm y_j^\phi\end{bmatrix},\qquad
\bm A=\begin{bmatrix}\bm D\bm U\\\bm D\bm U\end{bmatrix},\qquad
\bm B=\begin{bmatrix}\bm0\\\lambda\bm I_8\end{bmatrix}.
\]
The identity $\bm I_8$ and the zero block are both of order eight. This simulation consequently has $m=16$, $p=3$ and $n_z=8$.

The nominal undifferenced scales are $\sigma_s^\rho=\sigma_\rho/\max(\sin e_s,0.15)$ and $\sigma_s^\phi=\sigma_\phi/\max(\sin e_s,0.15)$, with $\sigma_\rho=0.05\,\mathrm m$ and $\sigma_\phi=0.015\,\mathrm m$. The laws in~\eqref{eq:noise-laws} are drawn independently at each satellite's corresponding scale before differencing. For observation type $a\in\{\rho,\phi\}$, the nominal differenced covariance block is $\bm D\operatorname{diag}((\sigma_1^a)^2,\ldots,(\sigma_9^a)^2)\bm D^{\mathsf T}$. We apply the inverse lower Cholesky factor of the combined block covariance to the observations and both designs. This is the working transform for all four laws.

For favourable geometry, draw each azimuth uniformly from $[0,2\pi)$ and each elevation uniformly from $[25^\circ,80^\circ]$. Weak geometry uses $[25^\circ,145^\circ]$ and $[12^\circ,35^\circ]$, respectively. In the five-epoch bias experiment, the two lowest-elevation non-reference satellites receive a positive bias for two consecutive epochs each. The code and phase of each selected satellite share an affected interval, and starts are independent across selected satellites. The dimensionless $b$ multiplies the relevant nominal undifferenced scale. The $b=0$ condition adds no local bias. Bias is added before differencing and whitening. Each geometry--law--bias cell contains 1,000 paired problems.

\subsection{Recovery across noise, geometry and duration}\label{sec:noiselaw}
Table~\ref{tab:noiselaw} groups recovery by geometry and background error law. Read across $b$ within a method and down a column for comparisons on paired problems.
\begin{table}[htbp]
\caption{Full-vector recovery (\%) on paired five-epoch problems, by satellite geometry, background error law and two-epoch local bias.}\label{tab:noiselaw}
\centering\footnotesize
\renewcommand{\arraystretch}{1.12}
\setlength{\tabcolsep}{2pt}
\begin{tabular*}{\textwidth}{@{\extracolsep{\fill}}lrrrrrrrrrrrr@{}}\toprule
 & \multicolumn{12}{c}{Background error law} \\
 & \multicolumn{3}{c}{Gaussian} & \multicolumn{3}{c}{Student-$t_3$} & \multicolumn{3}{c}{Mixture} & \multicolumn{3}{c}{Cauchy} \\
\cmidrule(lr){2-4}\cmidrule(lr){5-7}\cmidrule(lr){8-10}\cmidrule(lr){11-13}
Method & $b=0$ & $b=5$ & $b=10$ & $b=0$ & $b=5$ & $b=10$ & $b=0$ & $b=5$ & $b=10$ & $b=0$ & $b=5$ & $b=10$ \\
\midrule
\multicolumn{13}{l}{\textit{Favourable geometry}} \\
WLS & 100.0 & 95.1 & 30.1 & 97.3 & 85.9 & 23.2 & 20.5 & 13.9 & 5.5 & 12.5 & 10.1 & 2.9 \\
Huber IRLS & 99.9 & 97.9 & 95.5 & 99.5 & 90.0 & 82.1 & 75.1 & 50.5 & 26.7 & 80.3 & 54.5 & 33.7 \\
Cauchy IRLS & 99.9 & 97.8 & 95.4 & 99.4 & 89.1 & 82.2 & 75.6 & 49.8 & 27.6 & 81.7 & 56.6 & 35.4 \\
IGG~III IRLS & 99.9 & 98.3 & 95.9 & 99.3 & 88.8 & 85.2 & 78.8 & 50.2 & 34.7 & 82.8 & 55.4 & 41.5 \\
Tukey IRLS & 99.9 & 97.8 & 96.2 & 99.4 & 89.1 & 84.4 & 77.9 & 51.3 & 28.8 & 83.5 & 57.7 & 36.7 \\
GMM, $K=2$ & 99.9 & 97.6 & 94.8 & 99.2 & 87.7 & 79.6 & 73.4 & 45.9 & 27.1 & 76.8 & 49.7 & 33.2 \\
\addlinespace
\multicolumn{13}{l}{\textit{Weak geometry}} \\
WLS & 93.9 & 19.6 & 0.0 & 64.4 & 13.0 & 0.0 & 0.6 & 0.4 & 0.0 & 1.0 & 0.7 & 0.0 \\
Huber IRLS & 89.8 & 55.3 & 37.7 & 79.5 & 36.2 & 26.0 & 19.2 & 5.7 & 2.7 & 29.1 & 8.5 & 3.7 \\
Cauchy IRLS & 88.1 & 54.1 & 42.8 & 77.7 & 36.0 & 28.7 & 20.1 & 5.7 & 3.2 & 32.8 & 9.8 & 5.0 \\
IGG~III IRLS & 85.5 & 59.5 & 52.2 & 76.2 & 40.3 & 39.7 & 24.5 & 6.5 & 8.7 & 35.4 & 11.0 & 10.5 \\
Tukey IRLS & 88.8 & 55.1 & 57.6 & 78.2 & 36.5 & 41.0 & 25.0 & 5.4 & 5.1 & 34.8 & 10.1 & 7.7 \\
GMM, $K=2$ & 83.9 & 52.5 & 37.8 & 73.8 & 34.0 & 25.1 & 19.6 & 4.9 & 3.2 & 28.8 & 9.3 & 4.7 \\
\addlinespace
\bottomrule\end{tabular*}
\par\smallskip\parbox{\textwidth}{\footnotesize Each cell contains 1,000 paired problems under nominal Cholesky whitening. The bias is $b$ times the nominal undifferenced scale of the relevant observation type, applied to two low-elevation non-reference satellites for two consecutive epochs. The $b=0$ cell adds no bias. Mixture denotes the Gaussian-mixture error law in Eq.~\eqref{eq:noise-laws}. Recovery requires full-vector agreement with simulated truth. All planned problems count, with incomplete searches counted as failures.}
\end{table}

For example, with favourable geometry and Gaussian background noise the zero-bias recovery is essentially 100\% for every method. At $b=10$, WLS falls to 30.1\% while the four IRLS configurations remain at 95.4--96.2\%. With favourable geometry and Gaussian-mixture errors without bias, WLS already recovers only 20.5\%, compared with 75.1--78.8\% for the IRLS configurations. These contrasts distinguish an advantage driven by localized bias from one present under a heavy-tailed background alone.

Under weak geometry and Gaussian-mixture noise at $b=10$, WLS recovers 0.0\% and IRLS recovers from 2.7\% to 8.7\%. With weak geometry and Gaussian noise at $b=0$, WLS recovers 93.9\% versus 85.5--89.8\% for IRLS. In weak geometry with Gaussian-mixture noise, IGG~III recovery is 6.5\% at $b=5$ and 8.7\% at $b=10$. These two estimates use separately generated problems. Pairing applies to methods within each cell.

\label{sec:duration}Table~\ref{tab:duration} holds favourable geometry and $b=5$ fixed. Let $d$ be the number of consecutive biased epochs per affected satellite in a record of $N$ epochs. The two lengths use separately generated problems. The four IRLS and three GMM configurations exceed WLS in every tested Gaussian-mixture and Cauchy duration cell. The smallest difference is 4.2 percentage points (pp), attained by each GMM configuration at $N=5$ and $d=4$ under both laws. The Gaussian controls give 70.7\% versus 30.9\% at $N=5,d=3$ and 76.4\% versus 39.5\% at $N=10,d=6$ for WLS and IGG~III, respectively. The Student-$t_3$ control at $N=5,d=3$ gives 59.6\% for WLS and 45.7\% for Huber. Tables~S3--S4 give all generating laws and paired WLS contrasts. Fig.~S1 shows the Gaussian and Student-$t_3$ controls.
\begin{table}[htbp]
\caption{Full-vector recovery (\%) by local-bias duration under Gaussian-mixture and Cauchy background noise. The panels use separately generated five-epoch and ten-epoch problems.}\label{tab:duration}
\centering\small
\renewcommand{\arraystretch}{1.06}
\setlength{\tabcolsep}{2pt}
\begin{tabular*}{\textwidth}{@{\extracolsep{\fill}}lrrrrrrrr@{}}\toprule
Method & \multicolumn{4}{c}{Mixture} & \multicolumn{4}{c}{Cauchy} \\
\midrule
\multicolumn{9}{l}{\textit{$N=5$ epochs}} \\
\shortstack[l]{Consecutive biased\\epochs per satellite, $d$} & 1 & 2 & 3 & 4 & 1 & 2 & 3 & 4 \\
WLS & 19.0 & 15.1 & 10.5 & 4.8 & 12.2 & 10.0 & 6.9 & 3.9 \\
Huber IRLS & 66.5 & 49.4 & 24.2 & 11.4 & 71.0 & 52.6 & 23.6 & 11.2 \\
Cauchy IRLS & 66.3 & 49.1 & 22.7 & 10.3 & 73.6 & 53.5 & 23.9 & 10.0 \\
IGG~III IRLS & 67.2 & 52.1 & 23.3 & 10.7 & 74.8 & 54.9 & 23.8 & 9.5 \\
Tukey IRLS & 67.4 & 51.4 & 24.4 & 10.8 & 75.7 & 55.3 & 25.8 & 9.9 \\
GMM, $K=2$ & 63.4 & 45.5 & 19.2 & 9.0 & 67.3 & 47.8 & 18.3 & 8.1 \\
\addlinespace
\multicolumn{9}{l}{\textit{$N=10$ epochs}} \\
\shortstack[l]{Consecutive biased\\epochs per satellite, $d$} & 2 & 3 & 4 & 6 & 2 & 3 & 4 & 6 \\
WLS & 50.0 & 41.9 & 34.4 & 23.8 & 12.9 & 11.6 & 9.3 & 6.4 \\
Huber IRLS & 93.4 & 87.8 & 76.1 & 41.1 & 96.0 & 91.1 & 80.5 & 43.4 \\
Cauchy IRLS & 93.3 & 87.9 & 76.2 & 38.4 & 95.9 & 91.5 & 80.7 & 41.2 \\
IGG~III IRLS & 93.8 & 88.3 & 77.1 & 38.6 & 96.3 & 91.5 & 80.9 & 41.8 \\
Tukey IRLS & 93.9 & 88.6 & 76.9 & 40.2 & 96.3 & 92.4 & 80.9 & 42.5 \\
GMM, $K=2$ & 93.0 & 87.4 & 75.6 & 35.6 & 89.8 & 84.1 & 74.7 & 34.9 \\
\addlinespace
\bottomrule\end{tabular*}
\par\smallskip\parbox{\textwidth}{\footnotesize Each cell uses 1,000 paired problems with favourable geometry. The two lowest-elevation non-reference satellites each receive a bias of five nominal scales. $N$ is the window length and $d$ counts consecutive biased epochs per affected satellite. Mixture denotes Gaussian-mixture background noise. Recovery compares the complete vector with simulated truth over all planned problems. Tables~S3--S4 give the full outcomes and paired contrasts. Fig.~S1 shows the Gaussian and Student-$t_3$ controls.}
\end{table}

The duration table should first be read within a record-length panel. For five-epoch mixture records, Huber decreases from 66.5\% at one affected epoch to 11.4\% at four, while WLS decreases from 19.0\% to 4.8\%. At ten epochs under Cauchy, Huber falls from 96.0\% at duration two to 43.4\% at duration six. Thus a gain can remain positive while its absolute usefulness deteriorates sharply. The panels use different generated problems. Within a panel, recovery generally declines as the bias affects more of the record.

\label{sec:bound}A separate Gaussian grid assesses a fixed-weight WLS bound in 160 design cells. The grid crosses two geometries, one or two affected satellites, four bias durations, five amplitudes and two elevation classes, with 1,000 trials per cell. For a fixed design and simulated bias, $p_{\min}$ is the smallest probability that the true vector beats both members of a tested symmetric pair of integer alternatives. It upper-bounds the conditional full-vector WLS success probability under nominal Gaussian errors. Across cells, the rank correlation between the cell-mean bound and observed WLS failure is $-0.99$. The derivation is given under ``Fixed-weight Gaussian bound'' in the Supplementary Information. The broader simulation design in Section~\ref{sec:integrity} evaluates acceptance.
\FloatBarrier

\section{Semi-synthetic comparisons on measured observations}\label{sec:realdata}\label{sec:realdata-design}
\enlargethispage{7pt}
We add controlled errors to short windows from measured GNSS observations and evaluate integer recovery and baseline error. Sections~\ref{sec:realdata-injection}--\ref{sec:realdata-heavy} add sparse innovations in the original-order Cholesky-whitened coordinates. Section~\ref{sec:extensions} adds Markov errors to raw double differences or deterministic biases to consecutive epochs. All designs use the same code-and-phase model. Section~\ref{sec:whitening} applies alternative whitening transforms to the same contaminated raw data.

Public observations of the 20.210\,m YARR--YAR3 baseline in Western Australia on 31 March and 30 April 2023 supply measured noise and time-varying geometry. Both the base (YARR) and rover (YAR3) use Septentrio PolaRx5 receivers. GPS C1C code and L1C phase are sampled every 30\,s, giving a grid of 2,880 epochs per station-day. The preprocessing excludes a satellite when its mean elevation over the two receivers is below $10^\circ$. Observations and broadcast ephemerides are available from Geoscience Australia\footnote{\url{https://data.gnss.ga.gov.au/}; source hashes are archived.}. A nonzero broadcast health indicator excludes G22. At this sampling interval, the time from first to last epoch is 120\,s for $N=5$ and 270\,s for $N=10$. The tables use 03-31 and 04-30 for 31 March and 30 April 2023.

The unknown $\bm x=\bm r_{\rm rover}-\bm r_{\rm base}$ is the three-dimensional Earth-centred, Earth-fixed (ECEF) baseline. At each epoch, we first subtract the rover observation from the base observation for the same satellite, then subtract the reference-satellite value from each selected non-reference value. This receiver-and-satellite double difference (DD) gives the code vector $\bm P_j^{DD}$ and the carrier-phase vector $\bm\Phi_j^{DD}$; phase is converted from cycles to metres before fitting. The reference is the satellite with the highest mean of the two receivers' elevations at the first epoch. An eligible window keeps that reference and the satellite set throughout.

To establish the sign of the geometry term, let $\rho_{{\rm base},j}^s$ and $\rho_{{\rm rover},j}^s$ be geometric distances to satellite $s$, and let $\bm u_j^s$ point from the rover to that satellite. Expanding the base distance about the rover position, displaced by $-\bm x$, gives $\rho_{{\rm base},j}^s-\rho_{{\rm rover},j}^s=(\bm u_j^s)^{\mathsf T}\bm x+O(\|\bm x\|^2/\rho_j^s)$, where $\rho_j^s=\rho_{{\rm rover},j}^s$. The positive sign thus corresponds to the base-minus-rover observation convention. Let $\bm U_j$ stack the directions for the reference and selected non-reference satellites, and let $\bm D_j$ contain their target-minus-reference rows. For $n_z$ non-reference satellites, these matrices have sizes $(n_z+1)\times3$ and $n_z\times(n_z+1)$, respectively. The unwhitened model is
\begin{equation}
\begin{bmatrix}\bm P_j^{DD}\\\bm\Phi_j^{DD}\end{bmatrix}
=\begin{bmatrix}\bm D_j\bm U_j\\\bm D_j\bm U_j\end{bmatrix}\bm x
+\begin{bmatrix}\bm0\\\lambda\bm I_{n_z}\end{bmatrix}\bm z+\bm e_j^{\rm raw}.
\label{eq:real-dd}
\end{equation}
Both DD observation vectors are in metres. Each ambiguity is the target-minus-reference integer in cycles under the receiver-difference convention. We use the GPS L1 wavelength $\lambda=0.1903\,\mathrm m$. There are $n_z$ retained non-reference satellites, so $m=2n_z$ and $p=3$. The integer dimension can vary across windows. The stacked left-hand side is $\bm y_j^{\rm raw}$ and the coefficient blocks are $\bm A_j^{\rm raw}$ and $\bm B_j^{\rm raw}$. The fit estimates the baseline directly from the measured DD observations. The exact header-coordinate geometry enters the construction of the internal reference. Linearization remainder and unmodelled differential effects remain in $\bm e_j^{\rm raw}$.

The nominal single-receiver scales are $0.30\,\mathrm m/\max(\sin e_{ij},0.15)$ for code and $0.020\,\mathrm m/\max(\sin e_{ij},0.15)$ for phase, with $e_{ij}$ the mean elevation over the two receivers. For each observation type $a\in\{\rho,\phi\}$, the nominal DD covariance block is $2\bm D_j\operatorname{diag}((\sigma_{ij}^{a})^2)\bm D_j^{\mathsf T}$. The factor two comes from differencing across receivers in the nominal independent equal-scale model. The two blocks form $\bm\Sigma_j=\bm L_j\bm L_j^{\mathsf T}$, where $\bm L_j$ is lower triangular of order $2n_z$. Original-order Cholesky whitening applies $\bm L_j^{-1}$ to the observations and both designs. Code rows precede phase rows and the target-satellite order is retained within each block. Both dates informed the working scales.

For each target--reference pair, subtract the exact double-difference geometric range evaluated at the RINEX header coordinates from the measured double-difference carrier phase. Express the residual in cycles and split its time series at gaps and at discontinuities flagged as suspected cycle slips. On an arc that passes the prespecified stability screens, round the median residual to obtain the pairwise reference integer. The vector $\bm z_{\rm ref}$ collects the reference integers from accepted arcs covering the central epoch of the window, using the later central epoch for even $N$. We require at least four qualifying target satellites. The arc-derived integers are used for scoring. Stable common bias can affect them. Supplementary Section 7.3 and Tables S12--S13 give the screening rules and counts.

Within each date and window length, select eligible windows in chronological order without reusing epochs. The resulting five-epoch counts are 535 on 31 March and 532 on 30 April. The corresponding ten-epoch counts are 247 and 246. The two lengths share selected observations, and the injection conditions reuse the windows.

For each date, $\bm x_{\rm ref}$ is the rover-minus-base ECEF baseline formed from the approximate station coordinates in the Receiver Independent Exchange Format (RINEX) observation headers. For a completed fit in window $i$, $\hat{\bm x}_i=\hat{\bm x}(\check{\bm z}_i)$ is the baseline obtained under that method's final working metric or retained rows. Its baseline error is $e_{{\rm pos},i}=\|\hat{\bm x}_i-\bm x_{\rm ref}\|_2$. Among $n$ planned windows, $n_{10}$ counts completed solutions within 10\,cm of the coordinate benchmark, and $n_{a,10}$ counts those that also pass acceptance. Define $F_{10}=n_{10}/n$ and $V_{10}=n_{a,10}/n$. Table percentages report 100 times these fractions. The position metrics use the coordinate benchmark, whereas recovery uses the arc-derived integers.

\subsection{Sparse innovation design}\label{sec:realdata-injection}
On the measured windows without added errors, the lowest WLS recovery is 528/532 (99.2\%) on 30 April at $N=5$. The lowest across weighting configurations is 240/246 (97.6\%) for each GMM component count on 30 April at $N=10$. All three GMM settings have incomplete searches on the same three ten-epoch windows on that date. These clean-data outcomes are in the electronic recovery summaries listed in the Supplementary Information. The sparse injections probe recovery under specified errors on this high-recovery measured baseline.

In the original-order Cholesky-whitened coordinates, let $\bm c_j$ be a binary mask with independent entries $c_{ij}$ across components and epochs. The added-error vector $\bm\epsilon_j$ is drawn independently of the mask. The symbol $\odot$ denotes componentwise multiplication. The altered whitened observation is $\widetilde{\bm y}_j=\bm y_j+\bm c_j\odot\bm\epsilon_j$. The corresponding raw DD observation is
\begin{equation}
\widetilde{\bm y}_j^{\rm raw}=\bm y_j^{\rm raw}+\bm L_j(\bm c_j\odot\bm\epsilon_j),
\qquad c_{ij}\sim\operatorname{Bernoulli}(0.2).
\label{eq:real-sparse}
\end{equation}
The selection probability 0.2 applies separately to each whitened component. A non-diagonal $\bm L_j$ can spread an added error from one selected component across raw DD rows.

The added-error laws are $0.8\mathcal N(0,1)+0.2\mathcal N(0,10^2)$, Student-$t_3$ and Cauchy, with the latter two scaled to unit Gaussian MAD by~\eqref{eq:noise-scales} with $\sigma=1$. For the Gaussian-mixture injection, the marginal probability of a high-variance draw is $0.2\times0.2=0.04$ per whitened coordinate. Including the mask, the added-error law is $0.8\delta_0+0.16\mathcal N(0,1)+0.04\mathcal N(0,10^2)$, with $\delta_0$ a point mass at zero. Methods share the altered observations and the evaluation reference within each date--length--law cell. Pilot exploration preceded the fixed design.

\subsection{Full-vector recovery and positioning}\label{sec:realdata-heavy}
The sparse-injection study has 4,680 window--law inputs: 1,560 selected windows at the two lengths, each evaluated under three injection laws. The eight weighting configurations complete all 37,440 corresponding fits. IDS--ILS completes a further 4,680 fits on the same inputs. In Table~\ref{tab:real-heavy}, Huber and Cauchy IRLS generally have the largest recovery gains over WLS. Pooling both dates at $N=10$ under Gaussian-mixture injection, Huber and Cauchy recover 475/493 (96.3\%) and 474/493 (96.1\%), a difference of one vector. WLS recovers 48/493 (9.7\%) and IDS--ILS 454/493 (92.1\%). At $N=5$, pooling both dates by window count, IDS--ILS recovers 73.8\% under Gaussian-mixture injection and 69.4\% under Cauchy injection. GMM at $K=2$ ranges from 40.9\% to 51.4\% across the eight displayed cells. Tables~S9--S10 describe short-record implementation sensitivities.
\begin{table}[htbp]
\caption{Full-vector agreement with the long-arc reference integers (\%) after sparse injection in original-order Cholesky-whitened coordinates. Huber and Cauchy IRLS recover over 84\% in all eight cells.}\label{tab:real-heavy}
\centering\small
\renewcommand{\arraystretch}{1.12}
\setlength{\tabcolsep}{2pt}
\begin{tabular*}{\textwidth}{@{\extracolsep{\fill}}lrrrrrrrr@{}}\toprule
 & \multicolumn{4}{c}{Mixture} & \multicolumn{4}{c}{Cauchy} \\
 & \multicolumn{2}{c}{$N=5$} & \multicolumn{2}{c}{$N=10$} & \multicolumn{2}{c}{$N=5$} & \multicolumn{2}{c}{$N=10$} \\
Method & 03-31 & 04-30 & 03-31 & 04-30 & 03-31 & 04-30 & 03-31 & 04-30 \\
Planned windows, $n$ & 535 & 532 & 247 & 246 & 535 & 532 & 247 & 246 \\
\midrule
WLS & 12.1 & 12.0 & 8.5 & 11.0 & 16.8 & 19.9 & 17.8 & 17.5 \\
Huber IRLS & 89.2 & 85.5 & 95.1 & 97.6 & 89.5 & 87.6 & 99.2 & 97.6 \\
Cauchy IRLS & 89.0 & 84.8 & 95.5 & 96.7 & 90.8 & 88.0 & 98.4 & 98.4 \\
IGG~III IRLS & 71.6 & 69.5 & 85.8 & 85.4 & 75.7 & 71.6 & 85.8 & 87.0 \\
Tukey IRLS & 82.2 & 77.1 & 89.9 & 89.4 & 85.8 & 78.8 & 93.9 & 94.3 \\
IDS--ILS & 74.2 & 73.3 & 91.1 & 93.1 & 66.9 & 71.8 & 93.5 & 90.2 \\
GMM, $K=2$ & 43.6 & 41.0 & 40.9 & 42.7 & 49.5 & 45.7 & 51.4 & 50.0 \\
\bottomrule\end{tabular*}
\par\smallskip\parbox{\textwidth}{\footnotesize 03-31 and 04-30 denote 31 March and 30 April 2023. $N$ is epochs per window. Each column divides recovered full vectors by the stated number of planned windows. All searches shown here complete. IDS--ILS is iterative data snooping followed by integer least squares, which retains the full ambiguity dimension. Mixture refers to the Gaussian-mixture amplitude law in Section~\ref{sec:realdata-injection}.}
\end{table}

Under Student-$t_3$ injection, each IRLS configuration improves on WLS in all four date--length cells. Each GMM component count scores below WLS in three cells (Table~S15 and the electronic recovery summaries). At least two of the three GMM settings return different vectors on 13 of the 4,680 inputs (Supplementary Section 7.2).

Paired block resampling uses the same resampled windows for both methods. The four IRLS and three GMM differences from WLS span 25.8 to 88.5\,pp across the Gaussian-mixture and Cauchy cells and the archived 30-, 60- and 120-minute block and boundary settings. The smallest lower endpoint is 19.5\,pp for the 30 April five-epoch Cauchy GMM comparisons with 120-minute blocks and crossing windows retained. The minimum is 20.0\,pp with 60-minute blocks and the same boundary rule (Table~S15). Table~\ref{tab:paired-recovery} shows the four IRLS comparisons at 60 minutes. For example, Huber improves on WLS by 77.0\,pp in the five-epoch mixture windows on 31 March, with an interval of 73.9 to 80.1\,pp. The 95\% percentile intervals use 2,000 resamples and describe contrasts within the selected dates. Dependence can persist across adjacent blocks and windows within the same long arc.

\begin{table}[htbp]
\caption{Paired differences in full-vector recovery from WLS on the sparse-injection windows, in percentage points. Each entry gives the observed difference with its 95\% block-bootstrap percentile interval.}\label{tab:paired-recovery}
\centering\small
\renewcommand{\arraystretch}{1.12}
\setlength{\tabcolsep}{2pt}
\begin{tabular*}{\textwidth}{@{\extracolsep{\fill}}lrrrrr@{}}\toprule
 & & \multicolumn{4}{c}{IRLS weight function} \\
Date & Injection law & Huber & Cauchy & IGG~III & Tukey \\
\midrule
\multicolumn{6}{l}{\textit{$N=5$ epochs}} \\
03-31 & Mixture & \shortstack[r]{$77.0$\\$[73.9,80.1]$} & \shortstack[r]{$76.8$\\$[73.8,79.9]$} & \shortstack[r]{$59.4$\\$[55.5,63.4]$} & \shortstack[r]{$70.1$\\$[66.9,73.5]$} \\
\addlinespace[3pt]
04-30 & Mixture & \shortstack[r]{$73.5$\\$[70.2,76.6]$} & \shortstack[r]{$72.7$\\$[68.6,76.9]$} & \shortstack[r]{$57.5$\\$[51.6,62.7]$} & \shortstack[r]{$65.0$\\$[59.2,70.4]$} \\
\addlinespace[3pt]
03-31 & Cauchy & \shortstack[r]{$72.7$\\$[68.2,76.7]$} & \shortstack[r]{$74.0$\\$[69.4,78.3]$} & \shortstack[r]{$58.9$\\$[54.4,63.1]$} & \shortstack[r]{$69.0$\\$[64.5,73.2]$} \\
\addlinespace[3pt]
04-30 & Cauchy & \shortstack[r]{$67.7$\\$[64.0,71.6]$} & \shortstack[r]{$68.0$\\$[63.8,72.6]$} & \shortstack[r]{$51.7$\\$[47.0,56.1]$} & \shortstack[r]{$58.8$\\$[54.9,62.9]$} \\
\addlinespace[3pt]
\addlinespace
\multicolumn{6}{l}{\textit{$N=10$ epochs}} \\
03-31 & Mixture & \shortstack[r]{$86.6$\\$[82.7,90.6]$} & \shortstack[r]{$87.0$\\$[83.4,90.8]$} & \shortstack[r]{$77.3$\\$[72.7,82.0]$} & \shortstack[r]{$81.4$\\$[76.4,86.1]$} \\
\addlinespace[3pt]
04-30 & Mixture & \shortstack[r]{$86.6$\\$[81.9,91.1]$} & \shortstack[r]{$85.8$\\$[80.7,90.6]$} & \shortstack[r]{$74.4$\\$[68.7,80.4]$} & \shortstack[r]{$78.5$\\$[73.4,83.0]$} \\
\addlinespace[3pt]
03-31 & Cauchy & \shortstack[r]{$81.4$\\$[75.6,86.3]$} & \shortstack[r]{$80.6$\\$[75.2,85.3]$} & \shortstack[r]{$68.0$\\$[59.3,76.3]$} & \shortstack[r]{$76.1$\\$[69.9,82.0]$} \\
\addlinespace[3pt]
04-30 & Cauchy & \shortstack[r]{$80.1$\\$[74.3,85.1]$} & \shortstack[r]{$80.9$\\$[75.3,85.7]$} & \shortstack[r]{$69.5$\\$[63.7,75.3]$} & \shortstack[r]{$76.8$\\$[70.7,82.5]$} \\
\addlinespace[3pt]
\addlinespace
\bottomrule\end{tabular*}
\par\smallskip\parbox{\textwidth}{\footnotesize Each pair uses the same window and injection under original-order Cholesky whitening. Windows are assigned to 60-minute blocks by start time and crossing windows are retained. The 2,000 resamples use common block indices across methods within each date and cell. Differences use unrounded counts and are computed as IRLS minus WLS. Positive values favour IRLS. See Table~\ref{tab:real-heavy} for counts and Table~S15 for the extended 60-minute contrasts.}
\end{table}

Table~\ref{tab:position} reports the median, 95th percentile (P95) and $F_{10}$ for baseline errors before ratio screening on the same inputs. Across the eight Gaussian-mixture and Cauchy date--length--law cells, the IGG~III medians range from 4.2 to 4.6\,cm and the WLS medians from 70.6 to 97.1\,cm (Table~S14). Pooling the individual windows across dates within each length and law gives median ranges of 4.3 to 4.6\,cm and 76.5 to 95.1\,cm, respectively. The electronic positioning summaries retain empirical RMSE as a sample statistic.
\begin{table}[htbp]
\caption{Three-dimensional baseline error before ratio screening, relative to the header-coordinate benchmark. Windows from both dates are pooled at each record length.}\label{tab:position}
\centering\small
\renewcommand{\arraystretch}{1.12}
\setlength{\tabcolsep}{2pt}
\begin{tabular*}{\textwidth}{@{\extracolsep{\fill}}lrrrrrr@{}}\toprule
 & \multicolumn{3}{c}{Mixture} & \multicolumn{3}{c}{Cauchy} \\
Method & Median (cm) & P95 (cm) & $F_{10}$ (\%) & Median (cm) & P95 (cm) & $F_{10}$ (\%) \\
\midrule
\multicolumn{7}{l}{\textit{$N=5$ epochs; 1,067 planned windows}} \\
WLS & 95.1 & 293.8 & 14.9 & 76.5 & 759.7 & 22.2 \\
Huber IRLS & 4.5 & 56.8 & 89.5 & 4.4 & 58.1 & 90.9 \\
Cauchy IRLS & 4.4 & 60.5 & 89.5 & 4.4 & 53.8 & 91.8 \\
IGG~III IRLS & 4.6 & 74.3 & 79.9 & 4.5 & 70.9 & 81.9 \\
Tukey IRLS & 4.5 & 64.7 & 85.8 & 4.4 & 62.7 & 89.2 \\
IDS--ILS & 5.5 & 95.3 & 73.6 & 5.5 & 95.5 & 70.3 \\
GMM, $K=2$ & 13.8 & 103.4 & 48.5 & 7.1 & 96.2 & 54.1 \\
\addlinespace
\multicolumn{7}{l}{\textit{$N=10$ epochs; 493 planned windows}} \\
WLS & 79.8 & 198.0 & 11.0 & 76.9 & 537.3 & 19.9 \\
Huber IRLS & 4.4 & 6.3 & 96.8 & 4.3 & 6.1 & 98.6 \\
Cauchy IRLS & 4.4 & 6.1 & 96.8 & 4.3 & 5.9 & 99.0 \\
IGG~III IRLS & 4.3 & 31.0 & 94.3 & 4.3 & 28.2 & 93.3 \\
Tukey IRLS & 4.3 & 6.9 & 95.9 & 4.3 & 6.0 & 97.2 \\
IDS--ILS & 4.5 & 52.7 & 92.3 & 4.5 & 45.2 & 91.5 \\
GMM, $K=2$ & 21.9 & 90.7 & 46.5 & 7.1 & 83.2 & 57.4 \\
\addlinespace
\bottomrule\end{tabular*}
\par\smallskip\parbox{\textwidth}{\footnotesize P95 is the 95th percentile of individual baseline errors pooled across dates. Each error uses the baseline fitted at the estimated integer vector. $F_{10}$ is the fraction of all planned windows with a completed solution and baseline error of at most 10\,cm. The table reports $100F_{10}$. Table~\ref{tab:sparse-acceptance} gives $V_{10}$ after screening. Table~S14 gives date-specific errors. Mixture denotes Gaussian-mixture injection.}
\end{table}

The median alone hides important differences. In the pooled ten-epoch mixture sample, Huber and IGG~III medians are close (4.4 and 4.3\,cm), but their P95 errors are 6.3 and 31.0\,cm and their $F_{10}$ values are 96.8\% and 94.3\%. IDS--ILS has a 4.5\,cm median but a 52.7\,cm P95. The tail and threshold columns therefore qualify comparisons that would look similar by median. At five epochs the Huber mixture P95 is 56.8\,cm despite a 4.5\,cm median; a typical accurate baseline coexists with a poorer upper tail. Pooling uses the individual baseline errors from both dates. Table~S14 gives the date-specific results.

In the four Gaussian-mixture and Cauchy length groups of the 48-case audit, IGG~III's mean within-case percentage of final working weights at or below $10^{-11}$ ranges from 15.5\% to 19.0\% (Table~S11). Changing the augmentation response from the weighted mean $\bar y_i$ to zero changes the full vector in five of 48 cases for IGG~III and three of 48 for Tukey. The principal results use the mean-preserving augmentation in~\eqref{eq:mean-floor}.

\subsection{Dependent contamination and persistent bias}\label{sec:extensions}
The Markov design adds errors before whitening to raw DD rows of the two retained non-reference satellites with the lowest mean elevations over the window. The code-only condition changes only code rows. The code-and-phase condition changes both corresponding row types. Every fit uses both code and phase observations.

For each selected satellite $i$, let $\xi_{ij}$ be its binary contamination state at epoch $j$. States evolve independently across satellites, with $\Pr(\xi_{i1}=1)=0.2$, $\Pr(\xi_{ij}=1\mid \xi_{i,j-1}=1)=0.8$ and $\Pr(\xi_{ij}=1\mid \xi_{i,j-1}=0)=0.05$. The initial probability 0.2 is stationary for these transitions. A state of 1 activates the innovation at that epoch. For an active observation, multiply the nominal DD standard deviation (the square root of the covariance diagonal entry) by a draw from the unit-scale Gaussian-mixture law in~\eqref{eq:noise-laws}, or by a standard Cauchy draw scaled by $\Phi^{-1}(0.75)\approx0.6745$. In the code-and-phase condition, the two observation types share the satellite state and use independent amplitude draws conditional on that state.

The chains provide the temporal dependence. Methods share the seeds and the added raw DD errors within each cell.

The design has 6,240 window--law--composition inputs and 56,160 method records across nine configurations. Each GMM configuration has four incomplete searches, counted as failures. The other configurations complete every search. The four IRLS and three GMM configurations exceed WLS in every date--length--law--composition cell. Table~\ref{tab:markov} pools dates by planned window count. At $N=10$ with code-only Gaussian-mixture injection, Huber and IDS--ILS both recover 99.2\% of the reference vectors and WLS recovers 87.4\%. With code-and-phase Gaussian-mixture injection, IDS--ILS, Huber, Cauchy and WLS recover 474, 460, 456 and 330 of 493 vectors (96.1\%, 93.3\%, 92.5\% and 66.9\%). Huber ranks above IDS--ILS in the ten-epoch sparse Gaussian-mixture comparison. The designs differ in injection coordinates, affected rows and temporal structure.

For the ten-epoch code-and-phase mixture cell on 30 April, IDS--ILS exceeds WLS by 32.1\,pp with a 60-minute paired-block interval from 26.6 to 38.3\,pp. Huber minus IDS--ILS is $-4.1$\,pp, with interval $[-7.0,-1.6]$\,pp. The corresponding March interval includes zero (Table~S20). Tables~S17--S18 give all date-specific recovery rates.
\begin{table}[htbp]
\caption{Full-vector agreement with the internal reference (\%) under Markov-dependent injection in raw double-difference (DD) observations.}\label{tab:markov}
\centering\small
\renewcommand{\arraystretch}{1.12}
\setlength{\tabcolsep}{2pt}
\begin{tabular*}{\textwidth}{@{\extracolsep{\fill}}lrrrrrrrr@{}}\toprule
 & \multicolumn{8}{c}{Raw DD rows receiving added errors} \\
 & \multicolumn{4}{c}{Code only} & \multicolumn{4}{c}{Code and phase} \\
 & \multicolumn{2}{c}{Mixture} & \multicolumn{2}{c}{Cauchy} & \multicolumn{2}{c}{Mixture} & \multicolumn{2}{c}{Cauchy} \\
Method & $N=5$ & $N=10$ & $N=5$ & $N=10$ & $N=5$ & $N=10$ & $N=5$ & $N=10$ \\
\midrule
WLS & 82.1 & 87.4 & 84.5 & 90.5 & 67.7 & 66.9 & 71.1 & 70.0 \\
Huber IRLS & 97.0 & 99.2 & 97.9 & 99.2 & 86.1 & 93.3 & 89.2 & 95.9 \\
Cauchy IRLS & 97.2 & 99.0 & 97.9 & 99.0 & 85.7 & 92.5 & 88.8 & 95.3 \\
IGG~III IRLS & 96.5 & 99.2 & 96.9 & 99.4 & 82.6 & 86.6 & 84.1 & 89.7 \\
Tukey IRLS & 97.2 & 99.0 & 98.0 & 99.4 & 82.7 & 87.8 & 85.7 & 91.5 \\
IDS--ILS & 94.8 & 99.2 & 95.8 & 99.2 & 87.4 & 96.1 & 87.5 & 96.8 \\
GMM, $K=2$ & 94.0 & 98.4 & 94.7 & 98.6 & 75.7 & 77.3 & 78.1 & 77.3 \\
\bottomrule\end{tabular*}
\par\smallskip\parbox{\textwidth}{\footnotesize All fits use code and phase. For each composition and amplitude law, recovery pools recovered full vectors over both dates, with 1,067 planned windows at $N=5$ and 493 at $N=10$. Recovery requires full-vector agreement with the internal reference. Incomplete searches count as failures. Methods share perturbed inputs within a cell. Mixture denotes Gaussian-mixture amplitudes. Date-specific rates are in Tables~S17--S18 and paired differences are in Table~S20.}
\end{table}

A separate experiment adds a positive 1.0\,m code bias and 0.050\,m phase bias to each of the two lowest-elevation target satellites for $d$ consecutive epochs within the short window. With four of five epochs affected, WLS recovers 67.5\% on 31 March and 68.6\% on 30 April, compared with 53.8\% and 56.2\% for IGG~III. Supplementary Section 7.3 gives the deterministic start rule and the single-observation-type controls.
\FloatBarrier

\section{Accepted-solution performance and applicability}\label{sec:conditions}
Recovery is the fraction of planned problems that produce a correct integer vector. Acceptance determines which first-ranked candidates are reported as fixed solutions. We evaluate frozen thresholds in simulation and across dates, then examine sensitivity to whitening and the computational cost of the estimators.

\subsection{Frozen-threshold simulation evaluation}\label{sec:integrity}
For each error law, the calibration and evaluation sets each contain 26,000 five-epoch problems with disjoint seeds within that law. Each set includes 2,000 problems without added bias, 1,000 per geometry. The other 24,000 problems cross two geometries, two durations, three bias amplitudes and two elevation classes, with 1,000 problems per combination. Bias affects two satellites for one or two epochs. The amplitudes $1\sigma$, $3\sigma$ and $5\sigma$ use the nominal undifferenced scale of the affected observation type. The low- and high-elevation classes select the two lowest and two highest non-reference satellites, respectively. Supplementary Section 2 documents the seed audit. Table~\ref{tab:integrity} reports full-vector recovery before screening on the same evaluation problems used in Table~\ref{tab:frozen}.
\begin{table}[htbp]
\caption{Full-vector recovery (\%) before ratio screening on the simulation evaluation sets used in Table 9. Each error law has 26,000 planned evaluation problems.}\label{tab:integrity}
\centering
\small
\setlength{\tabcolsep}{4pt}
\renewcommand{\arraystretch}{1.12}
\begin{tabular*}{\textwidth}{@{\extracolsep{\fill}}lrrrr@{}}
\toprule
Method & Gaussian & Student-\(t_3\) & Gaussian mixture & Cauchy \\
\midrule
WLS & 87.8 & 73.3 & 10.0 & 7.5 \\
Huber IRLS & 89.8 & 81.8 & 40.8 & 46.5 \\
Cauchy IRLS & 88.7 & 80.8 & 41.5 & 48.8 \\
IGG~III IRLS & 87.9 & 80.3 & 44.3 & 50.7 \\
Tukey IRLS & 89.1 & 81.1 & 43.5 & 50.3 \\
GMM, \(K=2\) & 86.5 & 78.7 & 39.4 & 44.4 \\
\bottomrule
\end{tabular*}
\par\smallskip\parbox{\textwidth}{\footnotesize Each evaluation set combines 2,000 problems without added bias and 24,000 locally biased problems across two geometries. Recovery requires agreement in every integer component with simulated truth. All planned problems count. Incomplete searches count as failures.}
\end{table}

For each method and nominal calibration target $a_{\rm target}\in\{0.05,0.10,0.25,0.50\}$, let $f$ be the fraction of planned calibration problems with valid ratios. Let $Q_R$ be the linearly interpolated sample quantile of the valid ratios. Set
\[
\tau=Q_R\bigl(1-\min\{1,a_{\rm target}/f\}\bigr)
\quad\text{for }f>0,
\]
and set $\tau=+\infty$ if $f=0$. Each method--target pair is calibrated separately. The fitted threshold is applied unchanged in evaluation. Ties and unattainable targets are recorded. Table~\ref{tab:frozen} reports actual evaluation $A$ and $q$ together.
\begin{table}[htbp]
\caption{Actual acceptance and conditional wrong-fix frequency on separate simulation evaluation problems, by nominal calibration target. For both laws, every alternative shown has lower $q$ than WLS at each target. All displayed $q$ values exceed 0.5\%.}\label{tab:frozen}
\centering\small
\renewcommand{\arraystretch}{1.12}
\setlength{\tabcolsep}{2pt}
\begin{tabular*}{\textwidth}{@{\extracolsep{\fill}}lrrrrrrrr@{}}\toprule
 & \multicolumn{8}{c}{Nominal calibration target} \\
Method & \multicolumn{2}{c}{5\%} & \multicolumn{2}{c}{10\%} & \multicolumn{2}{c}{25\%} & \multicolumn{2}{c}{50\%} \\
 & $A$ (\%) & $q$ (\%) & $A$ (\%) & $q$ (\%) & $A$ (\%) & $q$ (\%) & $A$ (\%) & $q$ (\%) \\
\midrule
\multicolumn{9}{l}{\textit{Mixture}} \\
WLS & 4.98 & 69.60 & 10.00 & 74.82 & 25.18 & 81.75 & 50.24 & 86.19 \\
Huber IRLS & 5.14 & 7.41 & 10.31 & 11.87 & 25.08 & 23.14 & 49.93 & 39.44 \\
Cauchy IRLS & 5.17 & 5.21 & 10.25 & 9.12 & 25.24 & 18.79 & 50.03 & 35.64 \\
IGG~III IRLS & 4.83 & 4.30 & 9.96 & 6.64 & 24.77 & 13.93 & 49.88 & 30.32 \\
Tukey IRLS & 5.25 & 4.47 & 10.10 & 7.01 & 25.10 & 16.15 & 49.69 & 32.63 \\
GMM, $K=2$ & 4.82 & 11.64 & 9.57 & 19.45 & 24.67 & 32.80 & 49.75 & 49.68 \\
\addlinespace
\multicolumn{9}{l}{\textit{Cauchy}} \\
WLS & 5.15 & 67.61 & 10.20 & 74.35 & 25.15 & 83.00 & 50.11 & 88.25 \\
Huber IRLS & 5.01 & 2.46 & 9.95 & 5.33 & 25.13 & 13.97 & 50.58 & 30.54 \\
Cauchy IRLS & 4.85 & 2.54 & 9.95 & 3.63 & 25.13 & 9.23 & 50.60 & 24.15 \\
IGG~III IRLS & 4.82 & 2.64 & 9.75 & 3.67 & 24.72 & 8.09 & 49.92 & 21.24 \\
Tukey IRLS & 4.68 & 2.30 & 9.75 & 3.08 & 24.92 & 8.26 & 50.30 & 22.24 \\
GMM, $K=2$ & 4.94 & 10.04 & 9.95 & 16.35 & 25.26 & 29.77 & 50.32 & 45.07 \\
\addlinespace
\bottomrule\end{tabular*}
\par\smallskip\parbox{\textwidth}{\footnotesize Each law has 26,000 calibration and 26,000 evaluation problems, with disjoint seeds within that law. Thresholds are fitted separately by method and target and applied unchanged in evaluation. $A=n_a/n$ uses all planned evaluation problems and $q=n_w/n_a$ uses accepted problems. A fix is wrong if at least one ambiguity differs from simulated truth. A dash in $q$ would denote zero acceptance and an undefined conditional frequency. Mixture denotes Gaussian-mixture errors.}
\end{table}

At the nominal 10\% target, WLS has 1,946 wrong vectors among 2,601 accepted Gaussian-mixture problems and 1,971 among 2,651 accepted Cauchy problems. IGG~III has 172 among 2,589 and 93 among 2,536, respectively. Under Gaussian-mixture noise, $A$ and $q$ are 10.00\% and 74.82\% for WLS, and 9.96\% and 6.64\% for IGG~III. IGG~III has the lowest $q$ among the four IRLS rules under that law. Under Cauchy noise, Tukey has the lowest IRLS $q$ at 3.08\% with $A=9.75\%$. We use 0.5\% as a study reference for the conditional wrong-fix frequency. Within the Gaussian-mixture bias strata, the same IGG~III threshold accepts 1,036 of 8,000 problems at $1\sigma$ with 47 wrong, and 595 of 8,000 at $5\sigma$ with 68 wrong (Table~S8).

Let $S_{\rm rec}$ be the all-planned recovery proportion reported in Table~\ref{tab:integrity}. At most $nS_{\rm rec}$ of the $nA$ accepted vectors can be correct. Thus $q\ge\max(0,1-S_{\rm rec}/A)$ for $A>0$. At the nominal 50\% Gaussian-mixture target, IGG~III has 11,516 correct vectors and 12,969 acceptances among 26,000 problems. The lower bound is 11.20\% and the observed $q$ is 30.32\%. The bound depends on the number of correct candidates. The observed conditional frequency also depends on their ranking by $R$.

The second policy evaluates thresholds at the 0th through 99th percentiles of the valid calibration ratios. Among thresholds with at least 100 calibration acceptances and a two-sided 95\% Wilson upper endpoint for $q$ at or below 0.5\%, it selects the one with the most calibration acceptances. The selected threshold is frozen for evaluation. An empty qualifying set produces no acceptance. No calibration threshold qualifies under Gaussian-mixture or Cauchy noise for any of the eight weighting configurations. Table S28 shows Gaussian and Student-$t_3$ evaluation results, and Table S5 gives the complete outcomes. At zero observed wrong fixes, the nominal upper endpoint is $z_{0.975}^2/(n_a+z_{0.975}^2)$ under independent Bernoulli sampling. Here $z_{0.975}$ is the 0.975 standard normal quantile. The 0.5\% cutoff then requires at least 765 acceptances. The pointwise Wilson calculation is used for threshold selection. The frozen threshold is assessed on the separate evaluation set. Fixed-failure-rate integer-aperture designs target the unconditional wrong-fix frequency $u$ \citep{teunissenverhagen2009,verhagen2013}.

\subsection{Cross-date acceptance after sparse injection}\label{sec:realdata-ratio}
For each specified injection law, window length and configuration, calibrate a threshold on one date using the quantile rule in Section~\ref{sec:integrity} and apply it unchanged on the other date. Repeat with the dates reversed. Both dates informed the nominal working scales. Table~\ref{tab:sparse-acceptance} reports ten-epoch outcomes on 30 April using thresholds calibrated on 31 March. Table~S16 gives the reverse direction for Gaussian-mixture and Cauchy injection at both lengths. The electronic table \texttt{sparse\_acceptance.csv} contains all 432 policies: three laws, two lengths, two directions, four targets and nine configurations.
\begin{table}[htbp]
\caption{Cross-date acceptance on ten-epoch sparse-injection windows. Thresholds calibrated on 31 March 2023 are evaluated on 30 April 2023.}\label{tab:sparse-acceptance}
\centering\small
\renewcommand{\arraystretch}{1.12}
\setlength{\tabcolsep}{2pt}
\begin{tabular*}{\textwidth}{@{\extracolsep{\fill}}lrrrrrrrr@{}}\toprule
 & \multicolumn{8}{c}{Nominal calibration target} \\
Method & \multicolumn{4}{c}{10\%} & \multicolumn{4}{c}{50\%} \\
 & $A$ (\%) & $V_{10}$ (\%) & \shortstack{Wrong/\\accepted} & $q$ (\%) & $A$ (\%) & $V_{10}$ (\%) & \shortstack{Wrong/\\accepted} & $q$ (\%) \\
\midrule
\multicolumn{9}{l}{\textit{Mixture}} \\
WLS & 5.28 & 2.44 & 7/13 & 53.85 & 46.75 & 7.32 & 97/115 & 84.35 \\
Huber IRLS & 11.38 & 11.38 & 0/28 & 0.00 & 51.63 & 51.63 & 0/127 & 0.00 \\
Cauchy IRLS & 11.38 & 11.38 & 0/28 & 0.00 & 49.19 & 49.19 & 0/121 & 0.00 \\
IGG~III IRLS & 8.54 & 8.54 & 0/21 & 0.00 & 53.66 & 53.66 & 2/132 & 1.52 \\
Tukey IRLS & 8.94 & 8.94 & 0/22 & 0.00 & 49.59 & 49.59 & 1/122 & 0.82 \\
IDS--ILS & 14.23 & 14.23 & 0/35 & 0.00 & 58.94 & 58.54 & 1/145 & 0.69 \\
GMM, $K=2$ & 7.32 & 4.88 & 8/18 & 44.44 & 46.34 & 23.17 & 68/114 & 59.65 \\
\addlinespace
\multicolumn{9}{l}{\textit{Cauchy}} \\
WLS & 6.10 & 3.25 & 7/15 & 46.67 & 47.97 & 13.82 & 89/118 & 75.42 \\
Huber IRLS & 8.54 & 8.54 & 0/21 & 0.00 & 50.81 & 50.81 & 0/125 & 0.00 \\
Cauchy IRLS & 9.35 & 9.35 & 0/23 & 0.00 & 52.44 & 52.44 & 0/129 & 0.00 \\
IGG~III IRLS & 10.57 & 10.57 & 0/26 & 0.00 & 52.03 & 52.03 & 0/128 & 0.00 \\
Tukey IRLS & 10.98 & 10.98 & 0/27 & 0.00 & 54.88 & 54.88 & 0/135 & 0.00 \\
IDS--ILS & 8.13 & 8.13 & 0/20 & 0.00 & 45.93 & 45.53 & 1/113 & 0.88 \\
GMM, $K=2$ & 8.54 & 6.10 & 7/21 & 33.33 & 48.78 & 28.46 & 57/120 & 47.50 \\
\addlinespace
\bottomrule\end{tabular*}
\par\smallskip\parbox{\textwidth}{\footnotesize Each row has $n=246$ planned evaluation windows under original-order Cholesky whitening. $A=n_a/n$ and $V_{10}=n_{a,10}/n$ use all planned windows. $q=n_w/n_a$ uses the accepted windows. Here $n_a$ counts acceptances, $n_{a,10}$ counts accepted baselines within 10\,cm of the header-coordinate benchmark, and $n_w$ counts accepted vectors that differ from the internal reference. The rates are displayed in percent. Mixture denotes Gaussian-mixture injection. The electronic table \texttt{sparse\_acceptance.csv} gives accepted baseline-error quantiles.}
\end{table}

At the nominal 10\% target, WLS accepts 13 of the 246 planned Gaussian-mixture windows. Six accepted baselines lie within 10\,cm of the coordinate benchmark, and seven of the accepted vectors differ from the internal reference. Thus $A=13/246$, $V_{10}=6/246$ and $q=7/13$, displayed as 5.28\%, 2.44\% and 53.85\%.

The eight cells cross two calibration-to-evaluation directions, two lengths and two injection laws under original-order Cholesky whitening. At the nominal 10\% target, WLS $q$ ranges from 35.6\% to 68.6\%. Huber and Cauchy have zero observed wrong fixes in 16 method--cell records, with 20 to 67 acceptances each. The nominal Wilson upper endpoints under independent Bernoulli sampling range from 5.42\% to 16.11\% (Table~S26). The windows can be temporally dependent. At nominal 50\%, $q$ spans 0--0.36\% for Huber, 0--0.42\% for Cauchy, 0--14.03\% for IDS--ILS and 67.4--86.9\% for WLS.

At the nominal 50\% Gaussian-mixture target in Table~\ref{tab:sparse-acceptance}, Huber, Cauchy and IDS--ILS have 0/127, 0/121 and 1/145 wrong vectors among accepted windows, respectively. Their P95 baseline errors among accepted windows are 5.4, 5.6 and 7.4\,cm, respectively (\texttt{sparse\_acceptance.csv}, column \texttt{accepted\_p95\_error\_m}, converted to centimetres). IDS--ILS attains the largest displayed $V_{10}$ at 58.54\%, compared with 51.63\% for Huber and 49.19\% for Cauchy. WLS has 97 wrong vectors in 115 acceptances and $V_{10}=7.32\%$. IGG~III has 132 accepted baselines within 10\,cm, including two vectors that differ from their internal reference.

The simulation and measured-observation designs differ in background errors, geometry, injection location and duration (Table~S7). The simulation combines no-added-bias controls with one- or two-epoch satellite biases. In sparse Gaussian-mixture injection, the high-variance component has marginal probability 0.04 per whitened coordinate. Huber's pre-screening recovery ranges from 85.5\% to 97.6\% across the four sparse Gaussian-mixture cells, compared with 40.8\% in the pooled simulation evaluation (Tables~\ref{tab:real-heavy} and~\ref{tab:integrity}).

The archived Markov acceptance analysis contains 576 policies. At the nominal 50\% target with code-and-phase contamination, $q$ ranges from 3.89\% to 11.30\% for WLS, 0\% to 2.26\% for IDS--ILS, 0\% to 1.77\% for Huber and 0\% to 2.44\% for Cauchy (Supplementary Section 5, Table~S19).

For sparse injection, the threshold-learning sensitivity resamples calibration and evaluation blocks separately and applies the cutoff fitted in each calibration draw to the corresponding evaluation draw (Supplementary Section 4). Percentile intervals can remain zero when very few resamples contain a wrong accepted fix. Supplementary Section 6 compares nominal Wilson endpoints for the observed zero-wrong-fix cells under original-order whitening (Table~S26).

\subsection{Whitening as a condition on the gains}\label{sec:representation}\label{sec:whitening}
We hold the contaminated raw DD observations and internal reference vectors fixed. We supply alternative whitening representations to each complete estimator and recalibrate thresholds where specified. The representation can affect component-wise weighting, component medians used for initialization and sequential residual screening.

Write the original nominal covariance as $\bm\Sigma_j=\bm L_j\bm L_j^{\mathsf T}$, with lower Cholesky factor $\bm L_j$. Let $\bm P$ be the $m\times m$ orthogonal permutation matrix reversing the satellite-row order separately within the code and phase blocks, and let $\bm L_{P,j}$ be the lower Cholesky factor of $\bm P\bm\Sigma_j\bm P^{\mathsf T}$. For symmetric whitening, $\bm\Sigma_j^{-1/2}$ denotes the symmetric positive-definite inverse square root. The transformations from the original whitened coordinates to reversed-order Cholesky and symmetric whitening are, respectively,
\begin{equation}
\bm T_j^{\rm rev}=\bm L_{P,j}^{-1}\bm P\bm L_j,
\qquad
\bm T_j^{\rm sym}=\bm\Sigma_j^{-1/2}\bm L_j.
\label{eq:whitening}
\end{equation}
For either transform, the alternative fit uses $(\bm T_j\widetilde{\bm y}_j,\bm T_j\bm A_j,\bm T_j\bm B_j)$. Symmetric whitening preserves the original row order. Both transforms are $m\times m$ and orthogonal. Thus $\|\bm T_j\bm r_j\|_2^2=\|\bm r_j\|_2^2$ and the WLS objective is invariant. The observed WLS vectors agree across the tested representations.

Table~\ref{tab:representation-recovery} summarizes 3,120 window--law inputs in eight cells formed by two dates, two lengths and two laws. All 96 IRLS--WLS contrasts (four rules, three representations and eight cells) are positive. The smallest is 40.98\,pp: on 30 April at $N=5$, with Gaussian-mixture injection and reversed-order Cholesky whitening, Huber recovers 282/532 and WLS 64/532. Of the 64 altered-minus-original contrasts (four rules, two alternatives and eight cells), 60 are negative. Huber's reversed-order difference is $-19.80$\,pp with equal cell weights and $-22.53$\,pp pooled over the 3,120 window--law inputs.

For 31 March at $N=5$ under Gaussian-mixture injection, the Huber reversed-minus-original contrast is $-29.72$\,pp, with a 60-minute paired-block interval of $[-33.71,-25.82]$\,pp. On 30 April at the same length and law, Huber minus Tukey is $+8.46$, $-3.76$ and $-3.95$\,pp in the original, reversed and symmetric representations, respectively. The electronic summaries and interval table provide these individual cell results (Supplementary Section 6).
\begin{table}[htbp]
\caption{Full-vector agreement with the internal reference (\%) under three whitening representations of the same perturbed raw observations.}\label{tab:representation-recovery}
\centering\small
\renewcommand{\arraystretch}{1.12}
\setlength{\tabcolsep}{2pt}
\begin{tabular*}{\textwidth}{@{\extracolsep{\fill}}lrrrrr@{}}\toprule
Method & \shortstack{Original-order\\Cholesky} & \multicolumn{2}{c}{Reversed-order Cholesky} & \multicolumn{2}{c}{Symmetric whitening} \\
 & \shortstack{Recovery\\range (\%)} & \shortstack{Recovery\\range (\%)} & \shortstack{Mean change\\vs original (pp)} & \shortstack{Recovery\\range (\%)} & \shortstack{Mean change\\vs original (pp)} \\
\midrule
WLS & 8.5--19.9 & 8.5--19.9 & 0.0 & 8.5--19.9 & 0.0 \\
Huber IRLS & 85.5--99.2 & 53.0--90.7 & -19.8 & 67.5--95.1 & -10.4 \\
Cauchy IRLS & 84.8--98.4 & 55.1--91.9 & -17.7 & 70.7--96.4 & -8.3 \\
IGG~III IRLS & 69.5--87.0 & 53.4--86.2 & -8.3 & 61.8--88.3 & -5.1 \\
Tukey IRLS & 77.1--94.3 & 56.8--90.7 & -11.1 & 71.4--95.1 & -3.5 \\
IDS--ILS & 66.9--93.5 & 35.5--75.6 & -25.7 & 51.3--85.0 & -14.3 \\
\bottomrule\end{tabular*}
\par\smallskip\parbox{\textwidth}{\footnotesize The range spans eight cells formed by two dates, two lengths and two injection laws. Compute each recovery rate over all planned windows in its cell. The mean change gives equal weight to the eight differences between each alternative and the original-order rate, in pp. All representations use the same 3,120 window--law inputs. Reversal acts separately within the code and phase blocks.}
\end{table}

An alternative transform maps the sparse perturbation in~\eqref{eq:real-sparse} to $\bm T_j(\bm c_j\odot\bm\epsilon_j)$ and can spread the same added error across whitened components. In a separate Markov code-and-phase subset, select 64 evenly spaced windows in each date--length group and apply both injection laws, giving 512 window--law inputs. Huber changes by $-0.20$ and $-1.56$\,pp under reversed and symmetric whitening, respectively. Across methods and cells the changes span $-7.81$ to $+6.25$\,pp. These results compare different injection designs and samples.

With a separate threshold per representation at the nominal 10\% target, Huber and Cauchy have zero observed wrong fixes in all 48 records (two methods, three representations and eight cells), with 17 to 67 acceptances per record. At the nominal 50\% target, Table~\ref{tab:representation-example} shows the cell with the largest observed Huber $q$ spread across representations. The original and reversed representations have $A$ of 51.78\% and 51.40\%, $q$ of 0.00\% and 10.55\%, and $V_{10}$ of 51.78\% and 46.17\%. Table~S23 summarizes ranges for the nominal 10\% target. Table~S27 gives the five-epoch Gaussian-mixture results for both date directions and all four nominal targets, with thresholds calibrated separately for each whitening representation.

A separate comparison transfers the original-order threshold to the reversed-order representation. For five-epoch Cauchy windows calibrated on 30 April and evaluated on 31 March at the nominal 10\% target, reversed-order Huber accepts 49 of 535 windows with its own threshold and 2 of 535 with the original-order threshold. Both accepted subsets have zero observed wrong fixes and all accepted baselines are within 10\,cm. Their $V_{10}$ values are 9.16\% and 0.37\%, respectively.
\begin{table}[htbp]
\caption{Illustrative Huber acceptance after calibrating a separate threshold for each whitening representation. The five-epoch Gaussian-mixture windows use 30 April 2023 for calibration and 31 March 2023 for evaluation, at a nominal 50\% calibration target.}\label{tab:representation-example}
\centering\small
\renewcommand{\arraystretch}{1.12}
\setlength{\tabcolsep}{2pt}
\begin{tabular*}{\textwidth}{@{\extracolsep{\fill}}lrrrrr@{}}\toprule
Representation & $A$ (\%) & $V_{10}$ (\%) & \shortstack{Accepted within\\10 cm (count)} & \shortstack{Wrong/\\accepted} & $q$ (\%) \\
\midrule
Original & 51.78 & 51.78 & 277 & 0/277 & 0.00 \\
Reversed & 51.40 & 46.17 & 247 & 29/275 & 10.55 \\
Symmetric & 51.78 & 50.65 & 271 & 7/277 & 2.53 \\
\bottomrule\end{tabular*}
\par\smallskip\parbox{\textwidth}{\footnotesize Each row uses $n=535$ planned evaluation windows. $A=n_a/n$, $V_{10}=n_{a,10}/n$ and $q=n_w/n_a$ are reported in percent. The count $n_{a,10}$ requires acceptance and baseline error of at most 10\,cm relative to the header-coordinate benchmark. Wrong/accepted is $n_w/n_a$ relative to the internal reference. Original and Reversed refer to original-order and within-code/phase-block reversed-order Cholesky. The electronic whitening acceptance results give the eight Huber cells at this target. This cell shows the largest observed range of $q$ across representations.}
\end{table}

\subsection{Implementation cost and scope}\label{sec:cost}
The timing benchmark uses four evenly spaced windows in each of the four date--length groups, each with no added error, sparse Gaussian-mixture injection and sparse Cauchy injection. This gives 48 inputs. Each of the nine configurations receives two warm-ups and three measured repetitions with a fixed seeded order. The benchmark uses one process and one BLAS thread on an AMD Ryzen 9 8945HX, with Python 3.13.9, NumPy 2.4.1 and SciPy 1.17.0. Timings cover the solver and final-candidate diagnostics, with loading, injection and serialization outside the measured interval. All 1,296 measured invocations complete. Huber has a median paired runtime ratio of 4.04 to WLS at both lengths. The IDS--ILS ratios are 1.18 and 1.34, and the GMM $K=2$ ratios are 3.32 and 14.51, for $N=5$ and $N=10$, respectively (Table~\ref{tab:runtime}). Table~S21 gives the more detailed timing groups.
\begin{table}[htbp]
\caption{Single-process runtime on the selected timing cases, by record length.}\label{tab:runtime}
\centering\small
\renewcommand{\arraystretch}{1.12}
\setlength{\tabcolsep}{2pt}
\begin{tabular*}{\textwidth}{@{\extracolsep{\fill}}lrrrrrr@{}}\toprule
Method & \multicolumn{3}{c}{$N=5$} & \multicolumn{3}{c}{$N=10$} \\
 & Median (ms) & P95 (ms) & \shortstack{Median\\paired ratio\\to WLS} & Median (ms) & P95 (ms) & \shortstack{Median\\paired ratio\\to WLS} \\
\midrule
WLS & 3.45 & 5.66 & 1.00 & 3.44 & 5.33 & 1.00 \\
Huber IRLS & 14.57 & 25.77 & 4.04 & 15.09 & 20.44 & 4.04 \\
Cauchy IRLS & 15.63 & 24.45 & 4.34 & 16.00 & 23.66 & 4.54 \\
IGG~III IRLS & 13.19 & 26.15 & 3.55 & 13.62 & 21.18 & 4.03 \\
Tukey IRLS & 16.07 & 21.49 & 4.45 & 15.93 & 27.89 & 4.46 \\
IDS--ILS & 4.05 & 5.79 & 1.18 & 4.70 & 6.51 & 1.34 \\
GMM, $K=2$ & 10.62 & 50.27 & 3.32 & 49.40 & 74.90 & 14.51 \\
\bottomrule\end{tabular*}
\par\smallskip\parbox{\textwidth}{\footnotesize The paired ratio divides each measured runtime by the WLS runtime for the same input and repetition, then takes the median of those ratios. The time columns give the median and P95 in milliseconds. Each method and length has 24 inputs and three measured repetitions per input after warm-up. All 72 calls complete. Timings include the solver and final-candidate diagnostics and use one process and one BLAS thread. Table~S21 gives the detailed timing groups.}
\end{table}

The 48 selected probes include 16 inputs each without added error, with Gaussian-mixture injection and with Cauchy injection. The median-only estimate is the initialization fit in Section~\ref{sec:irls}. It recovers 7/16 Gaussian-mixture vectors and 5/16 Cauchy vectors, compared with 16/16 and 15/16 after the Huber updates. At the 20-update budget, changing the start to WLS changes the full vector in 0, 1, 10 and 5 of the 48 cases for Huber, Cauchy, IGG~III and Tukey, respectively. This is 16 of 192 paired initialization comparisons. Increasing the budget from 20 to 60 changes zero vectors in this probe set (Table~S24).

The label audit selects 64 windows per date and length, giving 256 windows. For each target--reference pair, use its accepted long arc covering the central epoch. Exclude the evaluation window and five minutes on either side. The two-sided rule requires both remaining sides to pass the stability checks and agree. The rule allowing one side also accepts a single passing side when the other is unavailable. Complete vectors are available for 52/256 windows under the two-sided rule and 245/256 under the one-side-allowed rule. Each reconstructed complete vector agrees with its original reference in the corresponding available cohort (Table~S25). Both rules use the same header-coordinate benchmark.

The measured-observation evidence comes from one 20.210\,m baseline on two dates. Both dates informed the working scales, and several injection laws and window lengths reuse the observations. The integer criterion uses internal long-arc references, and the position criterion uses a RINEX header-coordinate benchmark. The simulated laws and Markov parameters specify controlled errors. The searches and iterative fits use finite budgets, with incomplete searches counted as recovery failures. The measured cross-date comparisons describe threshold transfer under this shared working model. Small accepted samples and dependent windows limit the interpretation of observed wrong-fix frequencies and their nominal Wilson intervals. The results support performance comparisons within the tested designs. Operational risk would require additional evidence. Tables S9--S11 and S21--S25 and the electronic block-resampling tables document the implementation and sampling sensitivities.
\FloatBarrier

\section{Conclusions}\label{sec:conclusion}
With Gaussian-mixture or Cauchy errors, all four IRLS configurations improve full-vector recovery over nominal WLS across the tested simulation duration cells. Under Gaussian or Student-$t_3$ background noise, WLS recovers more vectors in some longer local-bias settings. With sparse Gaussian-mixture injection in the original-order Cholesky-whitened coordinates of 493 ten-epoch measured-observation windows, Huber and Cauchy IRLS recover 475 and 474 internal reference vectors (96.3\% and 96.1\%), compared with 48 for WLS (9.7\%) and 454 for IDS--ILS (92.1\%). In those windows Huber and WLS produce median baseline errors of 4.4 and 79.8\,cm and 95th-percentile errors of 6.3 and 198.0\,cm, respectively.

The estimator ordering depends on the injection and the evaluation endpoint. IDS--ILS recovers more vectors than Huber in the ten-epoch code-and-phase Markov Gaussian-mixture comparison. It also gives the highest observed $V_{10}$ among the methods in the sparse Gaussian-mixture example at the nominal 50\% calibration target. Changing the whitening representation alters IRLS recovery and acceptance performance. IGG~III and Tukey have hard cutoffs, and their integer estimates change in some initialization and augmentation probes. The short-record GMM results depend on the specified ridge and initialization.

The controlled designs, common search and matched inputs support comparisons of recovery, baseline error and acceptance for the tested cases. The measured evidence uses one short baseline, two dates, shared working scales and internal reference integers. Where Huber and Cauchy IRLS show zero wrong fixes in the small cross-date accepted samples at the nominal 10\% target, the nominal Wilson upper endpoints remain above the 0.5\% conditional-frequency reference.

\backmatter
\section*{Statements and Declarations}
\noindent\textbf{Funding.} This research is supported by the National Key R\&D Program of China (2024YFA1012401), the Science and Technology Commission of Shanghai Municipality of China (23JC1400501), and the National Natural Science Foundation of China (12241103).

\noindent\textbf{Competing interests.} The authors declare that they have no competing interests.

\noindent\textbf{Ethics approval and consent to participate.} Not applicable.

\noindent\textbf{Consent for publication.} Not applicable.

\noindent\textbf{Data availability.} The derived data, trial outcomes, configurations and numerical summaries supporting this study are recorded in the research archive prepared for journal peer review (\texttt{Review\_Data\_and\_Code.zip}). The file \texttt{REVIEW\_DATA\_INDEX.md} lists exact electronic-table paths relative to the archive's data directory. The file \texttt{REVIEW\_DATA\_MANIFEST.csv} lists file sizes and SHA-256 hashes. The original GNSS observations are publicly available from Geoscience Australia (\url{https://data.gnss.ga.gov.au/}) and NASA CDDIS (\url{https://cddis.nasa.gov/archive/gnss/data/daily/}); input file names, source URLs, hashes and derived internal labels are included in the review archive. Public deposition of the final research archive is planned.

\noindent\textbf{Code availability.} The custom experiment and analysis code, numerical environment information and source hashes are included in the same research archive, indexed by \texttt{REVIEW\_DATA\_INDEX.md}. Public deposition of the final code is planned.

\noindent\textbf{Declaration of generative AI and AI-assisted technologies.} The authors used AI tools for code checking and collaborative writing.

\noindent\textbf{Author Contributions.} R.~Zhang carried out the main research, implemented the code and numerical analyses, and drafted the manuscript. Y.~Zhang contributed to the study conception and design, discussed the results and supervised the work.

\bibliography{references}
\end{document}


\section*{Supplementary Information}
\textbf{Robust weighting for GNSS integer ambiguity resolution under heavy-tailed errors and localized biases}\\
Ruifeng Zhang$^{1}$ and Yunxin Zhang$^{1,2}$\\
$^{1}$School of Mathematical Sciences, Fudan University, Shanghai 200433, China.\\
$^{2}$Shanghai Key Laboratory for Contemporary Applied Mathematics, Laboratory of Mathematics for Nonlinear Science, Fudan University, Shanghai 200433, China.\\
Corresponding author: Yunxin Zhang, \texttt{xyz@fudan.edu.cn}.

\section{Principal protocol and complete simulation outcomes}
The principal IRLS protocol uses $\hat\sigma^{(t)}=\max(1.4826\,\mathrm{MAD}^{(t)},10^{-6})$, weights bounded below by $10^{-12}$, median initialization, a 20-update budget, real-parameter tolerance $10^{-6}$ and mean-preserving augmentation. The 20,000-node sphere-search budget applies after MLAMBDA-related candidate generation. Incomplete searches count as recovery failures and are ineligible for acceptance. Electronic records are stored in the archive \texttt{Review\_Data\_and\_Code.zip}. Its \texttt{REVIEW\_DATA\_INDEX.md} gives archive-relative paths and \texttt{REVIEW\_DATA\_MANIFEST.csv} gives file hashes.

\paragraph{Notation and displayed rates.}
Full-vector recovery requires agreement with the applicable simulated truth or internal reference in every integer component. A planned problem contributes to the denominator even when its search is incomplete. The common-completed analysis retains problems on which every compared configuration completes. A rescue is a correct alternative vector on a problem WLS gets wrong. A loss is a correct WLS vector on a problem the alternative gets wrong. The two rates use the same planned-pair denominator. Let $n$ be the number of planned problems, $n_a$ the accepted count and $n_w$ the accepted wrong-vector count. For positioning, $n_{10}$ counts completed solutions within 10\,cm of the header-coordinate benchmark and $n_{a,10}$ counts those that are also accepted. We use proportions
\[
A=\frac{n_a}{n},\qquad q=\frac{n_w}{n_a}\quad(n_a>0),\qquad u=\frac{n_w}{n},
\qquad F_{10}=\frac{n_{10}}n,\qquad V_{10}=\frac{n_{a,10}}n.
\]
When $n_a>0$, $u=Aq$. At zero acceptance, $A=u=0$ and $q$ is undefined. Rate columns marked \% display proportions multiplied by 100. Counts and percentage-point differences retain their stated units. Bare weight-function names in figures denote the corresponding IRLS rules.

The seven non-WLS weighting configurations comprise four IRLS rules and three component counts of one GMM family. Pointwise bootstrap intervals describe the selected designs. Section~8 gives the historical settings.

\paragraph{Electronic results.}
Paths begin in the data directory inside the archive. The sparse-injection recovery summaries are in \path{revision_20261102/results/competitive_correlated/summary.csv}, with window-level vectors in the adjacent \path{all_trials.csv}. Individual baseline-error summaries and block-length sensitivities are under \path{revision_20261102/results/sparse_acceptance/}.
\begin{itemize}
\item Main-text Section 5.2: \path{revision_20261103/results/sparse_acceptance.csv} (432 policies) and \path{revision_20261103/results/threshold_learning_intervals.csv} (2,016 records).
\item Main-text Section 5.3: under \path{revision_20261103/third_round/results/}, use \path{whitening_recovery_summary.csv} for per-cell rates, \path{whitening_recovery_intervals.csv} for 3,240 paired intervals, \path{whitening_threshold_learning.csv} for 3,456 threshold-learning results and \path{whitening_threshold_migration.csv} for threshold transfer.
\item Main-text Section 5.4: \path{revision_20261103/third_round/results/labels/acceptance_by_label.csv} (432 policies).
\end{itemize}
\footnotesize
\begingroup
\footnotesize
\setlength{\tabcolsep}{2pt}
\begin{longtable}{llrrrrl}
\caption{Random benchmark: full-vector recovery and paired changes relative to WLS. Change is $100(n_{\rm rescue}-n_{\rm loss})/n$ for $n$ planned pairs. Pointwise 95\% percentile intervals use 10,000 paired trial resamples. A dash in the interval denotes zero discordant pairs. The archive gives a common-completed sensitivity.}\label{tab:random-unified-full}\\
\toprule
Noise & Estimator & Correct/planned & \shortstack{Completed\\searches} & Rescue\% & Loss\% & \shortstack{Change vs WLS\\(pp) [95\% interval]} \\
\midrule
\endfirsthead
\multicolumn{7}{l}{\tablename~\thetable\ (continued)} \\
\toprule
Noise & Estimator & Correct/planned & \shortstack{Completed\\searches} & Rescue\% & Loss\% & \shortstack{Change vs WLS\\(pp) [95\% interval]} \\
\midrule
\endhead
Cauchy & WLS & 2/1000 & 1000 & 0.0 & 0.0 & 0.0 [--,--] \\
Cauchy & GMM, \(K=2\) & 490/1000 & 1000 & 48.8 & 0.0 & 48.8 [45.7,51.9] \\
Cauchy & GMM, \(K=3\) & 615/1000 & 1000 & 61.3 & 0.0 & 61.3 [58.3,64.3] \\
Cauchy & GMM, \(K=4\) & 672/1000 & 1000 & 67.0 & 0.0 & 67.0 [64.1,69.9] \\
Cauchy & Huber IRLS & 978/1000 & 1000 & 97.6 & 0.0 & 97.6 [96.6,98.5] \\
Cauchy & Cauchy IRLS & 981/1000 & 1000 & 97.9 & 0.0 & 97.9 [97.0,98.7] \\
Cauchy & IGG~III IRLS & 982/1000 & 1000 & 98.0 & 0.0 & 98.0 [97.1,98.8] \\
Cauchy & Tukey IRLS & 981/1000 & 1000 & 97.9 & 0.0 & 97.9 [97.0,98.8] \\
Gaussian & WLS & 984/1000 & 1000 & 0.0 & 0.0 & 0.0 [--,--] \\
Gaussian & GMM, \(K=2\) & 956/1000 & 1000 & 0.8 & 3.6 & -2.8 [-4.1,-1.5] \\
Gaussian & GMM, \(K=3\) & 955/1000 & 1000 & 0.8 & 3.7 & -2.9 [-4.2,-1.6] \\
Gaussian & GMM, \(K=4\) & 955/1000 & 1000 & 0.8 & 3.7 & -2.9 [-4.2,-1.6] \\
Gaussian & Huber IRLS & 977/1000 & 1000 & 0.7 & 1.4 & -0.7 [-1.6,0.2] \\
Gaussian & Cauchy IRLS & 972/1000 & 1000 & 0.8 & 2.0 & -1.2 [-2.3,-0.2] \\
Gaussian & IGG~III IRLS & 963/1000 & 1000 & 0.8 & 2.9 & -2.1 [-3.3,-0.9] \\
Gaussian & Tukey IRLS & 975/1000 & 1000 & 0.7 & 1.6 & -0.9 [-1.8,0.0] \\
Mixture & WLS & 41/1000 & 1000 & 0.0 & 0.0 & 0.0 [--,--] \\
Mixture & GMM, \(K=2\) & 756/1000 & 1000 & 71.9 & 0.4 & 71.5 [68.7,74.3] \\
Mixture & GMM, \(K=3\) & 874/1000 & 1000 & 83.4 & 0.1 & 83.3 [80.9,85.6] \\
Mixture & GMM, \(K=4\) & 874/1000 & 1000 & 83.4 & 0.1 & 83.3 [81.0,85.6] \\
Mixture & Huber IRLS & 898/1000 & 1000 & 85.7 & 0.0 & 85.7 [83.4,87.8] \\
Mixture & Cauchy IRLS & 905/1000 & 1000 & 86.5 & 0.1 & 86.4 [84.2,88.5] \\
Mixture & IGG~III IRLS & 943/1000 & 1000 & 90.3 & 0.1 & 90.2 [88.3,92.0] \\
Mixture & Tukey IRLS & 940/1000 & 1000 & 90.0 & 0.1 & 89.9 [88.0,91.7] \\
Student-\(t\) & WLS & 859/1000 & 1000 & 0.0 & 0.0 & 0.0 [--,--] \\
Student-\(t\) & GMM, \(K=2\) & 969/1000 & 1000 & 11.9 & 0.9 & 11.0 [8.9,13.1] \\
Student-\(t\) & GMM, \(K=3\) & 970/1000 & 1000 & 12.0 & 0.9 & 11.1 [9.0,13.2] \\
Student-\(t\) & GMM, \(K=4\) & 970/1000 & 1000 & 12.0 & 0.9 & 11.1 [9.0,13.2] \\
Student-\(t\) & Huber IRLS & 978/1000 & 1000 & 12.2 & 0.3 & 11.9 [9.9,14.1] \\
Student-\(t\) & Cauchy IRLS & 974/1000 & 1000 & 12.1 & 0.6 & 11.5 [9.5,13.6] \\
Student-\(t\) & IGG~III IRLS & 973/1000 & 1000 & 12.2 & 0.8 & 11.4 [9.2,13.6] \\
Student-\(t\) & Tukey IRLS & 977/1000 & 1000 & 12.3 & 0.5 & 11.8 [9.7,13.9] \\
\bottomrule
\end{longtable}
\endgroup

\begingroup
\footnotesize
\setlength{\tabcolsep}{2pt}
\begin{longtable}{llllrrrrl}
\caption{Between-satellite benchmark: full-vector recovery and paired changes relative to WLS by law, geometry and local bias. Fav. denotes favourable geometry. Bias is $b$ times the nominal undifferenced scale for each observation type. Change is $100(n_{\rm rescue}-n_{\rm loss})/n$ for $n$ planned pairs. Pointwise 95\% percentile intervals use 10,000 paired trial resamples. A dash in the interval denotes zero discordant pairs. The archive gives a common-completed sensitivity.}\label{tab:core-unified-full}\\
\toprule
Noise & Geom. & \(b\) & Estimator & Correct/planned & \shortstack{Completed\\searches} & Rescue\% & Loss\% & \shortstack{Change vs WLS\\(pp) [95\% interval]} \\
\midrule
\endfirsthead
\multicolumn{9}{l}{\tablename~\thetable\ (continued)} \\
\toprule
Noise & Geom. & \(b\) & Estimator & Correct/planned & \shortstack{Completed\\searches} & Rescue\% & Loss\% & \shortstack{Change vs WLS\\(pp) [95\% interval]} \\
\midrule
\endhead
Cauchy & Fav. & 0 & WLS & 125/1000 & 1000 & 0.0 & 0.0 & 0.0 [--,--] \\
Cauchy & Fav. & 0 & GMM, \(K=2\) & 768/1000 & 1000 & 64.8 & 0.5 & 64.3 [61.3,67.3] \\
Cauchy & Fav. & 0 & GMM, \(K=3\) & 770/1000 & 1000 & 65.0 & 0.5 & 64.5 [61.4,67.5] \\
Cauchy & Fav. & 0 & GMM, \(K=4\) & 769/1000 & 1000 & 64.9 & 0.5 & 64.4 [61.4,67.5] \\
Cauchy & Fav. & 0 & Huber IRLS & 803/1000 & 1000 & 68.2 & 0.4 & 67.8 [64.9,70.7] \\
Cauchy & Fav. & 0 & Cauchy IRLS & 817/1000 & 1000 & 69.6 & 0.4 & 69.2 [66.3,72.0] \\
Cauchy & Fav. & 0 & IGG~III IRLS & 828/1000 & 1000 & 70.7 & 0.4 & 70.3 [67.3,73.2] \\
Cauchy & Fav. & 0 & Tukey IRLS & 835/1000 & 1000 & 71.4 & 0.4 & 71.0 [68.1,73.8] \\
Cauchy & Fav. & 5 & WLS & 101/1000 & 1000 & 0.0 & 0.0 & 0.0 [--,--] \\
Cauchy & Fav. & 5 & GMM, \(K=2\) & 497/1000 & 1000 & 41.1 & 1.5 & 39.6 [36.3,42.7] \\
Cauchy & Fav. & 5 & GMM, \(K=3\) & 497/1000 & 1000 & 41.1 & 1.5 & 39.6 [36.4,42.8] \\
Cauchy & Fav. & 5 & GMM, \(K=4\) & 496/1000 & 1000 & 41.0 & 1.5 & 39.5 [36.3,42.7] \\
Cauchy & Fav. & 5 & Huber IRLS & 545/1000 & 1000 & 45.5 & 1.1 & 44.4 [41.2,47.6] \\
Cauchy & Fav. & 5 & Cauchy IRLS & 566/1000 & 1000 & 47.6 & 1.1 & 46.5 [43.3,49.7] \\
Cauchy & Fav. & 5 & IGG~III IRLS & 554/1000 & 1000 & 46.8 & 1.5 & 45.3 [41.9,48.5] \\
Cauchy & Fav. & 5 & Tukey IRLS & 577/1000 & 1000 & 48.9 & 1.3 & 47.6 [44.2,50.8] \\
Cauchy & Fav. & 10 & WLS & 29/1000 & 1000 & 0.0 & 0.0 & 0.0 [--,--] \\
Cauchy & Fav. & 10 & GMM, \(K=2\) & 332/1000 & 1000 & 31.0 & 0.7 & 30.3 [27.4,33.3] \\
Cauchy & Fav. & 10 & GMM, \(K=3\) & 332/1000 & 1000 & 31.0 & 0.7 & 30.3 [27.4,33.2] \\
Cauchy & Fav. & 10 & GMM, \(K=4\) & 332/1000 & 1000 & 31.0 & 0.7 & 30.3 [27.3,33.3] \\
Cauchy & Fav. & 10 & Huber IRLS & 337/1000 & 1000 & 31.2 & 0.4 & 30.8 [27.8,33.7] \\
Cauchy & Fav. & 10 & Cauchy IRLS & 354/1000 & 1000 & 32.9 & 0.4 & 32.5 [29.5,35.5] \\
Cauchy & Fav. & 10 & IGG~III IRLS & 415/1000 & 1000 & 39.1 & 0.5 & 38.6 [35.5,41.6] \\
Cauchy & Fav. & 10 & Tukey IRLS & 367/1000 & 1000 & 34.3 & 0.5 & 33.8 [30.8,36.8] \\
Cauchy & Weak & 0 & WLS & 10/1000 & 1000 & 0.0 & 0.0 & 0.0 [--,--] \\
Cauchy & Weak & 0 & GMM, \(K=2\) & 288/1000 & 1000 & 28.0 & 0.2 & 27.8 [25.0,30.6] \\
Cauchy & Weak & 0 & GMM, \(K=3\) & 289/1000 & 1000 & 28.1 & 0.2 & 27.9 [25.1,30.7] \\
Cauchy & Weak & 0 & GMM, \(K=4\) & 289/1000 & 1000 & 28.1 & 0.2 & 27.9 [25.1,30.7] \\
Cauchy & Weak & 0 & Huber IRLS & 291/1000 & 1000 & 28.2 & 0.1 & 28.1 [25.3,30.9] \\
Cauchy & Weak & 0 & Cauchy IRLS & 328/1000 & 1000 & 31.9 & 0.1 & 31.8 [28.9,34.7] \\
Cauchy & Weak & 0 & IGG~III IRLS & 354/1000 & 1000 & 34.6 & 0.2 & 34.4 [31.4,37.4] \\
Cauchy & Weak & 0 & Tukey IRLS & 348/1000 & 1000 & 34.1 & 0.3 & 33.8 [30.9,36.8] \\
Cauchy & Weak & 5 & WLS & 7/1000 & 1000 & 0.0 & 0.0 & 0.0 [--,--] \\
Cauchy & Weak & 5 & GMM, \(K=2\) & 93/1000 & 999 & 8.9 & 0.3 & 8.6 [6.8,10.4] \\
Cauchy & Weak & 5 & GMM, \(K=3\) & 93/1000 & 999 & 8.9 & 0.3 & 8.6 [6.8,10.4] \\
Cauchy & Weak & 5 & GMM, \(K=4\) & 93/1000 & 999 & 8.9 & 0.3 & 8.6 [6.9,10.4] \\
Cauchy & Weak & 5 & Huber IRLS & 85/1000 & 1000 & 7.9 & 0.1 & 7.8 [6.2,9.5] \\
Cauchy & Weak & 5 & Cauchy IRLS & 98/1000 & 1000 & 9.3 & 0.2 & 9.1 [7.3,11.0] \\
Cauchy & Weak & 5 & IGG~III IRLS & 110/1000 & 1000 & 10.5 & 0.2 & 10.3 [8.4,12.3] \\
Cauchy & Weak & 5 & Tukey IRLS & 101/1000 & 1000 & 9.6 & 0.2 & 9.4 [7.6,11.3] \\
Cauchy & Weak & 10 & WLS & 0/1000 & 1000 & 0.0 & 0.0 & 0.0 [--,--] \\
Cauchy & Weak & 10 & GMM, \(K=2\) & 47/1000 & 1000 & 4.7 & 0.0 & 4.7 [3.4,6.0] \\
Cauchy & Weak & 10 & GMM, \(K=3\) & 47/1000 & 1000 & 4.7 & 0.0 & 4.7 [3.4,6.1] \\
Cauchy & Weak & 10 & GMM, \(K=4\) & 47/1000 & 1000 & 4.7 & 0.0 & 4.7 [3.4,6.0] \\
Cauchy & Weak & 10 & Huber IRLS & 37/1000 & 1000 & 3.7 & 0.0 & 3.7 [2.6,4.9] \\
Cauchy & Weak & 10 & Cauchy IRLS & 50/1000 & 1000 & 5.0 & 0.0 & 5.0 [3.7,6.4] \\
Cauchy & Weak & 10 & IGG~III IRLS & 105/1000 & 1000 & 10.5 & 0.0 & 10.5 [8.7,12.5] \\
Cauchy & Weak & 10 & Tukey IRLS & 77/1000 & 1000 & 7.7 & 0.0 & 7.7 [6.1,9.4] \\
Gaussian & Fav. & 0 & WLS & 1000/1000 & 1000 & 0.0 & 0.0 & 0.0 [--,--] \\
Gaussian & Fav. & 0 & GMM, \(K=2\) & 999/1000 & 1000 & 0.0 & 0.1 & -0.1 [-0.3,0.0] \\
Gaussian & Fav. & 0 & GMM, \(K=3\) & 999/1000 & 1000 & 0.0 & 0.1 & -0.1 [-0.3,0.0] \\
Gaussian & Fav. & 0 & GMM, \(K=4\) & 999/1000 & 1000 & 0.0 & 0.1 & -0.1 [-0.3,0.0] \\
Gaussian & Fav. & 0 & Huber IRLS & 999/1000 & 1000 & 0.0 & 0.1 & -0.1 [-0.3,0.0] \\
Gaussian & Fav. & 0 & Cauchy IRLS & 999/1000 & 1000 & 0.0 & 0.1 & -0.1 [-0.3,0.0] \\
Gaussian & Fav. & 0 & IGG~III IRLS & 999/1000 & 1000 & 0.0 & 0.1 & -0.1 [-0.3,0.0] \\
Gaussian & Fav. & 0 & Tukey IRLS & 999/1000 & 1000 & 0.0 & 0.1 & -0.1 [-0.3,0.0] \\
Gaussian & Fav. & 5 & WLS & 951/1000 & 1000 & 0.0 & 0.0 & 0.0 [--,--] \\
Gaussian & Fav. & 5 & GMM, \(K=2\) & 976/1000 & 1000 & 3.6 & 1.1 & 2.5 [1.2,3.8] \\
Gaussian & Fav. & 5 & GMM, \(K=3\) & 976/1000 & 1000 & 3.6 & 1.1 & 2.5 [1.2,3.8] \\
Gaussian & Fav. & 5 & GMM, \(K=4\) & 976/1000 & 1000 & 3.6 & 1.1 & 2.5 [1.2,3.9] \\
Gaussian & Fav. & 5 & Huber IRLS & 979/1000 & 1000 & 3.6 & 0.8 & 2.8 [1.6,4.1] \\
Gaussian & Fav. & 5 & Cauchy IRLS & 978/1000 & 1000 & 3.6 & 0.9 & 2.7 [1.4,4.0] \\
Gaussian & Fav. & 5 & IGG~III IRLS & 983/1000 & 1000 & 4.0 & 0.8 & 3.2 [1.9,4.5] \\
Gaussian & Fav. & 5 & Tukey IRLS & 978/1000 & 1000 & 3.6 & 0.9 & 2.7 [1.4,4.0] \\
Gaussian & Fav. & 10 & WLS & 301/1000 & 1000 & 0.0 & 0.0 & 0.0 [--,--] \\
Gaussian & Fav. & 10 & GMM, \(K=2\) & 948/1000 & 1000 & 64.9 & 0.2 & 64.7 [61.7,67.7] \\
Gaussian & Fav. & 10 & GMM, \(K=3\) & 948/1000 & 1000 & 64.9 & 0.2 & 64.7 [61.7,67.7] \\
Gaussian & Fav. & 10 & GMM, \(K=4\) & 948/1000 & 1000 & 64.9 & 0.2 & 64.7 [61.7,67.7] \\
Gaussian & Fav. & 10 & Huber IRLS & 955/1000 & 1000 & 65.5 & 0.1 & 65.4 [62.4,68.3] \\
Gaussian & Fav. & 10 & Cauchy IRLS & 954/1000 & 1000 & 65.5 & 0.2 & 65.3 [62.2,68.3] \\
Gaussian & Fav. & 10 & IGG~III IRLS & 959/1000 & 1000 & 66.0 & 0.2 & 65.8 [62.7,68.8] \\
Gaussian & Fav. & 10 & Tukey IRLS & 962/1000 & 1000 & 66.4 & 0.3 & 66.1 [63.2,69.0] \\
Gaussian & Weak & 0 & WLS & 939/1000 & 1000 & 0.0 & 0.0 & 0.0 [--,--] \\
Gaussian & Weak & 0 & GMM, \(K=2\) & 839/1000 & 1000 & 1.5 & 11.5 & -10.0 [-12.2,-7.9] \\
Gaussian & Weak & 0 & GMM, \(K=3\) & 839/1000 & 1000 & 1.5 & 11.5 & -10.0 [-12.2,-7.9] \\
Gaussian & Weak & 0 & GMM, \(K=4\) & 839/1000 & 1000 & 1.5 & 11.5 & -10.0 [-12.2,-7.9] \\
Gaussian & Weak & 0 & Huber IRLS & 898/1000 & 1000 & 1.3 & 5.4 & -4.1 [-5.7,-2.5] \\
Gaussian & Weak & 0 & Cauchy IRLS & 881/1000 & 1000 & 1.4 & 7.2 & -5.8 [-7.6,-4.0] \\
Gaussian & Weak & 0 & IGG~III IRLS & 855/1000 & 1000 & 1.4 & 9.8 & -8.4 [-10.4,-6.4] \\
Gaussian & Weak & 0 & Tukey IRLS & 888/1000 & 1000 & 1.4 & 6.5 & -5.1 [-6.9,-3.4] \\
Gaussian & Weak & 5 & WLS & 196/1000 & 1000 & 0.0 & 0.0 & 0.0 [--,--] \\
Gaussian & Weak & 5 & GMM, \(K=2\) & 525/1000 & 1000 & 37.0 & 4.1 & 32.9 [29.5,36.4] \\
Gaussian & Weak & 5 & GMM, \(K=3\) & 525/1000 & 1000 & 37.0 & 4.1 & 32.9 [29.5,36.3] \\
Gaussian & Weak & 5 & GMM, \(K=4\) & 525/1000 & 1000 & 37.0 & 4.1 & 32.9 [29.5,36.3] \\
Gaussian & Weak & 5 & Huber IRLS & 553/1000 & 1000 & 38.2 & 2.5 & 35.7 [32.5,38.9] \\
Gaussian & Weak & 5 & Cauchy IRLS & 541/1000 & 1000 & 37.9 & 3.4 & 34.5 [31.1,37.8] \\
Gaussian & Weak & 5 & IGG~III IRLS & 595/1000 & 1000 & 43.8 & 3.9 & 39.9 [36.4,43.4] \\
Gaussian & Weak & 5 & Tukey IRLS & 551/1000 & 1000 & 39.2 & 3.7 & 35.5 [32.1,38.8] \\
Gaussian & Weak & 10 & WLS & 0/1000 & 1000 & 0.0 & 0.0 & 0.0 [--,--] \\
Gaussian & Weak & 10 & GMM, \(K=2\) & 378/1000 & 1000 & 37.8 & 0.0 & 37.8 [34.8,40.8] \\
Gaussian & Weak & 10 & GMM, \(K=3\) & 378/1000 & 1000 & 37.8 & 0.0 & 37.8 [34.8,40.8] \\
Gaussian & Weak & 10 & GMM, \(K=4\) & 378/1000 & 1000 & 37.8 & 0.0 & 37.8 [34.8,40.8] \\
Gaussian & Weak & 10 & Huber IRLS & 377/1000 & 1000 & 37.7 & 0.0 & 37.7 [34.6,40.7] \\
Gaussian & Weak & 10 & Cauchy IRLS & 428/1000 & 1000 & 42.8 & 0.0 & 42.8 [39.7,45.9] \\
Gaussian & Weak & 10 & IGG~III IRLS & 522/1000 & 1000 & 52.2 & 0.0 & 52.2 [49.2,55.3] \\
Gaussian & Weak & 10 & Tukey IRLS & 576/1000 & 1000 & 57.6 & 0.0 & 57.6 [54.5,60.7] \\
Mixture & Fav. & 0 & WLS & 205/1000 & 1000 & 0.0 & 0.0 & 0.0 [--,--] \\
Mixture & Fav. & 0 & GMM, \(K=2\) & 734/1000 & 1000 & 54.2 & 1.3 & 52.9 [49.7,56.1] \\
Mixture & Fav. & 0 & GMM, \(K=3\) & 734/1000 & 1000 & 54.2 & 1.3 & 52.9 [49.7,56.1] \\
Mixture & Fav. & 0 & GMM, \(K=4\) & 734/1000 & 1000 & 54.2 & 1.3 & 52.9 [49.6,56.2] \\
Mixture & Fav. & 0 & Huber IRLS & 751/1000 & 1000 & 55.2 & 0.6 & 54.6 [51.5,57.8] \\
Mixture & Fav. & 0 & Cauchy IRLS & 756/1000 & 1000 & 55.9 & 0.8 & 55.1 [51.9,58.3] \\
Mixture & Fav. & 0 & IGG~III IRLS & 788/1000 & 1000 & 59.2 & 0.9 & 58.3 [55.1,61.5] \\
Mixture & Fav. & 0 & Tukey IRLS & 779/1000 & 1000 & 58.3 & 0.9 & 57.4 [54.2,60.7] \\
Mixture & Fav. & 5 & WLS & 139/1000 & 1000 & 0.0 & 0.0 & 0.0 [--,--] \\
Mixture & Fav. & 5 & GMM, \(K=2\) & 459/1000 & 1000 & 35.1 & 3.1 & 32.0 [28.7,35.3] \\
Mixture & Fav. & 5 & GMM, \(K=3\) & 459/1000 & 1000 & 35.1 & 3.1 & 32.0 [28.7,35.3] \\
Mixture & Fav. & 5 & GMM, \(K=4\) & 459/1000 & 1000 & 35.1 & 3.1 & 32.0 [28.7,35.3] \\
Mixture & Fav. & 5 & Huber IRLS & 505/1000 & 1000 & 38.7 & 2.1 & 36.6 [33.4,39.8] \\
Mixture & Fav. & 5 & Cauchy IRLS & 498/1000 & 1000 & 38.1 & 2.2 & 35.9 [32.6,39.1] \\
Mixture & Fav. & 5 & IGG~III IRLS & 502/1000 & 1000 & 39.5 & 3.2 & 36.3 [33.0,39.7] \\
Mixture & Fav. & 5 & Tukey IRLS & 513/1000 & 1000 & 39.9 & 2.5 & 37.4 [34.1,40.7] \\
Mixture & Fav. & 10 & WLS & 55/1000 & 1000 & 0.0 & 0.0 & 0.0 [--,--] \\
Mixture & Fav. & 10 & GMM, \(K=2\) & 271/1000 & 1000 & 24.5 & 2.9 & 21.6 [18.6,24.5] \\
Mixture & Fav. & 10 & GMM, \(K=3\) & 271/1000 & 1000 & 24.5 & 2.9 & 21.6 [18.6,24.6] \\
Mixture & Fav. & 10 & GMM, \(K=4\) & 271/1000 & 1000 & 24.5 & 2.9 & 21.6 [18.6,24.6] \\
Mixture & Fav. & 10 & Huber IRLS & 267/1000 & 1000 & 23.5 & 2.3 & 21.2 [18.4,24.1] \\
Mixture & Fav. & 10 & Cauchy IRLS & 276/1000 & 1000 & 24.7 & 2.6 & 22.1 [19.2,25.1] \\
Mixture & Fav. & 10 & IGG~III IRLS & 347/1000 & 1000 & 32.1 & 2.9 & 29.2 [26.0,32.4] \\
Mixture & Fav. & 10 & Tukey IRLS & 288/1000 & 1000 & 26.0 & 2.7 & 23.3 [20.3,26.3] \\
Mixture & Weak & 0 & WLS & 6/1000 & 1000 & 0.0 & 0.0 & 0.0 [--,--] \\
Mixture & Weak & 0 & GMM, \(K=2\) & 196/1000 & 1000 & 19.2 & 0.2 & 19.0 [16.6,21.5] \\
Mixture & Weak & 0 & GMM, \(K=3\) & 196/1000 & 1000 & 19.2 & 0.2 & 19.0 [16.6,21.5] \\
Mixture & Weak & 0 & GMM, \(K=4\) & 196/1000 & 1000 & 19.2 & 0.2 & 19.0 [16.6,21.4] \\
Mixture & Weak & 0 & Huber IRLS & 192/1000 & 1000 & 18.7 & 0.1 & 18.6 [16.2,21.0] \\
Mixture & Weak & 0 & Cauchy IRLS & 201/1000 & 1000 & 19.6 & 0.1 & 19.5 [17.1,22.0] \\
Mixture & Weak & 0 & IGG~III IRLS & 245/1000 & 1000 & 24.1 & 0.2 & 23.9 [21.3,26.6] \\
Mixture & Weak & 0 & Tukey IRLS & 250/1000 & 1000 & 24.5 & 0.1 & 24.4 [21.8,27.0] \\
Mixture & Weak & 5 & WLS & 4/1000 & 1000 & 0.0 & 0.0 & 0.0 [--,--] \\
Mixture & Weak & 5 & GMM, \(K=2\) & 49/1000 & 1000 & 4.8 & 0.3 & 4.5 [3.2,5.9] \\
Mixture & Weak & 5 & GMM, \(K=3\) & 49/1000 & 1000 & 4.8 & 0.3 & 4.5 [3.2,5.9] \\
Mixture & Weak & 5 & GMM, \(K=4\) & 49/1000 & 1000 & 4.8 & 0.3 & 4.5 [3.1,5.9] \\
Mixture & Weak & 5 & Huber IRLS & 57/1000 & 1000 & 5.5 & 0.2 & 5.3 [3.9,6.8] \\
Mixture & Weak & 5 & Cauchy IRLS & 57/1000 & 1000 & 5.6 & 0.3 & 5.3 [3.8,6.8] \\
Mixture & Weak & 5 & IGG~III IRLS & 65/1000 & 1000 & 6.4 & 0.3 & 6.1 [4.6,7.7] \\
Mixture & Weak & 5 & Tukey IRLS & 54/1000 & 1000 & 5.3 & 0.3 & 5.0 [3.6,6.5] \\
Mixture & Weak & 10 & WLS & 0/1000 & 1000 & 0.0 & 0.0 & 0.0 [--,--] \\
Mixture & Weak & 10 & GMM, \(K=2\) & 32/1000 & 1000 & 3.2 & 0.0 & 3.2 [2.2,4.3] \\
Mixture & Weak & 10 & GMM, \(K=3\) & 32/1000 & 1000 & 3.2 & 0.0 & 3.2 [2.1,4.3] \\
Mixture & Weak & 10 & GMM, \(K=4\) & 32/1000 & 1000 & 3.2 & 0.0 & 3.2 [2.1,4.3] \\
Mixture & Weak & 10 & Huber IRLS & 27/1000 & 1000 & 2.7 & 0.0 & 2.7 [1.7,3.7] \\
Mixture & Weak & 10 & Cauchy IRLS & 32/1000 & 1000 & 3.2 & 0.0 & 3.2 [2.2,4.3] \\
Mixture & Weak & 10 & IGG~III IRLS & 87/1000 & 1000 & 8.7 & 0.0 & 8.7 [7.0,10.5] \\
Mixture & Weak & 10 & Tukey IRLS & 51/1000 & 1000 & 5.1 & 0.0 & 5.1 [3.8,6.5] \\
Student-\(t\) & Fav. & 0 & WLS & 973/1000 & 1000 & 0.0 & 0.0 & 0.0 [--,--] \\
Student-\(t\) & Fav. & 0 & GMM, \(K=2\) & 992/1000 & 1000 & 2.3 & 0.4 & 1.9 [0.9,2.9] \\
Student-\(t\) & Fav. & 0 & GMM, \(K=3\) & 992/1000 & 1000 & 2.3 & 0.4 & 1.9 [0.9,2.9] \\
Student-\(t\) & Fav. & 0 & GMM, \(K=4\) & 992/1000 & 1000 & 2.3 & 0.4 & 1.9 [0.9,2.9] \\
Student-\(t\) & Fav. & 0 & Huber IRLS & 995/1000 & 1000 & 2.3 & 0.1 & 2.2 [1.3,3.2] \\
Student-\(t\) & Fav. & 0 & Cauchy IRLS & 994/1000 & 1000 & 2.3 & 0.2 & 2.1 [1.2,3.1] \\
Student-\(t\) & Fav. & 0 & IGG~III IRLS & 993/1000 & 1000 & 2.3 & 0.3 & 2.0 [1.0,3.0] \\
Student-\(t\) & Fav. & 0 & Tukey IRLS & 994/1000 & 1000 & 2.3 & 0.2 & 2.1 [1.2,3.1] \\
Student-\(t\) & Fav. & 5 & WLS & 859/1000 & 1000 & 0.0 & 0.0 & 0.0 [--,--] \\
Student-\(t\) & Fav. & 5 & GMM, \(K=2\) & 877/1000 & 1000 & 8.6 & 6.8 & 1.8 [-0.6,4.2] \\
Student-\(t\) & Fav. & 5 & GMM, \(K=3\) & 877/1000 & 1000 & 8.6 & 6.8 & 1.8 [-0.6,4.3] \\
Student-\(t\) & Fav. & 5 & GMM, \(K=4\) & 877/1000 & 1000 & 8.6 & 6.8 & 1.8 [-0.6,4.2] \\
Student-\(t\) & Fav. & 5 & Huber IRLS & 900/1000 & 1000 & 9.0 & 4.9 & 4.1 [1.8,6.4] \\
Student-\(t\) & Fav. & 5 & Cauchy IRLS & 891/1000 & 1000 & 8.9 & 5.7 & 3.2 [0.9,5.6] \\
Student-\(t\) & Fav. & 5 & IGG~III IRLS & 888/1000 & 1000 & 9.1 & 6.2 & 2.9 [0.5,5.3] \\
Student-\(t\) & Fav. & 5 & Tukey IRLS & 891/1000 & 1000 & 8.9 & 5.7 & 3.2 [0.9,5.6] \\
Student-\(t\) & Fav. & 10 & WLS & 232/1000 & 1000 & 0.0 & 0.0 & 0.0 [--,--] \\
Student-\(t\) & Fav. & 10 & GMM, \(K=2\) & 796/1000 & 1000 & 57.9 & 1.5 & 56.4 [53.0,59.6] \\
Student-\(t\) & Fav. & 10 & GMM, \(K=3\) & 796/1000 & 1000 & 57.9 & 1.5 & 56.4 [53.1,59.6] \\
Student-\(t\) & Fav. & 10 & GMM, \(K=4\) & 796/1000 & 1000 & 57.9 & 1.5 & 56.4 [53.1,59.6] \\
Student-\(t\) & Fav. & 10 & Huber IRLS & 821/1000 & 1000 & 60.0 & 1.1 & 58.9 [55.7,62.1] \\
Student-\(t\) & Fav. & 10 & Cauchy IRLS & 822/1000 & 1000 & 60.0 & 1.0 & 59.0 [55.9,62.2] \\
Student-\(t\) & Fav. & 10 & IGG~III IRLS & 852/1000 & 1000 & 63.7 & 1.7 & 62.0 [58.7,65.1] \\
Student-\(t\) & Fav. & 10 & Tukey IRLS & 844/1000 & 1000 & 62.7 & 1.5 & 61.2 [57.9,64.4] \\
Student-\(t\) & Weak & 0 & WLS & 644/1000 & 1000 & 0.0 & 0.0 & 0.0 [--,--] \\
Student-\(t\) & Weak & 0 & GMM, \(K=2\) & 738/1000 & 1000 & 18.2 & 8.8 & 9.4 [6.2,12.6] \\
Student-\(t\) & Weak & 0 & GMM, \(K=3\) & 738/1000 & 1000 & 18.2 & 8.8 & 9.4 [6.2,12.6] \\
Student-\(t\) & Weak & 0 & GMM, \(K=4\) & 738/1000 & 1000 & 18.2 & 8.8 & 9.4 [6.3,12.6] \\
Student-\(t\) & Weak & 0 & Huber IRLS & 795/1000 & 1000 & 19.0 & 3.9 & 15.1 [12.3,17.9] \\
Student-\(t\) & Weak & 0 & Cauchy IRLS & 777/1000 & 1000 & 19.0 & 5.7 & 13.3 [10.3,16.2] \\
Student-\(t\) & Weak & 0 & IGG~III IRLS & 762/1000 & 1000 & 19.6 & 7.8 & 11.8 [8.7,14.9] \\
Student-\(t\) & Weak & 0 & Tukey IRLS & 782/1000 & 1000 & 19.5 & 5.7 & 13.8 [10.8,16.8] \\
Student-\(t\) & Weak & 5 & WLS & 130/1000 & 1000 & 0.0 & 0.0 & 0.0 [--,--] \\
Student-\(t\) & Weak & 5 & GMM, \(K=2\) & 340/1000 & 1000 & 24.9 & 3.9 & 21.0 [17.9,24.0] \\
Student-\(t\) & Weak & 5 & GMM, \(K=3\) & 340/1000 & 1000 & 24.9 & 3.9 & 21.0 [18.0,24.1] \\
Student-\(t\) & Weak & 5 & GMM, \(K=4\) & 340/1000 & 1000 & 24.9 & 3.9 & 21.0 [17.9,24.1] \\
Student-\(t\) & Weak & 5 & Huber IRLS & 362/1000 & 1000 & 26.2 & 3.0 & 23.2 [20.2,26.3] \\
Student-\(t\) & Weak & 5 & Cauchy IRLS & 360/1000 & 1000 & 26.1 & 3.1 & 23.0 [19.9,26.0] \\
Student-\(t\) & Weak & 5 & IGG~III IRLS & 403/1000 & 1000 & 31.2 & 3.9 & 27.3 [24.0,30.5] \\
Student-\(t\) & Weak & 5 & Tukey IRLS & 365/1000 & 1000 & 26.9 & 3.4 & 23.5 [20.5,26.5] \\
Student-\(t\) & Weak & 10 & WLS & 0/1000 & 1000 & 0.0 & 0.0 & 0.0 [--,--] \\
Student-\(t\) & Weak & 10 & GMM, \(K=2\) & 251/1000 & 1000 & 25.1 & 0.0 & 25.1 [22.4,27.8] \\
Student-\(t\) & Weak & 10 & GMM, \(K=3\) & 251/1000 & 1000 & 25.1 & 0.0 & 25.1 [22.5,27.8] \\
Student-\(t\) & Weak & 10 & GMM, \(K=4\) & 251/1000 & 1000 & 25.1 & 0.0 & 25.1 [22.4,27.8] \\
Student-\(t\) & Weak & 10 & Huber IRLS & 260/1000 & 1000 & 26.0 & 0.0 & 26.0 [23.4,28.7] \\
Student-\(t\) & Weak & 10 & Cauchy IRLS & 287/1000 & 1000 & 28.7 & 0.0 & 28.7 [26.0,31.5] \\
Student-\(t\) & Weak & 10 & IGG~III IRLS & 397/1000 & 1000 & 39.7 & 0.0 & 39.7 [36.7,42.7] \\
Student-\(t\) & Weak & 10 & Tukey IRLS & 410/1000 & 1000 & 41.0 & 0.0 & 41.0 [38.0,44.1] \\
\bottomrule
\end{longtable}
\endgroup

\begingroup
\footnotesize
\setlength{\tabcolsep}{2pt}
\begin{longtable}{lrrrr}
\caption{Full-vector recovery (\%) by local-bias duration in favourable geometry. The two selected satellites each receive a bias of five nominal scales for $d$ consecutive epochs. Each column uses 1,000 planned problems paired across configurations.}\label{tab:duration-unified-full}\\
\toprule
Estimator & Duration 1 & Duration 2 & Duration 3 & Duration 4 \\
\midrule
\endfirsthead
\multicolumn{5}{l}{\tablename~\thetable\ (continued)} \\
\toprule
Estimator & Duration 1 & Duration 2 & Duration 3 & Duration 4 \\
\midrule
\endhead
\multicolumn{5}{l}{Gaussian, \(N=5\); \(d=1,2,3,4\)} \\
WLS & 99.7 & 94.9 & 70.7 & 28.4 \\
GMM, \(K=2\) & 99.7 & 98.0 & 25.3 & 9.9 \\
GMM, \(K=3\) & 99.7 & 98.0 & 25.3 & 9.9 \\
GMM, \(K=4\) & 99.7 & 98.0 & 25.3 & 9.9 \\
Huber IRLS & 99.7 & 98.4 & 50.5 & 15.8 \\
Cauchy IRLS & 99.7 & 98.2 & 44.0 & 12.1 \\
IGG~III IRLS & 99.7 & 98.6 & 30.9 & 10.1 \\
Tukey IRLS & 99.7 & 98.4 & 45.4 & 13.1 \\
\midrule
\multicolumn{5}{l}{Gaussian, \(N=10\); \(d=2,3,4,6\)} \\
WLS & 99.9 & 99.4 & 97.5 & 76.4 \\
GMM, \(K=2\) & 99.9 & 99.9 & 99.5 & 34.9 \\
GMM, \(K=3\) & 99.9 & 99.9 & 99.5 & 34.9 \\
GMM, \(K=4\) & 99.9 & 99.9 & 99.5 & 34.9 \\
Huber IRLS & 99.9 & 99.9 & 99.4 & 59.7 \\
Cauchy IRLS & 99.9 & 99.9 & 99.5 & 50.8 \\
IGG~III IRLS & 99.9 & 100.0 & 99.7 & 39.5 \\
Tukey IRLS & 99.9 & 99.9 & 99.5 & 54.4 \\
\midrule
\multicolumn{5}{l}{Student-\(t_3\), \(N=5\); \(d=1,2,3,4\)} \\
WLS & 96.4 & 86.9 & 59.6 & 23.1 \\
GMM, \(K=2\) & 98.1 & 89.4 & 26.3 & 12.4 \\
GMM, \(K=3\) & 98.1 & 89.4 & 26.3 & 12.4 \\
GMM, \(K=4\) & 98.1 & 89.4 & 26.3 & 12.4 \\
Huber IRLS & 98.6 & 92.1 & 45.7 & 16.7 \\
Cauchy IRLS & 98.3 & 91.3 & 39.2 & 14.3 \\
IGG~III IRLS & 98.5 & 89.8 & 30.4 & 12.8 \\
Tukey IRLS & 98.2 & 90.6 & 42.0 & 14.8 \\
\midrule
\multicolumn{5}{l}{Student-\(t_3\), \(N=10\); \(d=2,3,4,6\)} \\
WLS & 99.0 & 97.8 & 94.3 & 69.6 \\
GMM, \(K=2\) & 99.9 & 99.5 & 97.9 & 33.8 \\
GMM, \(K=3\) & 99.9 & 99.5 & 97.9 & 33.8 \\
GMM, \(K=4\) & 99.9 & 99.5 & 97.9 & 33.8 \\
Huber IRLS & 99.9 & 99.6 & 98.3 & 51.6 \\
Cauchy IRLS & 99.9 & 99.6 & 98.1 & 44.4 \\
IGG~III IRLS & 99.9 & 99.5 & 98.4 & 39.4 \\
Tukey IRLS & 99.9 & 99.5 & 98.2 & 47.1 \\
\midrule
\multicolumn{5}{l}{Mixture, \(N=5\); \(d=1,2,3,4\)} \\
WLS & 19.0 & 15.1 & 10.5 & 4.8 \\
GMM, \(K=2\) & 63.4 & 45.5 & 19.2 & 9.0 \\
GMM, \(K=3\) & 63.4 & 45.5 & 19.2 & 9.0 \\
GMM, \(K=4\) & 63.4 & 45.5 & 19.2 & 9.0 \\
Huber IRLS & 66.5 & 49.4 & 24.2 & 11.4 \\
Cauchy IRLS & 66.3 & 49.1 & 22.7 & 10.3 \\
IGG~III IRLS & 67.2 & 52.1 & 23.3 & 10.7 \\
Tukey IRLS & 67.4 & 51.4 & 24.4 & 10.8 \\
\midrule
\multicolumn{5}{l}{Mixture, \(N=10\); \(d=2,3,4,6\)} \\
WLS & 50.0 & 41.9 & 34.4 & 23.8 \\
GMM, \(K=2\) & 93.0 & 87.4 & 75.6 & 35.6 \\
GMM, \(K=3\) & 93.0 & 87.4 & 75.6 & 35.6 \\
GMM, \(K=4\) & 93.0 & 87.4 & 75.6 & 35.6 \\
Huber IRLS & 93.4 & 87.8 & 76.1 & 41.1 \\
Cauchy IRLS & 93.3 & 87.9 & 76.2 & 38.4 \\
IGG~III IRLS & 93.8 & 88.3 & 77.1 & 38.6 \\
Tukey IRLS & 93.9 & 88.6 & 76.9 & 40.2 \\
\midrule
\multicolumn{5}{l}{Cauchy, \(N=5\); \(d=1,2,3,4\)} \\
WLS & 12.2 & 10.0 & 6.9 & 3.9 \\
GMM, \(K=2\) & 67.3 & 47.8 & 18.3 & 8.1 \\
GMM, \(K=3\) & 67.5 & 47.8 & 18.3 & 8.1 \\
GMM, \(K=4\) & 67.6 & 47.8 & 18.3 & 8.1 \\
Huber IRLS & 71.0 & 52.6 & 23.6 & 11.2 \\
Cauchy IRLS & 73.6 & 53.5 & 23.9 & 10.0 \\
IGG~III IRLS & 74.8 & 54.9 & 23.8 & 9.5 \\
Tukey IRLS & 75.7 & 55.3 & 25.8 & 9.9 \\
\midrule
\multicolumn{5}{l}{Cauchy, \(N=10\); \(d=2,3,4,6\)} \\
WLS & 12.9 & 11.6 & 9.3 & 6.4 \\
GMM, \(K=2\) & 89.8 & 84.1 & 74.7 & 34.9 \\
GMM, \(K=3\) & 90.5 & 84.6 & 75.2 & 35.4 \\
GMM, \(K=4\) & 90.7 & 85.0 & 75.6 & 35.4 \\
Huber IRLS & 96.0 & 91.1 & 80.5 & 43.4 \\
Cauchy IRLS & 95.9 & 91.5 & 80.7 & 41.2 \\
IGG~III IRLS & 96.3 & 91.5 & 80.9 & 41.8 \\
Tukey IRLS & 96.3 & 92.4 & 80.9 & 42.5 \\
\bottomrule
\end{longtable}
\endgroup

\begingroup
\footnotesize
\setlength{\tabcolsep}{2pt}
\begin{longtable}{lrrrr}
\caption{Paired changes in full-vector recovery relative to WLS, in percentage points. Brackets show pointwise 95\% percentile intervals from 10,000 paired trial bootstrap draws. Differences use all planned pairs. A dash indicates zero discordance.}\label{tab:duration-unified-paired}\\
\toprule
Estimator & Duration 1 & Duration 2 & Duration 3 & Duration 4 \\
\midrule
\endfirsthead
\multicolumn{5}{l}{\tablename~\thetable\ (continued)} \\
\toprule
Estimator & Duration 1 & Duration 2 & Duration 3 & Duration 4 \\
\midrule
\endhead
\multicolumn{5}{l}{Gaussian, \(N=5\); \(d=1,2,3,4\)} \\
GMM, \(K=2\) & 0.0 [-0.4,0.4] & 3.1 [1.8,4.4] & -45.4 [-48.5,-42.3] & -18.5 [-21.0,-16.1] \\
GMM, \(K=3\) & 0.0 [-0.4,0.4] & 3.1 [1.8,4.4] & -45.4 [-48.5,-42.3] & -18.5 [-21.0,-16.1] \\
GMM, \(K=4\) & 0.0 [-0.4,0.4] & 3.1 [1.8,4.4] & -45.4 [-48.5,-42.3] & -18.5 [-21.1,-16.1] \\
Huber IRLS & 0.0 [-0.4,0.4] & 3.5 [2.3,4.8] & -20.2 [-22.8,-17.7] & -12.6 [-14.8,-10.4] \\
Cauchy IRLS & 0.0 [-0.4,0.4] & 3.3 [2.0,4.6] & -26.7 [-29.6,-23.9] & -16.3 [-18.8,-13.9] \\
IGG~III IRLS & 0.0 [-0.4,0.4] & 3.7 [2.4,5.1] & -39.8 [-42.9,-36.6] & -18.3 [-20.8,-15.9] \\
Tukey IRLS & 0.0 [-0.4,0.4] & 3.5 [2.2,4.9] & -25.3 [-28.1,-22.6] & -15.3 [-17.7,-12.9] \\
\midrule
\multicolumn{5}{l}{Gaussian, \(N=10\); \(d=2,3,4,6\)} \\
GMM, \(K=2\) & 0.0 [--,--] & 0.5 [0.1,1.0] & 2.0 [1.1,3.0] & -41.5 [-44.6,-38.4] \\
GMM, \(K=3\) & 0.0 [--,--] & 0.5 [0.1,1.0] & 2.0 [1.1,3.0] & -41.5 [-44.5,-38.4] \\
GMM, \(K=4\) & 0.0 [--,--] & 0.5 [0.1,1.0] & 2.0 [1.1,3.0] & -41.5 [-44.6,-38.4] \\
Huber IRLS & 0.0 [--,--] & 0.5 [0.1,1.0] & 1.9 [1.0,2.9] & -16.7 [-19.1,-14.3] \\
Cauchy IRLS & 0.0 [--,--] & 0.5 [0.1,1.0] & 2.0 [1.1,3.0] & -25.6 [-28.4,-22.9] \\
IGG~III IRLS & 0.0 [--,--] & 0.6 [0.2,1.1] & 2.2 [1.2,3.2] & -36.9 [-40.0,-33.7] \\
Tukey IRLS & 0.0 [--,--] & 0.5 [0.1,1.0] & 2.0 [1.1,3.0] & -22.0 [-24.6,-19.3] \\
\midrule
\multicolumn{5}{l}{Student-\(t_3\), \(N=5\); \(d=1,2,3,4\)} \\
GMM, \(K=2\) & 1.7 [0.5,3.0] & 2.5 [0.1,4.9] & -33.3 [-36.5,-30.2] & -10.7 [-13.0,-8.4] \\
GMM, \(K=3\) & 1.7 [0.5,3.0] & 2.5 [0.1,4.9] & -33.3 [-36.5,-30.1] & -10.7 [-13.0,-8.3] \\
GMM, \(K=4\) & 1.7 [0.5,3.0] & 2.5 [0.1,4.9] & -33.3 [-36.5,-30.1] & -10.7 [-13.0,-8.4] \\
Huber IRLS & 2.2 [1.1,3.4] & 5.2 [3.0,7.4] & -13.9 [-16.6,-11.3] & -6.4 [-8.5,-4.2] \\
Cauchy IRLS & 1.9 [0.7,3.1] & 4.4 [2.1,6.7] & -20.4 [-23.3,-17.5] & -8.8 [-11.1,-6.5] \\
IGG~III IRLS & 2.1 [0.8,3.4] & 2.9 [0.5,5.4] & -29.2 [-32.5,-26.0] & -10.3 [-12.7,-8.0] \\
Tukey IRLS & 1.8 [0.6,3.0] & 3.7 [1.4,6.1] & -17.6 [-20.4,-14.8] & -8.3 [-10.5,-6.1] \\
\midrule
\multicolumn{5}{l}{Student-\(t_3\), \(N=10\); \(d=2,3,4,6\)} \\
GMM, \(K=2\) & 0.9 [0.4,1.5] & 1.7 [0.8,2.6] & 3.6 [2.3,5.0] & -35.8 [-38.9,-32.8] \\
GMM, \(K=3\) & 0.9 [0.4,1.5] & 1.7 [0.8,2.7] & 3.6 [2.3,5.0] & -35.8 [-38.9,-32.7] \\
GMM, \(K=4\) & 0.9 [0.4,1.5] & 1.7 [0.8,2.7] & 3.6 [2.3,5.0] & -35.8 [-38.9,-32.8] \\
Huber IRLS & 0.9 [0.4,1.5] & 1.8 [0.9,2.8] & 4.0 [2.7,5.4] & -18.0 [-20.6,-15.4] \\
Cauchy IRLS & 0.9 [0.4,1.5] & 1.8 [0.9,2.7] & 3.8 [2.5,5.2] & -25.2 [-28.2,-22.2] \\
IGG~III IRLS & 0.9 [0.4,1.5] & 1.7 [0.8,2.7] & 4.1 [2.8,5.5] & -30.2 [-33.4,-27.0] \\
Tukey IRLS & 0.9 [0.4,1.5] & 1.7 [0.8,2.6] & 3.9 [2.6,5.3] & -22.5 [-25.4,-19.7] \\
\midrule
\multicolumn{5}{l}{Mixture, \(N=5\); \(d=1,2,3,4\)} \\
GMM, \(K=2\) & 44.4 [41.0,47.8] & 30.4 [27.0,33.7] & 8.7 [6.0,11.4] & 4.2 [2.2,6.2] \\
GMM, \(K=3\) & 44.4 [41.0,47.7] & 30.4 [27.2,33.7] & 8.7 [6.0,11.4] & 4.2 [2.2,6.2] \\
GMM, \(K=4\) & 44.4 [41.1,47.7] & 30.4 [27.1,33.8] & 8.7 [6.0,11.4] & 4.2 [2.2,6.2] \\
Huber IRLS & 47.5 [44.2,50.8] & 34.3 [31.2,37.5] & 13.7 [11.0,16.5] & 6.6 [4.6,8.6] \\
Cauchy IRLS & 47.3 [44.0,50.5] & 34.0 [30.8,37.2] & 12.2 [9.5,14.9] & 5.5 [3.5,7.5] \\
IGG~III IRLS & 48.2 [44.8,51.5] & 37.0 [33.7,40.3] & 12.8 [9.9,15.7] & 5.9 [3.8,8.0] \\
Tukey IRLS & 48.4 [45.1,51.7] & 36.3 [33.0,39.6] & 13.9 [11.0,16.8] & 6.0 [4.0,8.0] \\
\midrule
\multicolumn{5}{l}{Mixture, \(N=10\); \(d=2,3,4,6\)} \\
GMM, \(K=2\) & 43.0 [39.9,46.1] & 45.5 [42.2,48.7] & 41.2 [37.7,44.6] & 11.8 [8.5,15.1] \\
GMM, \(K=3\) & 43.0 [39.8,46.1] & 45.5 [42.2,48.7] & 41.2 [37.8,44.7] & 11.8 [8.6,15.1] \\
GMM, \(K=4\) & 43.0 [39.9,46.1] & 45.5 [42.2,48.8] & 41.2 [37.8,44.6] & 11.8 [8.6,15.0] \\
Huber IRLS & 43.4 [40.3,46.5] & 45.9 [42.6,49.2] & 41.7 [38.4,45.0] & 17.3 [14.2,20.4] \\
Cauchy IRLS & 43.3 [40.2,46.4] & 46.0 [42.8,49.2] & 41.8 [38.5,45.2] & 14.6 [11.4,17.9] \\
IGG~III IRLS & 43.8 [40.7,46.9] & 46.4 [43.2,49.7] & 42.7 [39.2,46.1] & 14.8 [11.4,18.2] \\
Tukey IRLS & 43.9 [40.9,46.9] & 46.7 [43.5,49.9] & 42.5 [39.1,45.9] & 16.4 [13.0,19.7] \\
\midrule
\multicolumn{5}{l}{Cauchy, \(N=5\); \(d=1,2,3,4\)} \\
GMM, \(K=2\) & 55.1 [51.9,58.3] & 37.8 [34.6,41.0] & 11.4 [8.8,14.0] & 4.2 [2.4,6.1] \\
GMM, \(K=3\) & 55.3 [52.0,58.5] & 37.8 [34.6,40.9] & 11.4 [8.9,14.0] & 4.2 [2.4,6.1] \\
GMM, \(K=4\) & 55.4 [52.2,58.6] & 37.8 [34.7,41.0] & 11.4 [8.9,14.0] & 4.2 [2.4,6.1] \\
Huber IRLS & 58.8 [55.5,62.0] & 42.6 [39.5,45.8] & 16.7 [14.0,19.4] & 7.3 [5.3,9.3] \\
Cauchy IRLS & 61.4 [58.1,64.6] & 43.5 [40.3,46.7] & 17.0 [14.2,19.8] & 6.1 [4.2,8.1] \\
IGG~III IRLS & 62.6 [59.4,65.8] & 44.9 [41.6,48.2] & 16.9 [14.0,19.7] & 5.6 [3.7,7.6] \\
Tukey IRLS & 63.5 [60.3,66.6] & 45.3 [42.1,48.5] & 18.9 [16.1,21.7] & 6.0 [4.0,8.0] \\
\midrule
\multicolumn{5}{l}{Cauchy, \(N=10\); \(d=2,3,4,6\)} \\
GMM, \(K=2\) & 76.9 [74.3,79.5] & 72.5 [69.7,75.3] & 65.4 [62.4,68.4] & 28.5 [25.4,31.6] \\
GMM, \(K=3\) & 77.6 [74.9,80.1] & 73.0 [70.2,75.7] & 65.9 [62.9,68.9] & 29.0 [25.9,32.2] \\
GMM, \(K=4\) & 77.8 [75.2,80.4] & 73.4 [70.6,76.1] & 66.3 [63.3,69.3] & 29.0 [25.9,32.1] \\
Huber IRLS & 83.1 [80.8,85.4] & 79.5 [77.0,82.0] & 71.2 [68.3,74.0] & 37.0 [33.8,40.2] \\
Cauchy IRLS & 83.0 [80.6,85.2] & 79.9 [77.3,82.4] & 71.4 [68.5,74.3] & 34.8 [31.6,38.0] \\
IGG~III IRLS & 83.4 [81.0,85.7] & 79.9 [77.3,82.4] & 71.6 [68.6,74.4] & 35.4 [32.2,38.6] \\
Tukey IRLS & 83.4 [81.0,85.7] & 80.8 [78.3,83.3] & 71.6 [68.7,74.4] & 36.1 [32.9,39.3] \\
\bottomrule
\end{longtable}
\endgroup

\begin{figure}[!htbp]\centering
\includegraphics[width=.70\textwidth]{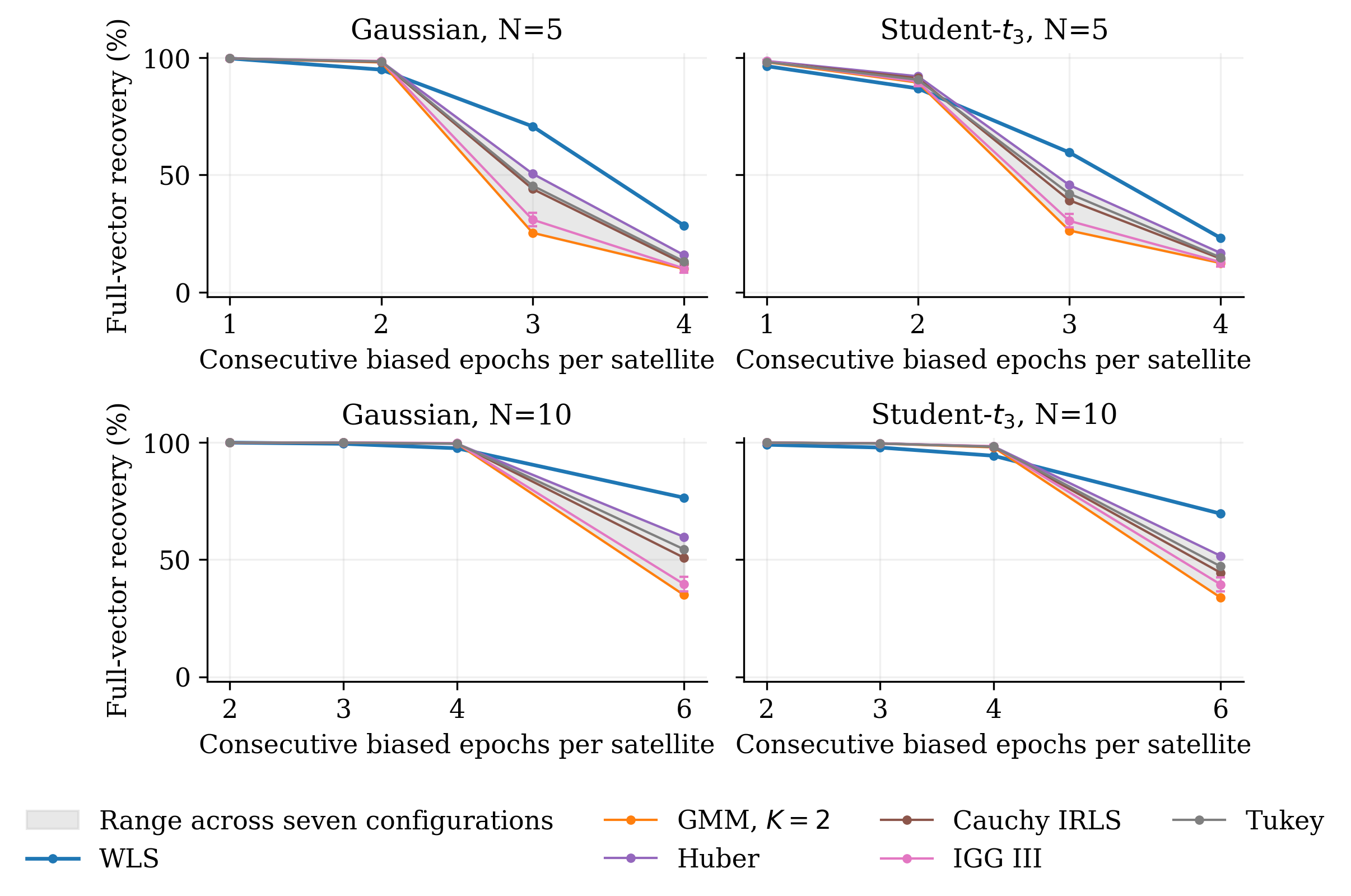}
\caption{Full-vector recovery in the duration controls. Columns show Gaussian and Student-$t_3$ backgrounds. Rows show $N=5$ and $N=10$. The horizontal axis gives the number of consecutive biased epochs per affected satellite. The shaded region spans the seven non-WLS weighting configurations. Curves show WLS, GMM with $K=2$, and the four IRLS rules. IGG~III error bars give pointwise Wilson intervals. WLS leads in several longer-bias settings.}
\label{fig:duration-controls}
\end{figure}
\FloatBarrier
\section{Simulation acceptance outcomes}
Within each error law, calibration and evaluation use disjoint seeds after repair of two legacy 32-bit Gaussian seed collisions. The later canonical keys receive unique 64-bit nonce seeds fixed without outcome information, and every method is recomputed for the repaired inputs. Tables S5 and S6 use all planned evaluation problems, under the Wilson-reference and nominal acceptance-target rules, respectively. The common-completed sensitivity uses its own denominator. The Wilson upper endpoint is computed from a two-sided nominal 95\% interval for the conditional wrong-fix frequency $q$. The 0.5\% reference applies to this conditional frequency. Table S28 presents the Gaussian and Student-$t_3$ results for the Wilson rule.
\begingroup
\footnotesize
\setlength{\tabcolsep}{2pt}
\begin{longtable}{lrrrrrr}
\caption{Evaluation of thresholds selected on calibration problems by the 0.5\% Wilson rule. A dash in $\tau$ means no qualifying threshold. At zero evaluation acceptance, $A=0$ and both $q$ and the upper endpoint are undefined.}\label{tab:integrity-unified-full}\\
\toprule
Estimator & Calibration $\tau$ & Accepted $n_a$ & Wrong $n_w$ & $A$ (\%) & $q$ (\%) & \shortstack{Wilson upper endpoint\\for $q$ (\%)} \\
\midrule
\endfirsthead
\multicolumn{7}{l}{\tablename~\thetable\ (continued)} \\
\toprule
Estimator & Calibration $\tau$ & Accepted $n_a$ & Wrong $n_w$ & $A$ (\%) & $q$ (\%) & \shortstack{Wilson upper endpoint\\for $q$ (\%)} \\
\midrule
\endhead
\multicolumn{7}{l}{Gaussian, all planned evaluation problems} \\
WLS & 2.7020 & 14422 & 39 & 55.5 & 0.270 & 0.369 \\
GMM, \(K=2\) & 10228251914.4876 & 1542 & 7 & 5.9 & 0.454 & 0.934 \\
GMM, \(K=3\) & -- & 0 & 0 & 0.0 & -- & -- \\
GMM, \(K=4\) & -- & 0 & 0 & 0.0 & -- & -- \\
Huber IRLS & 4.1638 & 14675 & 43 & 56.4 & 0.293 & 0.394 \\
Cauchy IRLS & 4.9591 & 14434 & 45 & 55.5 & 0.312 & 0.417 \\
IGG~III IRLS & 7.8585 & 10536 & 46 & 40.5 & 0.437 & 0.582 \\
Tukey IRLS & 5.5832 & 12851 & 40 & 49.4 & 0.311 & 0.424 \\
\multicolumn{7}{l}{Student-\(t_3\), all planned evaluation problems} \\
WLS & 4.6765 & 5006 & 19 & 19.3 & 0.380 & 0.592 \\
GMM, \(K=2\) & -- & 0 & 0 & 0.0 & -- & -- \\
GMM, \(K=3\) & -- & 0 & 0 & 0.0 & -- & -- \\
GMM, \(K=4\) & -- & 0 & 0 & 0.0 & -- & -- \\
Huber IRLS & 5.7215 & 9097 & 30 & 35.0 & 0.330 & 0.470 \\
Cauchy IRLS & 6.3567 & 9656 & 38 & 37.1 & 0.394 & 0.540 \\
IGG~III IRLS & 18.6436 & 2341 & 5 & 9.0 & 0.214 & 0.499 \\
Tukey IRLS & 9.3811 & 5911 & 15 & 22.7 & 0.254 & 0.418 \\
\multicolumn{7}{l}{Mixture, all planned evaluation problems} \\
WLS & -- & 0 & 0 & 0.0 & -- & -- \\
GMM, \(K=2\) & -- & 0 & 0 & 0.0 & -- & -- \\
GMM, \(K=3\) & -- & 0 & 0 & 0.0 & -- & -- \\
GMM, \(K=4\) & -- & 0 & 0 & 0.0 & -- & -- \\
Huber IRLS & -- & 0 & 0 & 0.0 & -- & -- \\
Cauchy IRLS & -- & 0 & 0 & 0.0 & -- & -- \\
IGG~III IRLS & -- & 0 & 0 & 0.0 & -- & -- \\
Tukey IRLS & -- & 0 & 0 & 0.0 & -- & -- \\
\multicolumn{7}{l}{Cauchy, all planned evaluation problems} \\
WLS & -- & 0 & 0 & 0.0 & -- & -- \\
GMM, \(K=2\) & -- & 0 & 0 & 0.0 & -- & -- \\
GMM, \(K=3\) & -- & 0 & 0 & 0.0 & -- & -- \\
GMM, \(K=4\) & -- & 0 & 0 & 0.0 & -- & -- \\
Huber IRLS & -- & 0 & 0 & 0.0 & -- & -- \\
Cauchy IRLS & -- & 0 & 0 & 0.0 & -- & -- \\
IGG~III IRLS & -- & 0 & 0 & 0.0 & -- & -- \\
Tukey IRLS & -- & 0 & 0 & 0.0 & -- & -- \\
\bottomrule
\end{longtable}
\endgroup
\begingroup
\footnotesize
\setlength{\tabcolsep}{2pt}
\begin{longtable}{lrrrrrr}
\caption{Actual evaluation outcomes by nominal calibration acceptance target. There are 26,000 planned evaluation problems per law. Correct and accepted (\%) is $100(n_a-n_w)/n$. A dash in $q$ indicates zero acceptances. The archive gives the common-completed sensitivity.}\label{tab:integrity-unified-quantiles}\\
\toprule
Estimator & Accepted $n_a$ & Wrong $n_w$ & $A$ (\%) & $q$ (\%) & $u$ (\%) & \shortstack{Correct and\\accepted (\%)} \\
\midrule
\endfirsthead
\multicolumn{7}{l}{\tablename~\thetable\ (continued)} \\
\toprule
Estimator & Accepted $n_a$ & Wrong $n_w$ & $A$ (\%) & $q$ (\%) & $u$ (\%) & \shortstack{Correct and\\accepted (\%)} \\
\midrule
\endhead
\multicolumn{7}{l}{Gaussian, nominal calibration acceptance 5\%} \\
WLS & 1285 & 0 & 4.9 & 0.000 & 0.000 & 4.9 \\
GMM, \(K=2\) & 1291 & 3 & 5.0 & 0.232 & 0.012 & 5.0 \\
GMM, \(K=3\) & 1316 & 23 & 5.1 & 1.748 & 0.088 & 5.0 \\
GMM, \(K=4\) & 1311 & 29 & 5.0 & 2.212 & 0.112 & 4.9 \\
Huber IRLS & 1294 & 0 & 5.0 & 0.000 & 0.000 & 5.0 \\
Cauchy IRLS & 1300 & 0 & 5.0 & 0.000 & 0.000 & 5.0 \\
IGG~III IRLS & 1342 & 0 & 5.2 & 0.000 & 0.000 & 5.2 \\
Tukey IRLS & 1309 & 0 & 5.0 & 0.000 & 0.000 & 5.0 \\
\multicolumn{7}{l}{Gaussian, nominal calibration acceptance 10\%} \\
WLS & 2577 & 0 & 9.9 & 0.000 & 0.000 & 9.9 \\
GMM, \(K=2\) & 2603 & 21 & 10.0 & 0.807 & 0.081 & 9.9 \\
GMM, \(K=3\) & 2561 & 65 & 9.8 & 2.538 & 0.250 & 9.6 \\
GMM, \(K=4\) & 2560 & 75 & 9.8 & 2.930 & 0.288 & 9.6 \\
Huber IRLS & 2599 & 0 & 10.0 & 0.000 & 0.000 & 10.0 \\
Cauchy IRLS & 2617 & 0 & 10.1 & 0.000 & 0.000 & 10.1 \\
IGG~III IRLS & 2523 & 2 & 9.7 & 0.079 & 0.008 & 9.7 \\
Tukey IRLS & 2650 & 1 & 10.2 & 0.038 & 0.004 & 10.2 \\
\multicolumn{7}{l}{Gaussian, nominal calibration acceptance 25\%} \\
WLS & 6432 & 0 & 24.7 & 0.000 & 0.000 & 24.7 \\
GMM, \(K=2\) & 6448 & 198 & 24.8 & 3.071 & 0.762 & 24.0 \\
GMM, \(K=3\) & 6427 & 351 & 24.7 & 5.461 & 1.350 & 23.4 \\
GMM, \(K=4\) & 6424 & 365 & 24.7 & 5.682 & 1.404 & 23.3 \\
Huber IRLS & 6454 & 1 & 24.8 & 0.015 & 0.004 & 24.8 \\
Cauchy IRLS & 6459 & 1 & 24.8 & 0.015 & 0.004 & 24.8 \\
IGG~III IRLS & 6442 & 9 & 24.8 & 0.140 & 0.035 & 24.7 \\
Tukey IRLS & 6379 & 4 & 24.5 & 0.063 & 0.015 & 24.5 \\
\multicolumn{7}{l}{Gaussian, nominal calibration acceptance 50\%} \\
WLS & 12933 & 17 & 49.7 & 0.131 & 0.065 & 49.7 \\
GMM, \(K=2\) & 13024 & 1158 & 50.1 & 8.891 & 4.454 & 45.6 \\
GMM, \(K=3\) & 12953 & 1293 & 49.8 & 9.982 & 4.973 & 44.8 \\
GMM, \(K=4\) & 13043 & 1330 & 50.2 & 10.197 & 5.115 & 45.0 \\
Huber IRLS & 12906 & 18 & 49.6 & 0.139 & 0.069 & 49.6 \\
Cauchy IRLS & 12893 & 21 & 49.6 & 0.163 & 0.081 & 49.5 \\
IGG~III IRLS & 12884 & 101 & 49.6 & 0.784 & 0.388 & 49.2 \\
Tukey IRLS & 12851 & 40 & 49.4 & 0.311 & 0.154 & 49.3 \\
\multicolumn{7}{l}{Student-\(t_3\), nominal calibration acceptance 5\%} \\
WLS & 1333 & 2 & 5.1 & 0.150 & 0.008 & 5.1 \\
GMM, \(K=2\) & 1277 & 10 & 4.9 & 0.783 & 0.038 & 4.9 \\
GMM, \(K=3\) & 1255 & 32 & 4.8 & 2.550 & 0.123 & 4.7 \\
GMM, \(K=4\) & 1302 & 44 & 5.0 & 3.379 & 0.169 & 4.8 \\
Huber IRLS & 1273 & 0 & 4.9 & 0.000 & 0.000 & 4.9 \\
Cauchy IRLS & 1268 & 0 & 4.9 & 0.000 & 0.000 & 4.9 \\
IGG~III IRLS & 1281 & 1 & 4.9 & 0.078 & 0.004 & 4.9 \\
Tukey IRLS & 1335 & 3 & 5.1 & 0.225 & 0.012 & 5.1 \\
\multicolumn{7}{l}{Student-\(t_3\), nominal calibration acceptance 10\%} \\
WLS & 2631 & 7 & 10.1 & 0.266 & 0.027 & 10.1 \\
GMM, \(K=2\) & 2470 & 57 & 9.5 & 2.308 & 0.219 & 9.3 \\
GMM, \(K=3\) & 2582 & 115 & 9.9 & 4.454 & 0.442 & 9.5 \\
GMM, \(K=4\) & 2574 & 155 & 9.9 & 6.022 & 0.596 & 9.3 \\
Huber IRLS & 2630 & 1 & 10.1 & 0.038 & 0.004 & 10.1 \\
Cauchy IRLS & 2614 & 1 & 10.1 & 0.038 & 0.004 & 10.1 \\
IGG~III IRLS & 2582 & 6 & 9.9 & 0.232 & 0.023 & 9.9 \\
Tukey IRLS & 2654 & 6 & 10.2 & 0.226 & 0.023 & 10.2 \\
\multicolumn{7}{l}{Student-\(t_3\), nominal calibration acceptance 25\%} \\
WLS & 6651 & 42 & 25.6 & 0.631 & 0.162 & 25.4 \\
GMM, \(K=2\) & 6468 & 493 & 24.9 & 7.622 & 1.896 & 23.0 \\
GMM, \(K=3\) & 6549 & 709 & 25.2 & 10.826 & 2.727 & 22.5 \\
GMM, \(K=4\) & 6488 & 786 & 25.0 & 12.115 & 3.023 & 21.9 \\
Huber IRLS & 6497 & 12 & 25.0 & 0.185 & 0.046 & 24.9 \\
Cauchy IRLS & 6486 & 10 & 24.9 & 0.154 & 0.038 & 24.9 \\
IGG~III IRLS & 6534 & 37 & 25.1 & 0.566 & 0.142 & 25.0 \\
Tukey IRLS & 6389 & 17 & 24.6 & 0.266 & 0.065 & 24.5 \\
\multicolumn{7}{l}{Student-\(t_3\), nominal calibration acceptance 50\%} \\
WLS & 12906 & 466 & 49.6 & 3.611 & 1.792 & 47.8 \\
GMM, \(K=2\) & 13060 & 1932 & 50.2 & 14.793 & 7.431 & 42.8 \\
GMM, \(K=3\) & 12938 & 2288 & 49.8 & 17.684 & 8.800 & 41.0 \\
GMM, \(K=4\) & 12906 & 2368 & 49.6 & 18.348 & 9.108 & 40.5 \\
Huber IRLS & 13014 & 132 & 50.1 & 1.014 & 0.508 & 49.5 \\
Cauchy IRLS & 12953 & 134 & 49.8 & 1.035 & 0.515 & 49.3 \\
IGG~III IRLS & 12890 & 260 & 49.6 & 2.017 & 1.000 & 48.6 \\
Tukey IRLS & 12954 & 178 & 49.8 & 1.374 & 0.685 & 49.1 \\
\multicolumn{7}{l}{Mixture, nominal calibration acceptance 5\%} \\
WLS & 1296 & 902 & 5.0 & 69.599 & 3.469 & 1.5 \\
GMM, \(K=2\) & 1254 & 146 & 4.8 & 11.643 & 0.562 & 4.3 \\
GMM, \(K=3\) & 1334 & 186 & 5.1 & 13.943 & 0.715 & 4.4 \\
GMM, \(K=4\) & 1247 & 228 & 4.8 & 18.284 & 0.877 & 3.9 \\
Huber IRLS & 1336 & 99 & 5.1 & 7.410 & 0.381 & 4.8 \\
Cauchy IRLS & 1343 & 70 & 5.2 & 5.212 & 0.269 & 4.9 \\
IGG~III IRLS & 1255 & 54 & 4.8 & 4.303 & 0.208 & 4.6 \\
Tukey IRLS & 1365 & 61 & 5.2 & 4.469 & 0.235 & 5.0 \\
\multicolumn{7}{l}{Mixture, nominal calibration acceptance 10\%} \\
WLS & 2601 & 1946 & 10.0 & 74.817 & 7.485 & 2.5 \\
GMM, \(K=2\) & 2488 & 484 & 9.6 & 19.453 & 1.862 & 7.7 \\
GMM, \(K=3\) & 2513 & 537 & 9.7 & 21.369 & 2.065 & 7.6 \\
GMM, \(K=4\) & 2476 & 620 & 9.5 & 25.040 & 2.385 & 7.1 \\
Huber IRLS & 2680 & 318 & 10.3 & 11.866 & 1.223 & 9.1 \\
Cauchy IRLS & 2665 & 243 & 10.2 & 9.118 & 0.935 & 9.3 \\
IGG~III IRLS & 2589 & 172 & 10.0 & 6.643 & 0.662 & 9.3 \\
Tukey IRLS & 2625 & 184 & 10.1 & 7.010 & 0.708 & 9.4 \\
\multicolumn{7}{l}{Mixture, nominal calibration acceptance 25\%} \\
WLS & 6548 & 5353 & 25.2 & 81.750 & 20.588 & 4.6 \\
GMM, \(K=2\) & 6415 & 2104 & 24.7 & 32.798 & 8.092 & 16.6 \\
GMM, \(K=3\) & 6333 & 2347 & 24.4 & 37.060 & 9.027 & 15.3 \\
GMM, \(K=4\) & 6381 & 2499 & 24.5 & 39.163 & 9.612 & 14.9 \\
Huber IRLS & 6522 & 1509 & 25.1 & 23.137 & 5.804 & 19.3 \\
Cauchy IRLS & 6562 & 1233 & 25.2 & 18.790 & 4.742 & 20.5 \\
IGG~III IRLS & 6441 & 897 & 24.8 & 13.926 & 3.450 & 21.3 \\
Tukey IRLS & 6527 & 1054 & 25.1 & 16.148 & 4.054 & 21.1 \\
\multicolumn{7}{l}{Mixture, nominal calibration acceptance 50\%} \\
WLS & 13063 & 11259 & 50.2 & 86.190 & 43.304 & 6.9 \\
GMM, \(K=2\) & 12935 & 6426 & 49.8 & 49.679 & 24.715 & 25.0 \\
GMM, \(K=3\) & 12931 & 6136 & 49.7 & 47.452 & 23.600 & 26.1 \\
GMM, \(K=4\) & 12896 & 6329 & 49.6 & 49.077 & 24.342 & 25.3 \\
Huber IRLS & 12982 & 5120 & 49.9 & 39.439 & 19.692 & 30.2 \\
Cauchy IRLS & 13009 & 4637 & 50.0 & 35.645 & 17.835 & 32.2 \\
IGG~III IRLS & 12969 & 3932 & 49.9 & 30.318 & 15.123 & 34.8 \\
Tukey IRLS & 12919 & 4216 & 49.7 & 32.634 & 16.215 & 33.5 \\
\multicolumn{7}{l}{Cauchy, nominal calibration acceptance 5\%} \\
WLS & 1340 & 906 & 5.2 & 67.612 & 3.485 & 1.7 \\
GMM, \(K=2\) & 1285 & 129 & 4.9 & 10.039 & 0.496 & 4.4 \\
GMM, \(K=3\) & 1294 & 164 & 5.0 & 12.674 & 0.631 & 4.3 \\
GMM, \(K=4\) & 1295 & 216 & 5.0 & 16.680 & 0.831 & 4.2 \\
Huber IRLS & 1302 & 32 & 5.0 & 2.458 & 0.123 & 4.9 \\
Cauchy IRLS & 1261 & 32 & 4.8 & 2.538 & 0.123 & 4.7 \\
IGG~III IRLS & 1252 & 33 & 4.8 & 2.636 & 0.127 & 4.7 \\
Tukey IRLS & 1217 & 28 & 4.7 & 2.301 & 0.108 & 4.6 \\
\multicolumn{7}{l}{Cauchy, nominal calibration acceptance 10\%} \\
WLS & 2651 & 1971 & 10.2 & 74.349 & 7.581 & 2.6 \\
GMM, \(K=2\) & 2587 & 423 & 9.9 & 16.351 & 1.627 & 8.3 \\
GMM, \(K=3\) & 2600 & 507 & 10.0 & 19.500 & 1.950 & 8.1 \\
GMM, \(K=4\) & 2638 & 616 & 10.1 & 23.351 & 2.369 & 7.8 \\
Huber IRLS & 2587 & 138 & 9.9 & 5.334 & 0.531 & 9.4 \\
Cauchy IRLS & 2587 & 94 & 9.9 & 3.634 & 0.362 & 9.6 \\
IGG~III IRLS & 2536 & 93 & 9.8 & 3.667 & 0.358 & 9.4 \\
Tukey IRLS & 2535 & 78 & 9.8 & 3.077 & 0.300 & 9.4 \\
\multicolumn{7}{l}{Cauchy, nominal calibration acceptance 25\%} \\
WLS & 6540 & 5428 & 25.2 & 82.997 & 20.877 & 4.3 \\
GMM, \(K=2\) & 6567 & 1955 & 25.3 & 29.770 & 7.519 & 17.7 \\
GMM, \(K=3\) & 6532 & 2205 & 25.1 & 33.757 & 8.481 & 16.6 \\
GMM, \(K=4\) & 6398 & 2281 & 24.6 & 35.652 & 8.773 & 15.8 \\
Huber IRLS & 6534 & 913 & 25.1 & 13.973 & 3.512 & 21.6 \\
Cauchy IRLS & 6534 & 603 & 25.1 & 9.229 & 2.319 & 22.8 \\
IGG~III IRLS & 6428 & 520 & 24.7 & 8.090 & 2.000 & 22.7 \\
Tukey IRLS & 6478 & 535 & 24.9 & 8.259 & 2.058 & 22.9 \\
\multicolumn{7}{l}{Cauchy, nominal calibration acceptance 50\%} \\
WLS & 13028 & 11497 & 50.1 & 88.248 & 44.219 & 5.9 \\
GMM, \(K=2\) & 13083 & 5896 & 50.3 & 45.066 & 22.677 & 27.6 \\
GMM, \(K=3\) & 12952 & 5717 & 49.8 & 44.140 & 21.988 & 27.8 \\
GMM, \(K=4\) & 12832 & 5822 & 49.4 & 45.371 & 22.392 & 27.0 \\
Huber IRLS & 13152 & 4017 & 50.6 & 30.543 & 15.450 & 35.1 \\
Cauchy IRLS & 13157 & 3177 & 50.6 & 24.147 & 12.219 & 38.4 \\
IGG~III IRLS & 12978 & 2757 & 49.9 & 21.244 & 10.604 & 39.3 \\
Tukey IRLS & 13077 & 2908 & 50.3 & 22.238 & 11.185 & 39.1 \\
\bottomrule
\end{longtable}
\endgroup

\begin{table}[htbp]
\centering
\caption{Design features of the two acceptance evaluations. The full simulation strata and cross-date policies are available as electronic tables.}
\begin{tabular}{@{}p{.17\linewidth}p{.39\linewidth}p{.36\linewidth}@{}}
\toprule
Feature & Between-satellite simulation & Sparse measured-observation case \\
\midrule
Geometry and label & Favourable and weak simulated geometries, five epochs, generated true integer & YARR--YAR3 short baseline, five and ten epochs, internal screened long-arc integer \\
Contamination & 2,000 clean and 24,000 contaminated problems per law in each of calibration and evaluation. Local biases \(1\sigma,3\sigma,5\sigma\) on two satellites, one or two epochs & Measured background plus random whitened innovations. The added Gaussian component with standard deviation 10 on the whitened scale is drawn with probability 0.04 per coordinate. Cauchy additions use MAD scaling \\
Calibration & 26,000 problems per law in each set, with disjoint calibration and evaluation seeds within law & 246--535 planned windows per date, length and law. Fit on one date and evaluate on the other under the shared working model \\
Mixture recovery & WLS 10.0\%, Huber IRLS 40.8\%, IGG~III IRLS 44.3\% in the pooled evaluation & WLS 8.5--12.1\% and Huber IRLS 85.5--97.6\% across the four mixture cells \\
\bottomrule
\end{tabular}
\end{table}

\begin{longtable}{llrrr}
\caption{IGG~III acceptance under Gaussian-mixture errors at the nominal 10\% target. At each bias level, each geometry has 4,000 planned problems and the pooled row has 8,000. The same frozen threshold applies to every row. At $5\sigma$, 29 of 44 accepted vectors are wrong in weak geometry.}\label{tab:igg-geometry}\\
\toprule
Bias & Geometry & Accepted & Wrong & $q$ among accepted (\%) \\
\midrule\endfirsthead
\toprule
Bias & Geometry & Accepted & Wrong & $q$ among accepted (\%) \\
\midrule\endhead
1$\sigma$ & Favourable & 978 & 26 & 2.66 \\
1$\sigma$ & Weak & 58 & 21 & 36.21 \\
1$\sigma$ & Both geometries & 1036 & 47 & 4.54 \\
5$\sigma$ & Favourable & 551 & 39 & 7.08 \\
5$\sigma$ & Weak & 44 & 29 & 65.91 \\
5$\sigma$ & Both geometries & 595 & 68 & 11.43 \\
\bottomrule
\end{longtable}

\FloatBarrier
\section{Short-record weighting diagnostics}
Forty-eight cases were selected to include extreme and ordinary ratios and correct and incorrect fixes. The 156 original-setting replays comprise 144 GMM fits on 48 inputs at three values of $K$, plus 12 WLS and IGG~III controls. All reproduce the stored integer vectors, with maximum baseline-coordinate difference below $10^{-12}\,\mathrm m$. In the GMM jobs, $J_1$ ranges from $6.59\times10^{-6}$ to 3.75, above the denominator guard. The independent QR calculation agrees with the original ratio to about $1.5\times10^{-6}$ relative. Ridge-supported directions account for a median 99.87\% of $J_2$. WLS initialization changes four of six selected vectors in the parameter checks. The archive retains the auxiliary residual-objective ratio, component spectra, responsibilities, final arrays and iteration traces. The frozen policies use the ambiguity-distance ratio.
\begin{longtable}{lllrlllll}
\caption{Purposively selected GMM replays. Each law--length--configuration row has eight inputs, four per date. Ridge share is the median percentage of $J_2$ attributable to directions with effectively zero scatter eigenvalues in the saved working metric. The last column compares the original $R$ with a QR-based recomputation of the candidate distances.}\\
\toprule
Law & $N$ & Configuration & $n$ & Min $J_1$ & Med. $R$ & Ridge share (\%) & \shortstack{Count with\\$J_1<10^{-15}$} & \shortstack{Max relative\\difference in $R$ (QR)} \\
\midrule
\endfirsthead
\multicolumn{9}{l}{\tablename~\thetable\ (continued)} \\
\toprule
Law & $N$ & Configuration & $n$ & Min $J_1$ & Med. $R$ & Ridge share (\%) & \shortstack{Count with\\$J_1<10^{-15}$} & \shortstack{Max relative\\difference in $R$ (QR)} \\
\midrule
\endhead
Cauchy & 5 & GMM, $K=2$ & 8 & 6.59e-06 & 4.22e+08 & 99.91 & 0 & 6.0e-07 \\
Cauchy & 5 & GMM, $K=3$ & 8 & 7.04e-06 & 1.26e+08 & 99.21 & 0 & 5.9e-07 \\
Cauchy & 5 & GMM, $K=4$ & 8 & 6.85e-06 & 6.32e+07 & 99.47 & 0 & 1.2e-06 \\
Cauchy & 10 & GMM, $K=2$ & 8 & 5.09e-04 & 2.48e+07 & 92.59 & 0 & 1.1e-07 \\
Cauchy & 10 & GMM, $K=3$ & 8 & 3.78e-04 & 1.29e+08 & 92.59 & 0 & 1.7e-07 \\
Cauchy & 10 & GMM, $K=4$ & 8 & 3.78e-04 & 1.31e+08 & 92.53 & 0 & 2.3e-07 \\
No added error & 5 & GMM, $K=2$ & 8 & 2.40e-05 & 1.46e+08 & 99.98 & 0 & 1.4e-06 \\
No added error & 5 & GMM, $K=3$ & 8 & 2.40e-05 & 1.46e+08 & 99.98 & 0 & 2.8e-07 \\
No added error & 5 & GMM, $K=4$ & 8 & 2.40e-05 & 1.46e+08 & 99.98 & 0 & 2.9e-07 \\
No added error & 10 & GMM, $K=2$ & 8 & 1.35e-04 & 4.00e+07 & 90.64 & 0 & 9.3e-07 \\
No added error & 10 & GMM, $K=3$ & 8 & 1.35e-04 & 4.53e+07 & 76.23 & 0 & 1.5e-06 \\
No added error & 10 & GMM, $K=4$ & 8 & 1.35e-04 & 4.64e+07 & 76.05 & 0 & 1.2e-06 \\
Mixture & 5 & GMM, $K=2$ & 8 & 1.09e-05 & 1.18e+09 & 100.00 & 0 & 1.3e-06 \\
Mixture & 5 & GMM, $K=3$ & 8 & 1.09e-05 & 1.21e+09 & 99.99 & 0 & 5.1e-07 \\
Mixture & 5 & GMM, $K=4$ & 8 & 1.11e-05 & 1.26e+09 & 99.98 & 0 & 1.1e-06 \\
Mixture & 10 & GMM, $K=2$ & 8 & 9.17e-05 & 4.81e+07 & 99.90 & 0 & 1.3e-07 \\
Mixture & 10 & GMM, $K=3$ & 8 & 1.26e-04 & 1.07e+08 & 99.84 & 0 & 1.8e-07 \\
Mixture & 10 & GMM, $K=4$ & 8 & 1.22e-04 & 1.12e+08 & 99.93 & 0 & 7.5e-08 \\
\bottomrule
\end{longtable}

\begin{longtable}{lrrrr}
\caption{Sensitivity of GMM ($K=2$) on six high-ratio inputs selected on 31 March, one per law--length combination. Vector changes are relative to the original median-initialized configuration of main-text Section 2.3. WLS initialization changes four of six vectors. Both tested ridge changes retain the six vectors and shift the range of $R$.}\\
\toprule
Variant & $n$ & Vectors changed & Min $R$ & Max $R$ \\
\midrule
\endfirsthead
\toprule
Variant & $n$ & Vectors changed & Min $R$ & Max $R$ \\
\midrule
\endhead
Initial covariance $\times0.1$ & 6 & 0 & 9.68e+06 & 4.30e+10 \\
Initial covariance $\times10$ & 6 & 0 & 7.92e+05 & 4.30e+10 \\
WLS initialization & 6 & 4 & 5.77e+06 & 4.33e+10 \\
Ridge coefficient $10^{-2}$ & 6 & 0 & 4.04e+04 & 8.86e+06 \\
Ridge coefficient $10^{-6}$ & 6 & 0 & 3.03e+09 & 4.30e+14 \\
\bottomrule
\end{longtable}

\begingroup
\renewcommand{\arraystretch}{0.92}
\setlength{\LTcapwidth}{\linewidth}
\begin{longtable}{llllrrr}
\caption{48-case weight audit, with eight cases per input law and length. Tiny (\%) is the mean of the within-case percentages of final working weights at or below $10^{-11}$. Floor cases (\%) is the percentage of cases with at least one row-augmentation activation. Vectors changed counts complete vectors that differ when the augmentation response is set to zero instead of the weighted mean.}\\
\toprule
Law & $N$ & Method & $n$ & Tiny (\%) & Floor cases (\%) & \shortstack{Vectors\\changed (count)} \\
\midrule
\endfirsthead
\multicolumn{7}{l}{\tablename~\thetable\ (continued)} \\
\toprule
Law & $N$ & Method & $n$ & Tiny (\%) & Floor cases (\%) & \shortstack{Vectors\\changed (count)} \\
\midrule
\endhead
Cauchy & 5 & Cauchy IRLS & 8 & 0.0 & 0.0 & 0 \\
Cauchy & 5 & Huber IRLS & 8 & 0.0 & 0.0 & 0 \\
Cauchy & 5 & IGG~III IRLS & 8 & 15.5 & 12.5 & 0 \\
Cauchy & 5 & Tukey IRLS & 8 & 10.0 & 0.0 & 0 \\
Cauchy & 10 & Cauchy IRLS & 8 & 0.0 & 0.0 & 0 \\
Cauchy & 10 & Huber IRLS & 8 & 0.0 & 0.0 & 0 \\
Cauchy & 10 & IGG~III IRLS & 8 & 17.9 & 12.5 & 1 \\
Cauchy & 10 & Tukey IRLS & 8 & 10.5 & 12.5 & 1 \\
No added error & 5 & Cauchy IRLS & 8 & 0.0 & 0.0 & 0 \\
No added error & 5 & Huber IRLS & 8 & 0.0 & 0.0 & 0 \\
No added error & 5 & IGG~III IRLS & 8 & 11.5 & 12.5 & 0 \\
No added error & 5 & Tukey IRLS & 8 & 5.5 & 0.0 & 0 \\
No added error & 10 & Cauchy IRLS & 8 & 0.0 & 0.0 & 0 \\
No added error & 10 & Huber IRLS & 8 & 0.0 & 0.0 & 0 \\
No added error & 10 & IGG~III IRLS & 8 & 13.2 & 0.0 & 0 \\
No added error & 10 & Tukey IRLS & 8 & 5.3 & 0.0 & 0 \\
Mixture & 5 & Cauchy IRLS & 8 & 0.0 & 0.0 & 0 \\
Mixture & 5 & Huber IRLS & 8 & 0.0 & 0.0 & 0 \\
Mixture & 5 & IGG~III IRLS & 8 & 18.8 & 37.5 & 3 \\
Mixture & 5 & Tukey IRLS & 8 & 14.6 & 25.0 & 2 \\
Mixture & 10 & Cauchy IRLS & 8 & 0.0 & 0.0 & 0 \\
Mixture & 10 & Huber IRLS & 8 & 0.0 & 0.0 & 0 \\
Mixture & 10 & IGG~III IRLS & 8 & 19.0 & 12.5 & 1 \\
Mixture & 10 & Tukey IRLS & 8 & 12.8 & 0.0 & 0 \\
\bottomrule
\end{longtable}
\endgroup

\FloatBarrier
\section{Sparse-injection case study on measured observations}
Full-vector recovery is evaluated against internal long-arc reference integers, and baseline error is measured relative to the RINEX header-coordinate benchmark. Windows are selected without overlap within each date and record length. The record lengths and injection laws reuse measured observations. Cross-date evaluation transfers the threshold fitted on the other date under the shared working covariance model.
\begin{table}[htbp]\centering
\caption{Long-arc and window screening by date. Each date has 2,880 usable epochs. Each sliding start is counted under its first failing window rule. Shared windows at $N=5$ or $10$ each share at least one epoch with a window of the other length. Shared epochs count distinct epochs selected at both lengths.}
\begin{tabular}{lrrrrrr}\toprule Date & $N$ & \shortstack{Candidate\\starts} & Sat. change & Ref. change & \shortstack{Eligible\\windows} & \shortstack{Selected non-\\overlapping windows}\\\midrule
2023-03-31 & 5 & 2876 & 328 & 70 & 2478 & 535\\
2023-03-31 & 10 & 2871 & 689 & 121 & 2061 & 247\\
2023-04-30 & 5 & 2876 & 320 & 62 & 2494 & 532\\
2023-04-30 & 10 & 2871 & 681 & 116 & 2074 & 246\\
\bottomrule\end{tabular}
\par\medskip
\begin{tabular}{lrrrrrr}\toprule Date & Arcs & Accepted & Rejected & \shortstack{Satellite pairs\\with accepted arcs} & Shared epochs & \shortstack{Shared windows\\$N=5/10$}\\\midrule
2023-03-31 & 742 & 706 & 36 & 614 & 2470 & 494/247\\
2023-04-30 & 772 & 724 & 48 & 620 & 2460 & 492/246\\
\bottomrule\end{tabular}
\end{table}
\begin{table}[htbp]\centering
\caption{Long-arc screening failures by criterion and date. Each entry counts arcs triggering the listed criterion. The same arc can trigger several criteria.}
\begin{tabular}{lrr}\toprule Criterion & 03-31 & 04-30\\\midrule
At least 30 epochs & 32 & 44\\
Duration at least 900 s & 32 & 44\\
MAD at most 0.1 cycle & 4 & 4\\
At least 95\% within 0.25 cycle of median & 0 & 0\\
Median within 0.25 cycle of candidate integer & 0 & 0\\
Both half-medians select and lie within 0.25 cycle of candidate & 26 & 38\\
\bottomrule\end{tabular}\end{table}

\FloatBarrier
\begin{longtable}{llllrrr}
\caption{Baseline error before ratio screening relative to the RINEX header-coordinate benchmark, by date, length and injection law, for the eight weighting configurations. $F_{10}$ is the fraction of all planned windows in a cell with a completed solution within 10\,cm. Main-text Table~6 pools individual windows across dates and reports errors in centimetres.}\\
\toprule
Date & $N$ & Law & Method & Median (m) & P95 (m) & $F_{10}$ (\%) \\
\midrule
\endfirsthead
\multicolumn{7}{l}{\tablename~\thetable\ (continued)} \\
\toprule
Date & $N$ & Law & Method & Median (m) & P95 (m) & $F_{10}$ (\%) \\
\midrule
\endhead
03-31 & 5 & Cauchy & GMM, $K=2$ & 0.066 & 0.952 & 55.1 \\
03-31 & 5 & Cauchy & GMM, $K=3$ & 0.066 & 0.952 & 55.1 \\
03-31 & 5 & Cauchy & GMM, $K=4$ & 0.066 & 0.952 & 55.1 \\
03-31 & 5 & Cauchy & WLS & 0.811 & 7.361 & 20.7 \\
03-31 & 5 & Cauchy & Cauchy IRLS & 0.045 & 0.548 & 91.8 \\
03-31 & 5 & Cauchy & Huber IRLS & 0.046 & 0.580 & 90.8 \\
03-31 & 5 & Cauchy & IGG~III IRLS & 0.046 & 0.728 & 82.6 \\
03-31 & 5 & Cauchy & Tukey IRLS & 0.045 & 0.577 & 90.8 \\
03-31 & 5 & Mixture & GMM, $K=2$ & 0.107 & 1.033 & 49.9 \\
03-31 & 5 & Mixture & GMM, $K=3$ & 0.107 & 1.033 & 49.9 \\
03-31 & 5 & Mixture & GMM, $K=4$ & 0.107 & 1.033 & 49.9 \\
03-31 & 5 & Mixture & WLS & 0.930 & 2.971 & 14.6 \\
03-31 & 5 & Mixture & Cauchy IRLS & 0.045 & 0.501 & 91.2 \\
03-31 & 5 & Mixture & Huber IRLS & 0.045 & 0.523 & 90.8 \\
03-31 & 5 & Mixture & IGG~III IRLS & 0.046 & 0.757 & 80.9 \\
03-31 & 5 & Mixture & Tukey IRLS & 0.046 & 0.653 & 87.1 \\
03-31 & 5 & Student-$t$ & GMM, $K=2$ & 0.060 & 0.996 & 60.7 \\
03-31 & 5 & Student-$t$ & GMM, $K=3$ & 0.060 & 0.996 & 60.7 \\
03-31 & 5 & Student-$t$ & GMM, $K=4$ & 0.060 & 0.996 & 60.7 \\
03-31 & 5 & Student-$t$ & WLS & 0.068 & 1.194 & 56.3 \\
03-31 & 5 & Student-$t$ & Cauchy IRLS & 0.045 & 0.575 & 90.8 \\
03-31 & 5 & Student-$t$ & Huber IRLS & 0.045 & 0.574 & 91.0 \\
03-31 & 5 & Student-$t$ & IGG~III IRLS & 0.046 & 0.731 & 83.9 \\
03-31 & 5 & Student-$t$ & Tukey IRLS & 0.045 & 0.576 & 88.8 \\
03-31 & 10 & Cauchy & GMM, $K=2$ & 0.065 & 0.755 & 58.7 \\
03-31 & 10 & Cauchy & GMM, $K=3$ & 0.065 & 0.755 & 58.7 \\
03-31 & 10 & Cauchy & GMM, $K=4$ & 0.065 & 0.755 & 58.7 \\
03-31 & 10 & Cauchy & WLS & 0.805 & 5.634 & 19.8 \\
03-31 & 10 & Cauchy & Cauchy IRLS & 0.044 & 0.059 & 99.2 \\
03-31 & 10 & Cauchy & Huber IRLS & 0.044 & 0.060 & 99.2 \\
03-31 & 10 & Cauchy & IGG~III IRLS & 0.044 & 0.307 & 91.9 \\
03-31 & 10 & Cauchy & Tukey IRLS & 0.043 & 0.059 & 97.2 \\
03-31 & 10 & Mixture & GMM, $K=2$ & 0.265 & 0.955 & 43.7 \\
03-31 & 10 & Mixture & GMM, $K=3$ & 0.265 & 0.955 & 43.7 \\
03-31 & 10 & Mixture & GMM, $K=4$ & 0.265 & 0.955 & 43.7 \\
03-31 & 10 & Mixture & WLS & 0.817 & 2.043 & 9.7 \\
03-31 & 10 & Mixture & Cauchy IRLS & 0.044 & 0.059 & 96.4 \\
03-31 & 10 & Mixture & Huber IRLS & 0.045 & 0.063 & 95.5 \\
03-31 & 10 & Mixture & IGG~III IRLS & 0.044 & 0.280 & 94.7 \\
03-31 & 10 & Mixture & Tukey IRLS & 0.045 & 0.059 & 96.0 \\
03-31 & 10 & Student-$t$ & GMM, $K=2$ & 0.056 & 0.844 & 70.4 \\
03-31 & 10 & Student-$t$ & GMM, $K=3$ & 0.056 & 0.844 & 70.4 \\
03-31 & 10 & Student-$t$ & GMM, $K=4$ & 0.056 & 0.844 & 70.4 \\
03-31 & 10 & Student-$t$ & WLS & 0.050 & 0.738 & 81.0 \\
03-31 & 10 & Student-$t$ & Cauchy IRLS & 0.044 & 0.058 & 99.2 \\
03-31 & 10 & Student-$t$ & Huber IRLS & 0.044 & 0.058 & 100.0 \\
03-31 & 10 & Student-$t$ & IGG~III IRLS & 0.043 & 0.058 & 96.8 \\
03-31 & 10 & Student-$t$ & Tukey IRLS & 0.044 & 0.057 & 98.0 \\
04-30 & 5 & Cauchy & GMM, $K=2$ & 0.077 & 0.977 & 53.0 \\
04-30 & 5 & Cauchy & GMM, $K=3$ & 0.077 & 0.977 & 53.2 \\
04-30 & 5 & Cauchy & GMM, $K=4$ & 0.076 & 0.977 & 53.2 \\
04-30 & 5 & Cauchy & WLS & 0.753 & 7.645 & 23.7 \\
04-30 & 5 & Cauchy & Cauchy IRLS & 0.043 & 0.511 & 91.9 \\
04-30 & 5 & Cauchy & Huber IRLS & 0.043 & 0.580 & 91.0 \\
04-30 & 5 & Cauchy & IGG~III IRLS & 0.045 & 0.683 & 81.2 \\
04-30 & 5 & Cauchy & Tukey IRLS & 0.044 & 0.647 & 87.6 \\
04-30 & 5 & Mixture & GMM, $K=2$ & 0.219 & 1.022 & 47.2 \\
04-30 & 5 & Mixture & GMM, $K=3$ & 0.219 & 1.022 & 47.2 \\
04-30 & 5 & Mixture & GMM, $K=4$ & 0.219 & 1.022 & 47.2 \\
04-30 & 5 & Mixture & WLS & 0.971 & 2.815 & 15.2 \\
04-30 & 5 & Mixture & Cauchy IRLS & 0.044 & 0.642 & 87.8 \\
04-30 & 5 & Mixture & Huber IRLS & 0.045 & 0.618 & 88.2 \\
04-30 & 5 & Mixture & IGG~III IRLS & 0.045 & 0.686 & 78.8 \\
04-30 & 5 & Mixture & Tukey IRLS & 0.044 & 0.642 & 84.6 \\
04-30 & 5 & Student-$t$ & GMM, $K=2$ & 0.061 & 0.932 & 61.1 \\
04-30 & 5 & Student-$t$ & GMM, $K=3$ & 0.061 & 0.932 & 61.1 \\
04-30 & 5 & Student-$t$ & GMM, $K=4$ & 0.061 & 0.932 & 61.1 \\
04-30 & 5 & Student-$t$ & WLS & 0.076 & 1.243 & 54.3 \\
04-30 & 5 & Student-$t$ & Cauchy IRLS & 0.044 & 0.589 & 89.5 \\
04-30 & 5 & Student-$t$ & Huber IRLS & 0.045 & 0.579 & 90.4 \\
04-30 & 5 & Student-$t$ & IGG~III IRLS & 0.045 & 0.690 & 82.5 \\
04-30 & 5 & Student-$t$ & Tukey IRLS & 0.044 & 0.653 & 88.0 \\
04-30 & 10 & Cauchy & GMM, $K=2$ & 0.078 & 0.872 & 56.1 \\
04-30 & 10 & Cauchy & GMM, $K=3$ & 0.078 & 0.872 & 56.1 \\
04-30 & 10 & Cauchy & GMM, $K=4$ & 0.076 & 0.858 & 56.5 \\
04-30 & 10 & Cauchy & WLS & 0.706 & 5.049 & 19.9 \\
04-30 & 10 & Cauchy & Cauchy IRLS & 0.043 & 0.061 & 98.8 \\
04-30 & 10 & Cauchy & Huber IRLS & 0.043 & 0.065 & 98.0 \\
04-30 & 10 & Cauchy & IGG~III IRLS & 0.042 & 0.203 & 94.7 \\
04-30 & 10 & Cauchy & Tukey IRLS & 0.043 & 0.068 & 97.2 \\
04-30 & 10 & Mixture & GMM, $K=2$ & 0.146 & 0.821 & 49.2 \\
04-30 & 10 & Mixture & GMM, $K=3$ & 0.146 & 0.822 & 49.2 \\
04-30 & 10 & Mixture & GMM, $K=4$ & 0.146 & 0.822 & 49.2 \\
04-30 & 10 & Mixture & WLS & 0.777 & 1.883 & 12.2 \\
04-30 & 10 & Mixture & Cauchy IRLS & 0.043 & 0.066 & 97.2 \\
04-30 & 10 & Mixture & Huber IRLS & 0.043 & 0.065 & 98.0 \\
04-30 & 10 & Mixture & IGG~III IRLS & 0.042 & 0.291 & 93.9 \\
04-30 & 10 & Mixture & Tukey IRLS & 0.042 & 0.074 & 95.9 \\
04-30 & 10 & Student-$t$ & GMM, $K=2$ & 0.060 & 0.788 & 67.5 \\
04-30 & 10 & Student-$t$ & GMM, $K=3$ & 0.060 & 0.788 & 67.5 \\
04-30 & 10 & Student-$t$ & GMM, $K=4$ & 0.060 & 0.788 & 67.5 \\
04-30 & 10 & Student-$t$ & WLS & 0.049 & 0.720 & 82.9 \\
04-30 & 10 & Student-$t$ & Cauchy IRLS & 0.042 & 0.062 & 99.2 \\
04-30 & 10 & Student-$t$ & Huber IRLS & 0.043 & 0.062 & 99.2 \\
04-30 & 10 & Student-$t$ & IGG~III IRLS & 0.043 & 0.066 & 97.2 \\
04-30 & 10 & Student-$t$ & Tukey IRLS & 0.043 & 0.068 & 98.0 \\
\bottomrule
\end{longtable}

\begin{longtable}{llllrrrr}
\caption{Date-specific paired changes in full-vector recovery from WLS (pp). Each contrast uses the same $n$ windows. The 24 non-overlapping one-hour blocks start at midnight. Windows are assigned by start epoch and crossing windows are retained. The 95\% intervals use 2,000 paired block resamples. Other block sizes and boundary rules are in the archive.}\\
\toprule
Date & $N$ & Law & Method & $n$ & \shortstack{Change vs\\WLS (pp)} & \shortstack{95\% block-bootstrap\\interval} & Blocks \\
\midrule
\endfirsthead
\multicolumn{8}{l}{\tablename~\thetable\ (continued)} \\
\toprule
Date & $N$ & Law & Method & $n$ & \shortstack{Change vs\\WLS (pp)} & \shortstack{95\% block-bootstrap\\interval} & Blocks \\
\midrule
\endhead
03-31 & 5 & Mixture & GMM, $K=2$ & 535 & 31.4 & [27.3,35.4] & 24 \\
03-31 & 5 & Mixture & GMM, $K=3$ & 535 & 31.4 & [27.3,35.4] & 24 \\
03-31 & 5 & Mixture & GMM, $K=4$ & 535 & 31.4 & [27.3,35.4] & 24 \\
03-31 & 5 & Mixture & Huber IRLS & 535 & 77.0 & [73.9,80.1] & 24 \\
03-31 & 5 & Mixture & Cauchy IRLS & 535 & 76.8 & [73.8,79.9] & 24 \\
03-31 & 5 & Mixture & IGG~III IRLS & 535 & 59.4 & [55.5,63.4] & 24 \\
03-31 & 5 & Mixture & Tukey IRLS & 535 & 70.1 & [66.9,73.5] & 24 \\
03-31 & 10 & Mixture & GMM, $K=2$ & 247 & 32.4 & [26.1,39.0] & 24 \\
03-31 & 10 & Mixture & GMM, $K=3$ & 247 & 32.4 & [26.1,39.0] & 24 \\
03-31 & 10 & Mixture & GMM, $K=4$ & 247 & 32.4 & [26.1,39.0] & 24 \\
03-31 & 10 & Mixture & Huber IRLS & 247 & 86.6 & [82.7,90.6] & 24 \\
03-31 & 10 & Mixture & Cauchy IRLS & 247 & 87.0 & [83.4,90.8] & 24 \\
03-31 & 10 & Mixture & IGG~III IRLS & 247 & 77.3 & [72.7,82.0] & 24 \\
03-31 & 10 & Mixture & Tukey IRLS & 247 & 81.4 & [76.4,86.1] & 24 \\
04-30 & 5 & Mixture & GMM, $K=2$ & 532 & 28.9 & [22.9,34.9] & 24 \\
04-30 & 5 & Mixture & GMM, $K=3$ & 532 & 28.9 & [22.9,34.9] & 24 \\
04-30 & 5 & Mixture & GMM, $K=4$ & 532 & 28.9 & [22.9,34.9] & 24 \\
04-30 & 5 & Mixture & Huber IRLS & 532 & 73.5 & [70.2,76.6] & 24 \\
04-30 & 5 & Mixture & Cauchy IRLS & 532 & 72.7 & [68.6,76.9] & 24 \\
04-30 & 5 & Mixture & IGG~III IRLS & 532 & 57.5 & [51.6,62.7] & 24 \\
04-30 & 5 & Mixture & Tukey IRLS & 532 & 65.0 & [59.2,70.4] & 24 \\
04-30 & 10 & Mixture & GMM, $K=2$ & 246 & 31.7 & [25.1,38.7] & 24 \\
04-30 & 10 & Mixture & GMM, $K=3$ & 246 & 31.7 & [25.1,38.7] & 24 \\
04-30 & 10 & Mixture & GMM, $K=4$ & 246 & 31.7 & [25.1,38.7] & 24 \\
04-30 & 10 & Mixture & Huber IRLS & 246 & 86.6 & [81.9,91.1] & 24 \\
04-30 & 10 & Mixture & Cauchy IRLS & 246 & 85.8 & [80.7,90.6] & 24 \\
04-30 & 10 & Mixture & IGG~III IRLS & 246 & 74.4 & [68.7,80.4] & 24 \\
04-30 & 10 & Mixture & Tukey IRLS & 246 & 78.5 & [73.4,83.0] & 24 \\
03-31 & 5 & Student-$t$ & GMM, $K=2$ & 535 & -0.4 & [-6.5,4.6] & 24 \\
03-31 & 5 & Student-$t$ & GMM, $K=3$ & 535 & -0.4 & [-6.5,4.6] & 24 \\
03-31 & 5 & Student-$t$ & GMM, $K=4$ & 535 & -0.4 & [-6.5,4.6] & 24 \\
03-31 & 5 & Student-$t$ & Huber IRLS & 535 & 34.4 & [29.6,39.1] & 24 \\
03-31 & 5 & Student-$t$ & Cauchy IRLS & 535 & 33.3 & [28.4,38.1] & 24 \\
03-31 & 5 & Student-$t$ & IGG~III IRLS & 535 & 22.1 & [16.1,27.4] & 24 \\
03-31 & 5 & Student-$t$ & Tukey IRLS & 535 & 28.2 & [22.0,34.2] & 24 \\
03-31 & 10 & Student-$t$ & GMM, $K=2$ & 247 & -11.3 & [-17.3,-5.8] & 24 \\
03-31 & 10 & Student-$t$ & GMM, $K=3$ & 247 & -11.3 & [-17.3,-5.8] & 24 \\
03-31 & 10 & Student-$t$ & GMM, $K=4$ & 247 & -11.3 & [-17.3,-5.8] & 24 \\
03-31 & 10 & Student-$t$ & Huber IRLS & 247 & 19.8 & [15.9,24.2] & 24 \\
03-31 & 10 & Student-$t$ & Cauchy IRLS & 247 & 18.6 & [15.0,22.4] & 24 \\
03-31 & 10 & Student-$t$ & IGG~III IRLS & 247 & 13.4 & [8.7,17.8] & 24 \\
03-31 & 10 & Student-$t$ & Tukey IRLS & 247 & 15.4 & [11.6,19.3] & 24 \\
04-30 & 5 & Student-$t$ & GMM, $K=2$ & 532 & 2.3 & [-4.2,8.4] & 24 \\
04-30 & 5 & Student-$t$ & GMM, $K=3$ & 532 & 2.3 & [-4.2,8.4] & 24 \\
04-30 & 5 & Student-$t$ & GMM, $K=4$ & 532 & 2.3 & [-4.2,8.4] & 24 \\
04-30 & 5 & Student-$t$ & Huber IRLS & 532 & 34.0 & [27.1,40.7] & 24 \\
04-30 & 5 & Student-$t$ & Cauchy IRLS & 532 & 32.0 & [24.1,39.5] & 24 \\
04-30 & 5 & Student-$t$ & IGG~III IRLS & 532 & 22.7 & [16.2,29.5] & 24 \\
04-30 & 5 & Student-$t$ & Tukey IRLS & 532 & 27.8 & [20.0,35.4] & 24 \\
04-30 & 10 & Student-$t$ & GMM, $K=2$ & 246 & -19.5 & [-27.3,-12.0] & 24 \\
04-30 & 10 & Student-$t$ & GMM, $K=3$ & 246 & -19.5 & [-27.3,-12.0] & 24 \\
04-30 & 10 & Student-$t$ & GMM, $K=4$ & 246 & -19.5 & [-27.3,-12.0] & 24 \\
04-30 & 10 & Student-$t$ & Huber IRLS & 246 & 16.3 & [10.7,21.7] & 24 \\
04-30 & 10 & Student-$t$ & Cauchy IRLS & 246 & 16.3 & [10.7,21.7] & 24 \\
04-30 & 10 & Student-$t$ & IGG~III IRLS & 246 & 11.0 & [5.4,16.4] & 24 \\
04-30 & 10 & Student-$t$ & Tukey IRLS & 246 & 13.0 & [7.0,18.7] & 24 \\
03-31 & 5 & Cauchy & GMM, $K=2$ & 535 & 32.7 & [27.9,37.8] & 24 \\
03-31 & 5 & Cauchy & GMM, $K=3$ & 535 & 32.7 & [27.9,37.8] & 24 \\
03-31 & 5 & Cauchy & GMM, $K=4$ & 535 & 32.7 & [27.9,37.8] & 24 \\
03-31 & 5 & Cauchy & Huber IRLS & 535 & 72.7 & [68.2,76.7] & 24 \\
03-31 & 5 & Cauchy & Cauchy IRLS & 535 & 74.0 & [69.4,78.3] & 24 \\
03-31 & 5 & Cauchy & IGG~III IRLS & 535 & 58.9 & [54.4,63.1] & 24 \\
03-31 & 5 & Cauchy & Tukey IRLS & 535 & 69.0 & [64.5,73.2] & 24 \\
03-31 & 10 & Cauchy & GMM, $K=2$ & 247 & 33.6 & [27.0,39.8] & 24 \\
03-31 & 10 & Cauchy & GMM, $K=3$ & 247 & 33.6 & [27.0,39.8] & 24 \\
03-31 & 10 & Cauchy & GMM, $K=4$ & 247 & 33.6 & [27.0,39.8] & 24 \\
03-31 & 10 & Cauchy & Huber IRLS & 247 & 81.4 & [75.6,86.3] & 24 \\
03-31 & 10 & Cauchy & Cauchy IRLS & 247 & 80.6 & [75.2,85.3] & 24 \\
03-31 & 10 & Cauchy & IGG~III IRLS & 247 & 68.0 & [59.3,76.3] & 24 \\
03-31 & 10 & Cauchy & Tukey IRLS & 247 & 76.1 & [69.9,82.0] & 24 \\
04-30 & 5 & Cauchy & GMM, $K=2$ & 532 & 25.8 & [20.0,31.2] & 24 \\
04-30 & 5 & Cauchy & GMM, $K=3$ & 532 & 25.8 & [20.0,31.2] & 24 \\
04-30 & 5 & Cauchy & GMM, $K=4$ & 532 & 25.8 & [20.0,31.2] & 24 \\
04-30 & 5 & Cauchy & Huber IRLS & 532 & 67.7 & [64.0,71.6] & 24 \\
04-30 & 5 & Cauchy & Cauchy IRLS & 532 & 68.0 & [63.8,72.6] & 24 \\
04-30 & 5 & Cauchy & IGG~III IRLS & 532 & 51.7 & [47.0,56.1] & 24 \\
04-30 & 5 & Cauchy & Tukey IRLS & 532 & 58.8 & [54.9,62.9] & 24 \\
04-30 & 10 & Cauchy & GMM, $K=2$ & 246 & 32.5 & [24.5,39.9] & 24 \\
04-30 & 10 & Cauchy & GMM, $K=3$ & 246 & 32.5 & [24.5,39.9] & 24 \\
04-30 & 10 & Cauchy & GMM, $K=4$ & 246 & 32.9 & [25.1,40.2] & 24 \\
04-30 & 10 & Cauchy & Huber IRLS & 246 & 80.1 & [74.3,85.1] & 24 \\
04-30 & 10 & Cauchy & Cauchy IRLS & 246 & 80.9 & [75.3,85.7] & 24 \\
04-30 & 10 & Cauchy & IGG~III IRLS & 246 & 69.5 & [63.7,75.3] & 24 \\
04-30 & 10 & Cauchy & Tukey IRLS & 246 & 76.8 & [70.7,82.5] & 24 \\
\bottomrule
\end{longtable}

\FloatBarrier
The electronic file \path{revision_20261103/results/sparse_acceptance.csv} contains all 432 policies for nine configurations. It gives dates, length, law, configuration, nominal target, fitted threshold, evaluation counts, rates and accepted-error median and P95. The Wilson column is the upper endpoint of a two-sided nominal 95\% interval under independent Bernoulli sampling. The companion file \path{revision_20261103/results/threshold_learning_intervals.csv} resamples calibration and evaluation blocks separately, refits the threshold and records zero-event resamples.
\begingroup
\small
\setlength{\tabcolsep}{4pt}
\setlength{\LTcapwidth}{0.9\linewidth}
\begin{longtable}{@{}l*{8}{r}@{}}
\caption{Cross-date acceptance after sparse injection under original-order Cholesky whitening. Thresholds calibrated on 30 April 2023 are applied unchanged on 31 March 2023. $A=n_a/n$ and $q=n_w/n_a$, where $n$ counts planned evaluation windows, $n_a$ acceptances and $n_w$ wrong accepted vectors. Both rates are in percent. A displayed $q$ of 0.00 means zero observed wrong fixes. The study reference for $q$ is 0.5\%.}\label{tab:sparse-cross-date-acceptance}\\
\toprule
Estimator & \multicolumn{2}{c}{Nominal 5\%} & \multicolumn{2}{c}{Nominal 10\%} & \multicolumn{2}{c}{Nominal 25\%} & \multicolumn{2}{c}{Nominal 50\%} \\
\cmidrule(lr){2-3}\cmidrule(lr){4-5}\cmidrule(lr){6-7}\cmidrule(lr){8-9}
 & $A$ (\%) & $q$ (\%) & $A$ (\%) & $q$ (\%) & $A$ (\%) & $q$ (\%) & $A$ (\%) & $q$ (\%) \\
\midrule
\endfirsthead
\multicolumn{9}{@{}l}{\tablename~\thetable\ (continued)} \\
\toprule
Estimator & \multicolumn{2}{c}{Nominal 5\%} & \multicolumn{2}{c}{Nominal 10\%} & \multicolumn{2}{c}{Nominal 25\%} & \multicolumn{2}{c}{Nominal 50\%} \\
\cmidrule(lr){2-3}\cmidrule(lr){4-5}\cmidrule(lr){6-7}\cmidrule(lr){8-9}
 & $A$ (\%) & $q$ (\%) & $A$ (\%) & $q$ (\%) & $A$ (\%) & $q$ (\%) & $A$ (\%) & $q$ (\%) \\
\midrule
\endhead
\bottomrule
\endlastfoot
\multicolumn{9}{@{}l}{\textbf{Gaussian mixture}, $N=5$, $n=535$} \\*
WLS & 7.29 & 46.15 & 11.40 & 50.82 & 25.98 & 71.94 & 48.79 & 80.08 \\*
GMM, $K=2$ & 6.36 & 20.59 & 10.47 & 23.21 & 25.79 & 39.13 & 47.10 & 44.44 \\*
GMM, $K=3$ & 4.86 & 11.54 & 11.59 & 27.42 & 25.23 & 40.00 & 46.54 & 47.79 \\*
GMM, $K=4$ & 4.49 & 8.33 & 11.78 & 25.40 & 25.42 & 41.18 & 46.73 & 47.20 \\*
Huber IRLS & 2.43 & 0.00 & 7.85 & 0.00 & 24.49 & 0.00 & 51.78 & 0.00 \\*
Cauchy IRLS & 4.67 & 0.00 & 9.53 & 0.00 & 25.05 & 0.00 & 55.14 & 0.34 \\*
IGG~III IRLS & 4.30 & 4.35 & 7.10 & 2.63 & 21.50 & 1.74 & 46.36 & 1.21 \\*
Tukey IRLS & 3.55 & 5.26 & 7.85 & 2.38 & 25.05 & 0.75 & 48.97 & 0.38 \\*
IDS--ILS & 5.23 & 0.00 & 10.65 & 0.00 & 27.85 & 3.36 & 52.34 & 8.21 \\
\addlinespace[4pt]
\multicolumn{9}{@{}l}{\textbf{Cauchy}, $N=5$, $n=535$} \\*
WLS & 3.36 & 38.89 & 11.59 & 53.23 & 28.97 & 65.16 & 55.70 & 75.50 \\*
GMM, $K=2$ & 4.67 & 4.00 & 8.97 & 14.58 & 24.30 & 31.54 & 46.17 & 34.41 \\*
GMM, $K=3$ & 5.42 & 6.90 & 8.60 & 15.22 & 25.05 & 35.82 & 49.72 & 40.60 \\*
GMM, $K=4$ & 6.36 & 8.82 & 9.16 & 18.37 & 25.42 & 33.82 & 50.65 & 39.48 \\*
Huber IRLS & 4.11 & 0.00 & 7.85 & 0.00 & 20.37 & 0.00 & 46.17 & 0.00 \\*
Cauchy IRLS & 4.11 & 0.00 & 7.66 & 0.00 & 19.81 & 0.00 & 48.22 & 0.00 \\*
IGG~III IRLS & 5.05 & 0.00 & 8.22 & 0.00 & 24.11 & 0.00 & 46.36 & 1.21 \\*
Tukey IRLS & 4.86 & 0.00 & 9.91 & 0.00 & 22.43 & 0.00 & 47.10 & 0.00 \\*
IDS--ILS & 4.11 & 4.55 & 9.72 & 1.92 & 23.55 & 4.76 & 51.96 & 14.03 \\
\addlinespace[4pt]
\multicolumn{9}{@{}l}{\textbf{Gaussian mixture}, $N=10$, $n=247$} \\*
WLS & 8.50 & 52.38 & 14.17 & 68.57 & 28.34 & 77.14 & 52.63 & 86.92 \\*
GMM, $K=2$ & 6.48 & 43.75 & 11.74 & 37.93 & 25.10 & 46.77 & 52.23 & 55.81 \\*
GMM, $K=3$ & 4.86 & 41.67 & 8.10 & 35.00 & 27.13 & 46.27 & 50.20 & 57.26 \\*
GMM, $K=4$ & 6.07 & 46.67 & 9.72 & 41.67 & 26.72 & 45.45 & 53.44 & 56.06 \\*
Huber IRLS & 5.67 & 0.00 & 8.91 & 0.00 & 25.51 & 0.00 & 47.77 & 0.00 \\*
Cauchy IRLS & 4.05 & 0.00 & 8.10 & 0.00 & 26.72 & 0.00 & 51.01 & 0.00 \\*
IGG~III IRLS & 4.05 & 0.00 & 11.74 & 0.00 & 25.10 & 0.00 & 47.37 & 0.00 \\*
Tukey IRLS & 4.45 & 0.00 & 10.93 & 0.00 & 24.29 & 0.00 & 50.61 & 0.00 \\*
IDS--ILS & 3.64 & 0.00 & 7.29 & 0.00 & 22.67 & 3.57 & 40.89 & 1.98 \\
\addlinespace[4pt]
\multicolumn{9}{@{}l}{\textbf{Cauchy}, $N=10$, $n=247$} \\*
WLS & 9.72 & 41.67 & 15.79 & 46.15 & 23.89 & 55.93 & 57.89 & 74.83 \\*
GMM, $K=2$ & 4.86 & 33.33 & 12.96 & 34.38 & 27.13 & 40.30 & 51.01 & 44.44 \\*
GMM, $K=3$ & 4.45 & 9.09 & 9.72 & 20.83 & 25.51 & 49.21 & 45.75 & 50.44 \\*
GMM, $K=4$ & 4.45 & 9.09 & 9.72 & 33.33 & 24.70 & 40.98 & 48.58 & 47.50 \\*
Huber IRLS & 7.29 & 0.00 & 11.74 & 0.00 & 26.72 & 0.00 & 49.39 & 0.00 \\*
Cauchy IRLS & 4.86 & 0.00 & 11.74 & 0.00 & 26.32 & 0.00 & 48.18 & 0.00 \\*
IGG~III IRLS & 5.26 & 0.00 & 9.31 & 0.00 & 25.91 & 0.00 & 47.37 & 0.85 \\*
Tukey IRLS & 5.26 & 0.00 & 8.91 & 0.00 & 25.51 & 1.59 & 45.75 & 0.88 \\*
IDS--ILS & 7.69 & 0.00 & 14.98 & 0.00 & 28.34 & 0.00 & 53.44 & 0.00 \\
\end{longtable}
\endgroup

\FloatBarrier
\section{Residual-screening and correlated-contamination comparisons}
All eight weighting configurations and IDS--ILS are evaluated on the same inputs within each cell. IDS--ILS screens supplied nominally whitened rows and retains the full ambiguity dimension. It returns no integer vector if rank, redundancy or search requirements fail. The Markov design adds controlled errors to rows of the raw double-difference stack as specified in Section 4.3 of the main text.
\begin{table}[htbp]
\centering\small
\caption{Full-vector recovery (\%) with Markov innovations in raw DD code rows. Entries give recovery for 31 March / 30 April 2023. Planned windows are 535 / 532 at $N=5$ and 247 / 246 at $N=10$. All fits use code and phase. Full-vector recovery uses the internal reference and counts incomplete searches as failures.}
\begin{tabular}{lrrrr}
\toprule
Method & Mixture, $N=5$ & Mixture, $N=10$ & Cauchy, $N=5$ & Cauchy, $N=10$ \\
\midrule
WLS & 82.6/81.6 & 86.6/88.2 & 85.0/84.0 & 91.1/89.8 \\
GMM, K=2 & 93.1/94.9 & 98.4/98.4 & 93.6/95.7 & 98.8/98.4 \\
GMM, K=3 & 93.1/94.9 & 98.4/98.4 & 93.6/95.7 & 98.8/98.4 \\
GMM, K=4 & 93.1/94.9 & 98.4/98.4 & 93.6/95.7 & 98.8/98.4 \\
Huber IRLS & 96.4/97.6 & 99.2/99.2 & 97.4/98.5 & 99.2/99.2 \\
Cauchy IRLS & 96.6/97.7 & 99.2/98.8 & 97.8/98.1 & 99.2/98.8 \\
IGG~III IRLS & 95.7/97.4 & 99.2/99.2 & 96.6/97.2 & 99.2/99.6 \\
Tukey IRLS & 96.3/98.1 & 99.2/98.8 & 98.1/97.9 & 99.6/99.2 \\
IDS--ILS & 94.4/95.1 & 99.2/99.2 & 95.0/96.6 & 98.8/99.6 \\
\bottomrule
\end{tabular}
\end{table}

\begin{table}[htbp]
\centering\small
\caption{Full-vector recovery (\%) with Markov innovations in raw DD code and phase rows. Entries give recovery for 31 March / 30 April 2023. Planned windows are 535 / 532 at $N=5$ and 247 / 246 at $N=10$. All fits use code and phase. Full-vector recovery uses the internal reference and counts incomplete searches as failures.}
\begin{tabular}{lrrrr}
\toprule
Method & Mixture, $N=5$ & Mixture, $N=10$ & Cauchy, $N=5$ & Cauchy, $N=10$ \\
\midrule
WLS & 67.3/68.0 & 69.2/64.6 & 69.3/72.9 & 71.3/68.7 \\
GMM, K=2 & 75.5/75.9 & 78.9/75.6 & 77.4/78.8 & 76.9/77.6 \\
GMM, K=3 & 75.5/75.9 & 78.9/75.6 & 77.4/78.8 & 77.3/77.6 \\
GMM, K=4 & 75.5/75.9 & 78.9/75.6 & 77.4/78.8 & 77.3/78.0 \\
Huber IRLS & 86.7/85.5 & 93.9/92.7 & 87.9/90.6 & 96.0/95.9 \\
Cauchy IRLS & 86.5/84.8 & 93.1/91.9 & 87.7/90.0 & 95.1/95.5 \\
IGG~III IRLS & 82.1/83.1 & 87.0/86.2 & 82.4/85.7 & 89.1/90.2 \\
Tukey IRLS & 83.2/82.1 & 88.3/87.4 & 83.7/87.6 & 91.1/91.9 \\
IDS--ILS & 87.5/87.4 & 95.5/96.7 & 86.0/89.1 & 98.0/95.5 \\
\bottomrule
\end{tabular}
\end{table}

\begin{longtable}{lllrrrrrr}
\caption{Cross-date Markov acceptance ranges at nominal calibration targets of 10\% and 50\%. Each range spans both date directions and lengths within the stated law and injection composition. Extrema of $n_a$, $q$ and $V_{10}$ can come from different cells. Evaluation rates are shown in percent. The archive contains full results for 576 policies.}\\
\toprule & & & \multicolumn{3}{c}{Nominal calibration target 10\%} & \multicolumn{3}{c}{Nominal calibration target 50\%}\\
\shortstack[l]{Raw DD rows\\injected} & Law & Method & \shortstack{Accepted\\count range} & $q$ (\%) & $V_{10}$ (\%) & \shortstack{Accepted\\count range} & $q$ (\%) & $V_{10}$ (\%)\\\midrule\endfirsthead
\multicolumn{9}{l}{\tablename~\thetable\ (continued)} \\
\toprule & & & \multicolumn{3}{c}{Nominal calibration target 10\%} & \multicolumn{3}{c}{Nominal calibration target 50\%}\\
\shortstack[l]{Raw DD rows\\injected} & Law & Method & \shortstack{Accepted\\count range} & $q$ (\%) & $V_{10}$ (\%) & \shortstack{Accepted\\count range} & $q$ (\%) & $V_{10}$ (\%)\\\midrule\endhead
\bottomrule\endfoot
Code only & Mixture & WLS & 24--64 & 0.00 & 8.8--12.0 & 122--272 & 0.00 & 47.9--51.1\\
Code only & Mixture & GMM, $K=2$ & 25--68 & 0.00 & 7.9--12.8 & 121--312 & 0.64--1.59 & 42.1--58.3\\
Code only & Mixture & GMM, $K=3$ & 25--69 & 0.00 & 7.9--13.0 & 120--311 & 0.32--0.88 & 41.9--58.3\\
Code only & Mixture & GMM, $K=4$ & 25--69 & 0.00 & 7.9--13.0 & 123--311 & 0.00--0.88 & 41.9--58.3\\
Code only & Mixture & Huber IRLS & 25--60 & 0.00 & 9.2--11.3 & 111--277 & 0.00 & 44.9--55.3\\
Code only & Mixture & Cauchy IRLS & 21--55 & 0.00 & 8.5--11.0 & 114--286 & 0.00 & 46.2--56.9\\
Code only & Mixture & IGG~III IRLS & 24--61 & 0.00 & 8.3--11.4 & 112--294 & 0.00--0.34 & 44.9--55.3\\
Code only & Mixture & Tukey IRLS & 25--55 & 0.00 & 9.8--10.3 & 117--292 & 0.00 & 43.9--54.9\\
Code only & Mixture & IDS--ILS & 17--67 & 0.00 & 6.9--13.4 & 121--273 & 0.00 & 46.7--51.6\\
\addlinespace
Code only & Cauchy & WLS & 10--57 & 0.00 & 4.1--14.6 & 115--295 & 0.00--0.68 & 43.6--54.8\\
Code only & Cauchy & GMM, $K=2$ & 23--57 & 0.00 & 9.3--11.8 & 117--303 & 0.00--1.27 & 43.7--56.4\\
Code only & Cauchy & GMM, $K=3$ & 16--56 & 0.00 & 6.5--13.8 & 118--296 & 0.00--1.26 & 43.9--55.3\\
Code only & Cauchy & GMM, $K=4$ & 16--58 & 0.00 & 6.5--13.0 & 119--292 & 0.00--1.26 & 44.1--54.5\\
Code only & Cauchy & Huber IRLS & 21--59 & 0.00 & 8.5--12.1 & 118--271 & 0.00 & 47.8--56.1\\
Code only & Cauchy & Cauchy IRLS & 21--57 & 0.00 & 8.5--10.7 & 118--278 & 0.00 & 47.8--53.3\\
Code only & Cauchy & IGG~III IRLS & 24--58 & 0.00 & 9.2--11.3 & 119--290 & 0.00 & 46.2--54.5\\
Code only & Cauchy & Tukey IRLS & 21--55 & 0.00 & 8.5--10.6 & 119--289 & 0.00 & 45.4--55.3\\
Code only & Cauchy & IDS--ILS & 19--66 & 0.00 & 7.7--14.2 & 122--279 & 0.00--0.36 & 48.5--52.0\\
\addlinespace
Code + phase & Mixture & WLS & 19--64 & 0.00--2.63 & 7.7--15.4 & 115--279 & 3.89--11.30 & 45.9--52.6\\
Code + phase & Mixture & GMM, $K=2$ & 23--58 & 0.00 & 8.6--11.7 & 117--306 & 6.30--8.50 & 43.0--54.3\\
Code + phase & Mixture & GMM, $K=3$ & 21--62 & 0.00--4.17 & 8.5--11.7 & 108--300 & 10.17--12.95 & 40.9--52.3\\
Code + phase & Mixture & GMM, $K=4$ & 25--62 & 0.00--4.08 & 9.0--11.7 & 113--292 & 10.27--13.43 & 41.3--50.9\\
Code + phase & Mixture & Huber IRLS & 21--55 & 0.00--1.96 & 8.5--11.8 & 106--282 & 0.00--1.77 & 42.9--54.9\\
Code + phase & Mixture & Cauchy IRLS & 19--57 & 0.00--1.75 & 7.7--11.8 & 111--277 & 0.73--2.17 & 44.9--55.7\\
Code + phase & Mixture & IGG~III IRLS & 23--54 & 0.00--4.35 & 9.3--11.0 & 111--282 & 0.90--3.55 & 44.9--54.9\\
Code + phase & Mixture & Tukey IRLS & 22--55 & 0.00--3.85 & 8.9--10.5 & 109--270 & 0.92--2.76 & 44.1--58.9\\
Code + phase & Mixture & IDS--ILS & 19--61 & 0.00--1.64 & 7.7--11.7 & 122--268 & 0.00--2.26 & 49.4--50.4\\
\addlinespace
Code + phase & Cauchy & WLS & 22--57 & 0.00--2.94 & 8.9--13.8 & 110--270 & 4.21--8.27 & 44.5--52.8\\
Code + phase & Cauchy & GMM, $K=2$ & 23--56 & 0.00 & 9.3--11.3 & 122--276 & 4.28--9.02 & 46.9--50.4\\
Code + phase & Cauchy & GMM, $K=3$ & 23--55 & 0.00--1.82 & 9.3--11.3 & 120--284 & 4.49--12.50 & 44.7--51.1\\
Code + phase & Cauchy & GMM, $K=4$ & 22--54 & 0.00--1.85 & 8.9--12.1 & 119--280 & 4.84--14.29 & 45.0--50.4\\
Code + phase & Cauchy & Huber IRLS & 20--55 & 0.00--1.82 & 8.1--13.4 & 111--283 & 0.00--1.77 & 44.9--54.5\\
Code + phase & Cauchy & Cauchy IRLS & 22--54 & 0.00--1.89 & 8.9--12.1 & 109--287 & 0.00--2.44 & 44.1--54.9\\
Code + phase & Cauchy & IGG~III IRLS & 20--63 & 0.00 & 8.1--12.1 & 113--312 & 0.88--2.88 & 44.3--58.6\\
Code + phase & Cauchy & Tukey IRLS & 20--54 & 0.00 & 8.1--10.9 & 118--320 & 0.74--2.81 & 43.7--60.2\\
Code + phase & Cauchy & IDS--ILS & 22--57 & 0.00--1.75 & 8.9--11.3 & 111--279 & 0.00--2.15 & 44.9--55.3\\
\addlinespace
\end{longtable}

\begin{longtable}{llllrll}
\caption{Markov full-vector recovery contrasts. Differences are the listed method minus the stated comparator, in percentage points. Positive values favour the listed method. The pointwise 95\% percentile intervals use 2,000 paired 60-minute block resamples within each date. Here 03-31 and 04-30 denote 31 March and 30 April 2023. Other block settings are in the archive.}\\
\toprule \shortstack[l]{Raw DD rows\\injected} & Law & Method & Comparator & $N$ & 03-31 & 04-30\\\midrule\endfirsthead
\multicolumn{7}{l}{\tablename~\thetable\ (continued)} \\
\toprule \shortstack[l]{Raw DD rows\\injected} & Law & Method & Comparator & $N$ & 03-31 & 04-30\\\midrule\endhead
\bottomrule\endfoot
Code only & Mixture & GMM, $K=2$ & WLS & 5 & 10.5 [7.6, 13.9] & 13.3 [10.6, 16.0]\\
Code only & Mixture & GMM, $K=2$ & WLS & 10 & 11.7 [8.1, 15.2] & 10.2 [5.9, 14.7]\\
Code only & Mixture & GMM, $K=3$ & WLS & 5 & 10.5 [7.6, 13.9] & 13.3 [10.6, 16.0]\\
Code only & Mixture & GMM, $K=3$ & WLS & 10 & 11.7 [8.1, 15.2] & 10.2 [5.9, 14.7]\\
Code only & Mixture & GMM, $K=4$ & WLS & 5 & 10.5 [7.6, 13.9] & 13.3 [10.6, 16.0]\\
Code only & Mixture & GMM, $K=4$ & WLS & 10 & 11.7 [8.1, 15.2] & 10.2 [5.9, 14.7]\\
Code only & Mixture & Huber IRLS & WLS & 5 & 13.8 [10.4, 17.6] & 16.0 [12.9, 18.8]\\
Code only & Mixture & Huber IRLS & WLS & 10 & 12.6 [8.3, 16.6] & 11.0 [6.9, 15.2]\\
Code only & Mixture & Cauchy IRLS & WLS & 5 & 14.0 [10.5, 17.8] & 16.2 [13.1, 19.0]\\
Code only & Mixture & Cauchy IRLS & WLS & 10 & 12.6 [8.3, 16.6] & 10.6 [6.5, 15.0]\\
Code only & Mixture & IGG~III IRLS & WLS & 5 & 13.1 [9.9, 16.7] & 15.8 [12.7, 18.7]\\
Code only & Mixture & IGG~III IRLS & WLS & 10 & 12.6 [8.3, 16.6] & 11.0 [6.9, 15.2]\\
Code only & Mixture & Tukey IRLS & WLS & 5 & 13.6 [10.3, 17.4] & 16.5 [13.4, 19.5]\\
Code only & Mixture & Tukey IRLS & WLS & 10 & 12.6 [8.3, 16.6] & 10.6 [6.5, 15.0]\\
Code only & Mixture & IDS--ILS & WLS & 5 & 11.8 [9.5, 14.1] & 13.5 [10.4, 16.6]\\
Code only & Mixture & IDS--ILS & WLS & 10 & 12.6 [8.3, 16.6] & 11.0 [7.3, 14.7]\\
Code only & Mixture & WLS & IDS--ILS & 5 & -11.8 [-14.1, -9.5] & -13.5 [-16.6, -10.4]\\
Code only & Mixture & WLS & IDS--ILS & 10 & -12.6 [-16.6, -8.3] & -11.0 [-14.7, -7.3]\\
Code only & Mixture & GMM, $K=2$ & IDS--ILS & 5 & -1.3 [-3.0, 0.6] & -0.2 [-2.3, 2.0]\\
Code only & Mixture & GMM, $K=2$ & IDS--ILS & 10 & -0.8 [-2.0, 0.0] & -0.8 [-2.8, 1.6]\\
Code only & Mixture & GMM, $K=3$ & IDS--ILS & 5 & -1.3 [-3.0, 0.6] & -0.2 [-2.3, 2.0]\\
Code only & Mixture & GMM, $K=3$ & IDS--ILS & 10 & -0.8 [-2.0, 0.0] & -0.8 [-2.8, 1.6]\\
Code only & Mixture & GMM, $K=4$ & IDS--ILS & 5 & -1.3 [-3.0, 0.6] & -0.2 [-2.3, 2.0]\\
Code only & Mixture & GMM, $K=4$ & IDS--ILS & 10 & -0.8 [-2.0, 0.0] & -0.8 [-2.8, 1.6]\\
Code only & Mixture & Huber IRLS & IDS--ILS & 5 & 2.1 [0.0, 4.3] & 2.4 [0.7, 4.5]\\
Code only & Mixture & Huber IRLS & IDS--ILS & 10 & 0.0 [0.0, 0.0] & 0.0 [-1.6, 2.1]\\
Code only & Mixture & Cauchy IRLS & IDS--ILS & 5 & 2.2 [0.4, 4.3] & 2.6 [0.9, 4.7]\\
Code only & Mixture & Cauchy IRLS & IDS--ILS & 10 & 0.0 [0.0, 0.0] & -0.4 [-2.4, 2.0]\\
Code only & Mixture & IGG~III IRLS & IDS--ILS & 5 & 1.3 [-0.8, 3.6] & 2.3 [0.6, 4.1]\\
Code only & Mixture & IGG~III IRLS & IDS--ILS & 10 & 0.0 [0.0, 0.0] & 0.0 [-1.6, 2.1]\\
Code only & Mixture & Tukey IRLS & IDS--ILS & 5 & 1.9 [-0.2, 4.1] & 3.0 [1.3, 5.0]\\
Code only & Mixture & Tukey IRLS & IDS--ILS & 10 & 0.0 [0.0, 0.0] & -0.4 [-2.4, 2.0]\\
\addlinespace
Code only & Cauchy & GMM, $K=2$ & WLS & 5 & 8.6 [5.7, 11.5] & 11.7 [9.0, 14.5]\\
Code only & Cauchy & GMM, $K=2$ & WLS & 10 & 7.7 [4.4, 11.0] & 8.5 [5.7, 11.4]\\
Code only & Cauchy & GMM, $K=3$ & WLS & 5 & 8.6 [5.7, 11.5] & 11.7 [9.0, 14.5]\\
Code only & Cauchy & GMM, $K=3$ & WLS & 10 & 7.7 [4.4, 11.0] & 8.5 [5.7, 11.4]\\
Code only & Cauchy & GMM, $K=4$ & WLS & 5 & 8.6 [5.7, 11.5] & 11.7 [9.0, 14.5]\\
Code only & Cauchy & GMM, $K=4$ & WLS & 10 & 7.7 [4.4, 11.0] & 8.5 [5.7, 11.4]\\
Code only & Cauchy & Huber IRLS & WLS & 5 & 12.3 [9.1, 15.7] & 14.5 [11.6, 17.4]\\
Code only & Cauchy & Huber IRLS & WLS & 10 & 8.1 [4.8, 11.3] & 9.3 [6.0, 13.1]\\
Code only & Cauchy & Cauchy IRLS & WLS & 5 & 12.7 [9.2, 16.3] & 14.1 [11.3, 17.1]\\
Code only & Cauchy & Cauchy IRLS & WLS & 10 & 8.1 [4.8, 11.3] & 8.9 [5.8, 12.1]\\
Code only & Cauchy & IGG~III IRLS & WLS & 5 & 11.6 [8.0, 15.5] & 13.2 [10.7, 15.9]\\
Code only & Cauchy & IGG~III IRLS & WLS & 10 & 8.1 [4.8, 11.3] & 9.8 [6.5, 13.4]\\
Code only & Cauchy & Tukey IRLS & WLS & 5 & 13.1 [9.4, 17.0] & 13.9 [11.3, 16.7]\\
Code only & Cauchy & Tukey IRLS & WLS & 10 & 8.5 [4.9, 12.4] & 9.3 [6.0, 13.1]\\
Code only & Cauchy & IDS--ILS & WLS & 5 & 9.9 [7.2, 12.9] & 12.6 [10.0, 15.4]\\
Code only & Cauchy & IDS--ILS & WLS & 10 & 7.7 [4.4, 11.1] & 9.8 [6.5, 13.3]\\
Code only & Cauchy & WLS & IDS--ILS & 5 & -9.9 [-12.9, -7.2] & -12.6 [-15.4, -10.0]\\
Code only & Cauchy & WLS & IDS--ILS & 10 & -7.7 [-11.1, -4.4] & -9.8 [-13.3, -6.5]\\
Code only & Cauchy & GMM, $K=2$ & IDS--ILS & 5 & -1.3 [-3.5, 0.7] & -0.9 [-3.4, 0.9]\\
Code only & Cauchy & GMM, $K=2$ & IDS--ILS & 10 & 0.0 [-1.2, 1.2] & -1.2 [-2.8, 0.0]\\
Code only & Cauchy & GMM, $K=3$ & IDS--ILS & 5 & -1.3 [-3.5, 0.7] & -0.9 [-3.4, 0.9]\\
Code only & Cauchy & GMM, $K=3$ & IDS--ILS & 10 & 0.0 [-1.2, 1.2] & -1.2 [-2.8, 0.0]\\
Code only & Cauchy & GMM, $K=4$ & IDS--ILS & 5 & -1.3 [-3.5, 0.7] & -0.9 [-3.4, 0.9]\\
Code only & Cauchy & GMM, $K=4$ & IDS--ILS & 10 & 0.0 [-1.2, 1.2] & -1.2 [-2.8, 0.0]\\
Code only & Cauchy & Huber IRLS & IDS--ILS & 5 & 2.4 [1.1, 3.9] & 1.9 [0.2, 3.6]\\
Code only & Cauchy & Huber IRLS & IDS--ILS & 10 & 0.4 [-0.8, 1.7] & -0.4 [-1.2, 0.0]\\
Code only & Cauchy & Cauchy IRLS & IDS--ILS & 5 & 2.8 [1.1, 4.6] & 1.5 [0.0, 3.2]\\
Code only & Cauchy & Cauchy IRLS & IDS--ILS & 10 & 0.4 [-0.8, 1.7] & -0.8 [-2.0, 0.0]\\
Code only & Cauchy & IGG~III IRLS & IDS--ILS & 5 & 1.7 [-0.7, 4.2] & 0.6 [-1.0, 1.9]\\*
Code only & Cauchy & IGG~III IRLS & IDS--ILS & 10 & 0.4 [-0.8, 1.7] & 0.0 [-1.2, 1.2]\\*
Code only & Cauchy & Tukey IRLS & IDS--ILS & 5 & 3.2 [1.1, 5.5] & 1.3 [-0.2, 2.8]\\*
Code only & Cauchy & Tukey IRLS & IDS--ILS & 10 & 0.8 [-0.8, 2.4] & -0.4 [-1.2, 0.0]\\
\addlinespace
Code + phase & Mixture & GMM, $K=2$ & WLS & 5 & 8.2 [5.5, 11.2] & 7.9 [4.9, 10.6]\\
Code + phase & Mixture & GMM, $K=2$ & WLS & 10 & 9.7 [5.0, 14.3] & 11.0 [5.3, 16.7]\\
Code + phase & Mixture & GMM, $K=3$ & WLS & 5 & 8.2 [5.5, 11.2] & 7.9 [4.9, 10.6]\\
Code + phase & Mixture & GMM, $K=3$ & WLS & 10 & 9.7 [5.0, 14.3] & 11.0 [5.3, 16.7]\\
Code + phase & Mixture & GMM, $K=4$ & WLS & 5 & 8.2 [5.5, 11.2] & 7.9 [4.9, 10.6]\\
Code + phase & Mixture & GMM, $K=4$ & WLS & 10 & 9.7 [5.0, 14.3] & 11.0 [5.3, 16.7]\\
Code + phase & Mixture & Huber IRLS & WLS & 5 & 19.4 [15.8, 23.4] & 17.5 [14.8, 20.2]\\
Code + phase & Mixture & Huber IRLS & WLS & 10 & 24.7 [18.9, 30.4] & 28.0 [21.7, 34.7]\\
Code + phase & Mixture & Cauchy IRLS & WLS & 5 & 19.3 [15.6, 23.2] & 16.7 [14.2, 19.3]\\
Code + phase & Mixture & Cauchy IRLS & WLS & 10 & 23.9 [18.2, 29.5] & 27.2 [21.5, 33.1]\\
Code + phase & Mixture & IGG~III IRLS & WLS & 5 & 14.8 [11.6, 18.0] & 15.0 [11.7, 18.3]\\
Code + phase & Mixture & IGG~III IRLS & WLS & 10 & 17.8 [12.6, 23.4] & 21.5 [15.5, 28.3]\\
Code + phase & Mixture & Tukey IRLS & WLS & 5 & 15.9 [12.5, 19.3] & 14.1 [11.0, 17.3]\\
Code + phase & Mixture & Tukey IRLS & WLS & 10 & 19.0 [13.5, 24.6] & 22.8 [16.7, 28.9]\\
Code + phase & Mixture & IDS--ILS & WLS & 5 & 20.2 [16.8, 23.9] & 19.4 [15.9, 23.2]\\
Code + phase & Mixture & IDS--ILS & WLS & 10 & 26.3 [21.2, 31.3] & 32.1 [26.6, 38.3]\\
Code + phase & Mixture & WLS & IDS--ILS & 5 & -20.2 [-23.9, -16.8] & -19.4 [-23.2, -15.9]\\
Code + phase & Mixture & WLS & IDS--ILS & 10 & -26.3 [-31.3, -21.2] & -32.1 [-38.3, -26.6]\\
Code + phase & Mixture & GMM, $K=2$ & IDS--ILS & 5 & -12.0 [-15.7, -8.2] & -11.5 [-15.4, -7.6]\\
Code + phase & Mixture & GMM, $K=2$ & IDS--ILS & 10 & -16.6 [-20.8, -12.4] & -21.1 [-25.8, -16.5]\\
Code + phase & Mixture & GMM, $K=3$ & IDS--ILS & 5 & -12.0 [-15.7, -8.2] & -11.5 [-15.4, -7.6]\\
Code + phase & Mixture & GMM, $K=3$ & IDS--ILS & 10 & -16.6 [-20.8, -12.4] & -21.1 [-25.8, -16.5]\\
Code + phase & Mixture & GMM, $K=4$ & IDS--ILS & 5 & -12.0 [-15.7, -8.2] & -11.5 [-15.4, -7.6]\\
Code + phase & Mixture & GMM, $K=4$ & IDS--ILS & 10 & -16.6 [-20.8, -12.4] & -21.1 [-25.8, -16.5]\\
Code + phase & Mixture & Huber IRLS & IDS--ILS & 5 & -0.7 [-3.4, 2.2] & -1.9 [-4.5, 0.8]\\
Code + phase & Mixture & Huber IRLS & IDS--ILS & 10 & -1.6 [-4.5, 1.2] & -4.1 [-7.0, -1.6]\\
Code + phase & Mixture & Cauchy IRLS & IDS--ILS & 5 & -0.9 [-3.7, 1.8] & -2.6 [-5.4, 0.2]\\
Code + phase & Mixture & Cauchy IRLS & IDS--ILS & 10 & -2.4 [-5.2, 0.4] & -4.9 [-7.7, -2.0]\\
Code + phase & Mixture & IGG~III IRLS & IDS--ILS & 5 & -5.4 [-8.7, -2.3] & -4.3 [-7.3, -1.3]\\
Code + phase & Mixture & IGG~III IRLS & IDS--ILS & 10 & -8.5 [-12.2, -5.2] & -10.6 [-14.7, -6.9]\\
Code + phase & Mixture & Tukey IRLS & IDS--ILS & 5 & -4.3 [-7.3, -1.3] & -5.3 [-8.4, -2.1]\\
Code + phase & Mixture & Tukey IRLS & IDS--ILS & 10 & -7.3 [-10.7, -3.7] & -9.3 [-12.5, -6.4]\\
\addlinespace
Code + phase & Cauchy & GMM, $K=2$ & WLS & 5 & 8.0 [4.7, 11.2] & 5.8 [3.2, 8.5]\\
Code + phase & Cauchy & GMM, $K=2$ & WLS & 10 & 5.7 [-0.4, 12.2] & 8.9 [4.5, 13.7]\\
Code + phase & Cauchy & GMM, $K=3$ & WLS & 5 & 8.0 [4.7, 11.2] & 5.8 [3.2, 8.5]\\
Code + phase & Cauchy & GMM, $K=3$ & WLS & 10 & 6.1 [-0.4, 13.0] & 8.9 [4.5, 13.7]\\
Code + phase & Cauchy & GMM, $K=4$ & WLS & 5 & 8.0 [4.7, 11.2] & 5.8 [3.2, 8.5]\\
Code + phase & Cauchy & GMM, $K=4$ & WLS & 10 & 6.1 [-0.4, 13.0] & 9.3 [4.9, 14.2]\\
Code + phase & Cauchy & Huber IRLS & WLS & 5 & 18.5 [14.7, 22.1] & 17.7 [14.6, 20.8]\\
Code + phase & Cauchy & Huber IRLS & WLS & 10 & 24.7 [18.3, 31.3] & 27.2 [21.4, 33.6]\\
Code + phase & Cauchy & Cauchy IRLS & WLS & 5 & 18.3 [14.1, 22.4] & 17.1 [14.0, 20.4]\\
Code + phase & Cauchy & Cauchy IRLS & WLS & 10 & 23.9 [16.9, 31.2] & 26.8 [21.2, 32.8]\\
Code + phase & Cauchy & IGG~III IRLS & WLS & 5 & 13.1 [9.3, 16.9] & 12.8 [9.9, 15.8]\\
Code + phase & Cauchy & IGG~III IRLS & WLS & 10 & 17.8 [10.4, 26.1] & 21.5 [16.4, 26.9]\\
Code + phase & Cauchy & Tukey IRLS & WLS & 5 & 14.4 [10.7, 18.0] & 14.7 [11.8, 17.6]\\
Code + phase & Cauchy & Tukey IRLS & WLS & 10 & 19.8 [12.3, 27.9] & 23.2 [18.0, 28.7]\\
Code + phase & Cauchy & IDS--ILS & WLS & 5 & 16.6 [12.8, 20.5] & 16.2 [12.8, 19.4]\\
Code + phase & Cauchy & IDS--ILS & WLS & 10 & 26.7 [19.5, 34.3] & 26.8 [21.1, 33.5]\\
Code + phase & Cauchy & WLS & IDS--ILS & 5 & -16.6 [-20.5, -12.8] & -16.2 [-19.4, -12.8]\\
Code + phase & Cauchy & WLS & IDS--ILS & 10 & -26.7 [-34.3, -19.5] & -26.8 [-33.5, -21.1]\\
Code + phase & Cauchy & GMM, $K=2$ & IDS--ILS & 5 & -8.6 [-11.2, -6.2] & -10.3 [-13.4, -7.3]\\
Code + phase & Cauchy & GMM, $K=2$ & IDS--ILS & 10 & -21.1 [-26.4, -16.0] & -17.9 [-22.6, -13.5]\\
Code + phase & Cauchy & GMM, $K=3$ & IDS--ILS & 5 & -8.6 [-11.2, -6.2] & -10.3 [-13.4, -7.3]\\
Code + phase & Cauchy & GMM, $K=3$ & IDS--ILS & 10 & -20.6 [-26.1, -15.7] & -17.9 [-22.6, -13.5]\\
Code + phase & Cauchy & GMM, $K=4$ & IDS--ILS & 5 & -8.6 [-11.2, -6.2] & -10.3 [-13.4, -7.3]\\
Code + phase & Cauchy & GMM, $K=4$ & IDS--ILS & 10 & -20.6 [-26.1, -15.7] & -17.5 [-22.4, -12.8]\\
Code + phase & Cauchy & Huber IRLS & IDS--ILS & 5 & 1.9 [-0.4, 4.1] & 1.5 [-1.9, 5.1]\\
Code + phase & Cauchy & Huber IRLS & IDS--ILS & 10 & -2.0 [-4.8, 0.8] & 0.4 [-2.0, 2.8]\\
Code + phase & Cauchy & Cauchy IRLS & IDS--ILS & 5 & 1.7 [-0.7, 4.1] & 0.9 [-2.5, 4.5]\\
Code + phase & Cauchy & Cauchy IRLS & IDS--ILS & 10 & -2.8 [-6.0, 0.0] & 0.0 [-2.9, 2.8]\\
Code + phase & Cauchy & IGG~III IRLS & IDS--ILS & 5 & -3.6 [-6.6, -0.7] & -3.4 [-7.3, 0.0]\\
Code + phase & Cauchy & IGG~III IRLS & IDS--ILS & 10 & -8.9 [-12.4, -5.7] & -5.3 [-10.2, -0.8]\\
Code + phase & Cauchy & Tukey IRLS & IDS--ILS & 5 & -2.2 [-5.8, 1.1] & -1.5 [-5.1, 2.0]\\
Code + phase & Cauchy & Tukey IRLS & IDS--ILS & 10 & -6.9 [-10.6, -3.4] & -3.7 [-8.5, 0.4]\\
\end{longtable}

\begin{longtable}{lrlrrrr}
\caption{Timing by date, record length and injection law. Each row summarizes 12 measured invocations from four inputs and three measured repetitions. The completed count is shown for every row.}\label{tab:runtime_full}\\
\toprule
Date & $N$ & Input / method & Median (ms) & P95 (ms) & Completed & \shortstack{Median\\iterations} \\
\midrule\endfirsthead
\multicolumn{7}{l}{\tablename~\thetable\ (continued)} \\
\toprule
Date & $N$ & Input / method & Median (ms) & P95 (ms) & Completed & \shortstack{Median\\iterations} \\
\midrule\endhead
03-31 & 5 & Cauchy: WLS & 2.92 & 3.88 & 12/12 & 1 \\
03-31 & 5 & Cauchy: GMM, $K=2$ & 22.65 & 41.02 & 12/12 & 11 \\
03-31 & 5 & Cauchy: GMM, $K=3$ & 23.70 & 41.14 & 12/12 & 11 \\
03-31 & 5 & Cauchy: GMM, $K=4$ & 23.28 & 48.25 & 12/12 & 10 \\
03-31 & 5 & Cauchy: Huber IRLS & 13.38 & 15.55 & 12/12 & 6 \\
03-31 & 5 & Cauchy: Cauchy IRLS & 13.15 & 21.69 & 12/12 & 7 \\
03-31 & 5 & Cauchy: IGG~III IRLS & 12.13 & 16.21 & 12/12 & 4 \\
03-31 & 5 & Cauchy: Tukey IRLS & 15.15 & 18.81 & 12/12 & 6 \\
03-31 & 5 & Cauchy: IDS--ILS & 4.11 & 5.33 & 12/12 & 3 \\
\midrule
03-31 & 5 & No added error: WLS & 2.75 & 4.99 & 12/12 & 1 \\
03-31 & 5 & No added error: GMM, $K=2$ & 9.05 & 10.95 & 12/12 & 1 \\
03-31 & 5 & No added error: GMM, $K=3$ & 8.41 & 9.18 & 12/12 & 1 \\
03-31 & 5 & No added error: GMM, $K=4$ & 8.58 & 11.50 & 12/12 & 1 \\
03-31 & 5 & No added error: Huber IRLS & 8.06 & 15.32 & 12/12 & 3 \\
03-31 & 5 & No added error: Cauchy IRLS & 10.46 & 19.38 & 12/12 & 4 \\
03-31 & 5 & No added error: IGG~III IRLS & 8.24 & 13.36 & 12/12 & 3 \\
03-31 & 5 & No added error: Tukey IRLS & 9.36 & 15.36 & 12/12 & 4 \\
03-31 & 5 & No added error: IDS--ILS & 3.03 & 5.70 & 12/12 & 1 \\
\midrule
03-31 & 5 & Gaussian mixture: WLS & 3.77 & 5.96 & 12/12 & 1 \\
03-31 & 5 & Gaussian mixture: GMM, $K=2$ & 24.10 & 38.46 & 12/12 & 10 \\
03-31 & 5 & Gaussian mixture: GMM, $K=3$ & 23.73 & 51.65 & 12/12 & 10 \\
03-31 & 5 & Gaussian mixture: GMM, $K=4$ & 24.22 & 63.13 & 12/12 & 10 \\
03-31 & 5 & Gaussian mixture: Huber IRLS & 15.11 & 26.35 & 12/12 & 9 \\
03-31 & 5 & Gaussian mixture: Cauchy IRLS & 13.14 & 24.37 & 12/12 & 7 \\
03-31 & 5 & Gaussian mixture: IGG~III IRLS & 14.76 & 24.77 & 12/12 & 8 \\
03-31 & 5 & Gaussian mixture: Tukey IRLS & 17.90 & 20.38 & 12/12 & 8 \\
03-31 & 5 & Gaussian mixture: IDS--ILS & 3.93 & 4.42 & 12/12 & 4 \\
\midrule
03-31 & 10 & Cauchy: WLS & 3.16 & 4.09 & 12/12 & 1 \\
03-31 & 10 & Cauchy: GMM, $K=2$ & 50.91 & 57.46 & 12/12 & 20 \\
03-31 & 10 & Cauchy: GMM, $K=3$ & 56.85 & 62.64 & 12/12 & 20 \\
03-31 & 10 & Cauchy: GMM, $K=4$ & 38.54 & 67.01 & 12/12 & 13 \\
03-31 & 10 & Cauchy: Huber IRLS & 17.72 & 19.44 & 12/12 & 8 \\
03-31 & 10 & Cauchy: Cauchy IRLS & 17.45 & 24.42 & 12/12 & 8 \\
03-31 & 10 & Cauchy: IGG~III IRLS & 13.72 & 14.57 & 12/12 & 6 \\
03-31 & 10 & Cauchy: Tukey IRLS & 16.57 & 23.56 & 12/12 & 7 \\
03-31 & 10 & Cauchy: IDS--ILS & 4.98 & 7.37 & 12/12 & 6 \\
\midrule
03-31 & 10 & No added error: WLS & 3.40 & 4.09 & 12/12 & 1 \\
03-31 & 10 & No added error: GMM, $K=2$ & 17.44 & 53.61 & 12/12 & 4 \\
03-31 & 10 & No added error: GMM, $K=3$ & 16.61 & 62.90 & 12/12 & 4 \\
03-31 & 10 & No added error: GMM, $K=4$ & 19.60 & 63.08 & 12/12 & 4 \\
03-31 & 10 & No added error: Huber IRLS & 10.71 & 13.31 & 12/12 & 4 \\
03-31 & 10 & No added error: Cauchy IRLS & 12.18 & 16.11 & 12/12 & 5 \\
03-31 & 10 & No added error: IGG~III IRLS & 10.92 & 13.27 & 12/12 & 4 \\
03-31 & 10 & No added error: Tukey IRLS & 10.57 & 14.77 & 12/12 & 4 \\
03-31 & 10 & No added error: IDS--ILS & 3.71 & 5.66 & 12/12 & 1 \\
\midrule
03-31 & 10 & Gaussian mixture: WLS & 3.04 & 4.22 & 12/12 & 1 \\
03-31 & 10 & Gaussian mixture: GMM, $K=2$ & 49.27 & 67.19 & 12/12 & 20 \\
03-31 & 10 & Gaussian mixture: GMM, $K=3$ & 53.58 & 62.94 & 12/12 & 20 \\
03-31 & 10 & Gaussian mixture: GMM, $K=4$ & 58.28 & 68.22 & 12/12 & 20 \\
03-31 & 10 & Gaussian mixture: Huber IRLS & 15.36 & 18.58 & 12/12 & 7 \\
03-31 & 10 & Gaussian mixture: Cauchy IRLS & 15.15 & 16.60 & 12/12 & 6 \\
03-31 & 10 & Gaussian mixture: IGG~III IRLS & 14.22 & 16.52 & 12/12 & 6 \\
03-31 & 10 & Gaussian mixture: Tukey IRLS & 13.53 & 24.74 & 12/12 & 6 \\
03-31 & 10 & Gaussian mixture: IDS--ILS & 4.70 & 5.44 & 12/12 & 5 \\
\midrule
04-30 & 5 & Cauchy: WLS & 4.30 & 5.55 & 12/12 & 1 \\
04-30 & 5 & Cauchy: GMM, $K=2$ & 30.98 & 55.18 & 12/12 & 13 \\
04-30 & 5 & Cauchy: GMM, $K=3$ & 35.04 & 56.09 & 12/12 & 14 \\
04-30 & 5 & Cauchy: GMM, $K=4$ & 33.55 & 62.45 & 12/12 & 12 \\
04-30 & 5 & Cauchy: Huber IRLS & 18.46 & 21.21 & 12/12 & 7 \\
04-30 & 5 & Cauchy: Cauchy IRLS & 18.10 & 26.47 & 12/12 & 7 \\
04-30 & 5 & Cauchy: IGG~III IRLS & 16.82 & 22.21 & 12/12 & 6 \\
04-30 & 5 & Cauchy: Tukey IRLS & 17.03 & 36.59 & 12/12 & 6 \\
04-30 & 5 & Cauchy: IDS--ILS & 4.92 & 13.39 & 12/12 & 3 \\
\midrule
04-30 & 5 & No added error: WLS & 3.44 & 4.54 & 12/12 & 1 \\
04-30 & 5 & No added error: GMM, $K=2$ & 10.83 & 19.94 & 12/12 & 1 \\
04-30 & 5 & No added error: GMM, $K=3$ & 10.34 & 16.68 & 12/12 & 1 \\
04-30 & 5 & No added error: GMM, $K=4$ & 9.89 & 17.21 & 12/12 & 1 \\
04-30 & 5 & No added error: Huber IRLS & 11.88 & 14.51 & 12/12 & 4 \\
04-30 & 5 & No added error: Cauchy IRLS & 15.44 & 19.53 & 12/12 & 6 \\
04-30 & 5 & No added error: IGG~III IRLS & 11.42 & 13.01 & 12/12 & 3 \\
04-30 & 5 & No added error: Tukey IRLS & 13.19 & 16.53 & 12/12 & 4 \\
04-30 & 5 & No added error: IDS--ILS & 3.76 & 4.52 & 12/12 & 1 \\
\midrule
04-30 & 5 & Gaussian mixture: WLS & 3.88 & 6.26 & 12/12 & 1 \\
04-30 & 5 & Gaussian mixture: GMM, $K=2$ & 24.86 & 54.89 & 12/12 & 11 \\
04-30 & 5 & Gaussian mixture: GMM, $K=3$ & 11.74 & 49.87 & 12/12 & 2 \\
04-30 & 5 & Gaussian mixture: GMM, $K=4$ & 11.01 & 50.45 & 12/12 & 2 \\
04-30 & 5 & Gaussian mixture: Huber IRLS & 19.34 & 28.83 & 12/12 & 8 \\
04-30 & 5 & Gaussian mixture: Cauchy IRLS & 16.43 & 26.03 & 12/12 & 7 \\
04-30 & 5 & Gaussian mixture: IGG~III IRLS & 15.69 & 36.53 & 12/12 & 6 \\
04-30 & 5 & Gaussian mixture: Tukey IRLS & 19.38 & 21.43 & 12/12 & 8 \\
04-30 & 5 & Gaussian mixture: IDS--ILS & 4.52 & 7.18 & 12/12 & 4 \\
\midrule
04-30 & 10 & Cauchy: WLS & 4.13 & 5.47 & 12/12 & 1 \\
04-30 & 10 & Cauchy: GMM, $K=2$ & 60.47 & 81.44 & 12/12 & 20 \\
04-30 & 10 & Cauchy: GMM, $K=3$ & 38.53 & 64.70 & 12/12 & 11 \\
04-30 & 10 & Cauchy: GMM, $K=4$ & 40.22 & 72.69 & 12/12 & 11 \\
04-30 & 10 & Cauchy: Huber IRLS & 18.58 & 22.86 & 12/12 & 8 \\
04-30 & 10 & Cauchy: Cauchy IRLS & 18.73 & 22.83 & 12/12 & 7 \\
04-30 & 10 & Cauchy: IGG~III IRLS & 19.24 & 32.25 & 12/12 & 6 \\
04-30 & 10 & Cauchy: Tukey IRLS & 20.94 & 32.34 & 12/12 & 7 \\
04-30 & 10 & Cauchy: IDS--ILS & 5.02 & 5.99 & 12/12 & 5 \\
\midrule
04-30 & 10 & No added error: WLS & 3.95 & 6.14 & 12/12 & 1 \\
04-30 & 10 & No added error: GMM, $K=2$ & 17.34 & 22.44 & 12/12 & 1 \\
04-30 & 10 & No added error: GMM, $K=3$ & 18.05 & 26.77 & 12/12 & 1 \\
04-30 & 10 & No added error: GMM, $K=4$ & 17.18 & 25.06 & 12/12 & 1 \\
04-30 & 10 & No added error: Huber IRLS & 13.60 & 16.23 & 12/12 & 4 \\
04-30 & 10 & No added error: Cauchy IRLS & 16.69 & 21.59 & 12/12 & 6 \\
04-30 & 10 & No added error: IGG~III IRLS & 13.22 & 15.24 & 12/12 & 3 \\
04-30 & 10 & No added error: Tukey IRLS & 15.45 & 17.93 & 12/12 & 5 \\
04-30 & 10 & No added error: IDS--ILS & 4.42 & 6.12 & 12/12 & 1 \\
\midrule
04-30 & 10 & Gaussian mixture: WLS & 4.00 & 4.72 & 12/12 & 1 \\
04-30 & 10 & Gaussian mixture: GMM, $K=2$ & 53.42 & 75.84 & 12/12 & 20 \\
04-30 & 10 & Gaussian mixture: GMM, $K=3$ & 59.30 & 98.28 & 12/12 & 20 \\
04-30 & 10 & Gaussian mixture: GMM, $K=4$ & 65.68 & 97.67 & 12/12 & 20 \\
04-30 & 10 & Gaussian mixture: Huber IRLS & 13.54 & 21.65 & 12/12 & 5 \\
04-30 & 10 & Gaussian mixture: Cauchy IRLS & 15.79 & 23.13 & 12/12 & 6 \\
04-30 & 10 & Gaussian mixture: IGG~III IRLS & 15.43 & 19.42 & 12/12 & 6 \\
04-30 & 10 & Gaussian mixture: Tukey IRLS & 15.95 & 22.56 & 12/12 & 6 \\
04-30 & 10 & Gaussian mixture: IDS--ILS & 5.55 & 6.72 & 12/12 & 5 \\
\midrule
\bottomrule\end{longtable}

\FloatBarrier
\section{Representation and implementation sensitivity}
The whitening probes include 48 selected cases, including inputs without added errors. They reproduce all 240 independently checked reversed integer vectors. The maximum orthogonality error is \(6.0\times10^{-15}\), the maximum physical reconstruction error is \(1.3\times10^{-12}\,\mathrm m\), and the maximum relative WLS objective difference is \(5.2\times10^{-11}\). The reversed and cyclic permutations act within the code and phase blocks, as in Eq.~(26) of the main text for the reversed transform. Symmetric whitening uses \(\bm\Sigma^{-1/2}\bm L\) in the original order. Each comparison holds the physical injections and reference vectors fixed.

The electronic file \path{revision_20261103/third_round/results/whitening_recovery_intervals.csv} contains 3,240 intervals for both paired contrasts, three block lengths and both crossing-window rules. The 3,456 records in the adjacent \path{whitening_threshold_learning.csv} pair each estimate and denominator with its block setting. Thresholds are relearned by resampling calibration and evaluation dates separately. The records give the calibration counts, ties, unattainable targets and empty valid sets. The adjacent \path{labels/acceptance_by_label.csv} compares all nine configurations with fixed acceptance decisions and coordinates on windows with available reconstructed references. The archive also counts accepted windows excluded because a reconstructed reference is unavailable.
\begin{longtable}{p{20mm}p{28mm}rrr}
\caption{Paired full-vector recovery under alternative whitening. Vectors changed counts complete vectors that differ from the original-order Cholesky result among 3,120 matched window--law inputs. Recovery uses full-vector agreement with the internal reference. The mean change weights the eight date--length--law cells equally. A dash indicates the original reference representation.}\label{tab:representation-recovery}\\
\toprule
Whitening & Method & Recovery range (\%) & \shortstack{Vectors changed\\(count)} & \shortstack{Mean change\\vs original (pp)} \\
\midrule\endfirsthead
\multicolumn{5}{l}{\tablename~\thetable\ (continued)} \\
\toprule
Whitening & Method & Recovery range (\%) & \shortstack{Vectors changed\\(count)} & \shortstack{Mean change\\vs original (pp)} \\
\midrule\endhead
Original & WLS & 8.5--19.9 & 0 & -- \\
Original & Huber IRLS & 85.5--99.2 & 0 & -- \\
Original & Cauchy IRLS & 84.8--98.4 & 0 & -- \\
Original & IGG~III IRLS & 69.5--87.0 & 0 & -- \\
Original & Tukey IRLS & 77.1--94.3 & 0 & -- \\
Original & IDS--ILS & 66.9--93.5 & 0 & -- \\
Reversed & WLS & 8.5--19.9 & 0 & 0.0 \\
Reversed & Huber IRLS & 53.0--90.7 & 924 & -19.8 \\
Reversed & Cauchy IRLS & 55.1--91.9 & 883 & -17.7 \\
Reversed & IGG~III IRLS & 53.4--86.2 & 1215 & -8.3 \\
Reversed & Tukey IRLS & 56.8--90.7 & 998 & -11.1 \\
Reversed & IDS--ILS & 35.5--75.6 & 1362 & -25.7 \\
Symmetric & WLS & 8.5--19.9 & 0 & 0.0 \\
Symmetric & Huber IRLS & 67.5--95.1 & 566 & -10.4 \\
Symmetric & Cauchy IRLS & 70.7--96.4 & 555 & -8.3 \\
Symmetric & IGG~III IRLS & 61.8--88.3 & 999 & -5.1 \\
Symmetric & Tukey IRLS & 71.4--95.1 & 711 & -3.5 \\
Symmetric & IDS--ILS & 51.3--85.0 & 905 & -14.3 \\
\bottomrule
\end{longtable}

\begin{longtable}{lllrr}
\caption{Cross-date acceptance at the nominal 10\% target, with a separate threshold calibrated for each whitening representation on shared windows. Each range spans both date directions, lengths and injection laws. Extrema can come from different cells. The percentage displays are $100A=100n_a/n$, $100q=100n_w/n_a$ and $100V_{10}=100n_{a,10}/n$, where $n$ counts planned windows, $n_a$ acceptances, $n_w$ accepted wrong vectors and $n_{a,10}$ accepted baselines within 10\,cm. A displayed $q$ range of 0.0--0.0 means zero observed wrong fixes in all eight cells.}\label{tab:representation-acceptance}\\
\toprule
Whitening & Method & $A$ (\%) & $q$ (\%) & $V_{10}$ (\%) \\
\midrule\endfirsthead
\multicolumn{5}{l}{\tablename~\thetable\ (continued)} \\
\toprule
Whitening & Method & $A$ (\%) & $q$ (\%) & $V_{10}$ (\%) \\
\midrule\endhead
Original & WLS & 5.3--15.8 & 35.6--68.6 & 2.4--9.7 \\
Original & Huber IRLS & 7.9--12.6 & 0.0--0.0 & 7.9--12.6 \\
Original & Cauchy IRLS & 7.7--11.7 & 0.0--0.0 & 7.7--11.7 \\
Original & IGG~III IRLS & 7.1--11.8 & 0.0--2.6 & 7.1--11.8 \\
Original & Tukey IRLS & 7.9--11.0 & 0.0--2.4 & 7.9--11.0 \\
Original & IDS--ILS & 7.3--15.0 & 0.0--3.6 & 7.3--15.0 \\
Reversed & WLS & 5.3--15.8 & 35.6--68.6 & 2.4--9.7 \\
Reversed & Huber IRLS & 9.2--11.2 & 0.0--0.0 & 9.2--11.2 \\
Reversed & Cauchy IRLS & 8.1--12.1 & 0.0--0.0 & 8.1--12.1 \\
Reversed & IGG~III IRLS & 7.3--12.6 & 0.0--6.1 & 7.3--12.6 \\
Reversed & Tukey IRLS & 7.7--11.7 & 0.0--3.5 & 7.7--11.7 \\
Reversed & IDS--ILS & 5.7--16.6 & 0.0--18.8 & 5.7--15.8 \\
Symmetric & WLS & 5.3--15.8 & 35.6--68.6 & 2.4--9.7 \\
Symmetric & Huber IRLS & 8.1--12.6 & 0.0--0.0 & 8.1--12.6 \\
Symmetric & Cauchy IRLS & 6.9--13.4 & 0.0--0.0 & 6.9--13.4 \\
Symmetric & IGG~III IRLS & 5.3--15.4 & 0.0--0.0 & 5.3--15.4 \\
Symmetric & Tukey IRLS & 7.3--15.0 & 0.0--0.0 & 7.3--15.0 \\
Symmetric & IDS--ILS & 8.6--11.5 & 0.0--13.1 & 7.9--10.9 \\
\bottomrule
\end{longtable}

\begin{longtable}{lrrrrr}
\caption{IRLS sensitivity to initialization and iteration budget on 48 selected probes, 16 per input condition. Each rule contributes $48\times2\times2=192$ fits across cases, initializations and budgets. Row-floor activations are summed over weighted solves. At the 20-update budget, changing initialization alters 16 of the $48\times4=192$ paired case--rule vectors. Increasing the budget from 20 to 60 changes no vector in these probes.}\label{tab:implementation-ablation}\\
\toprule
Method & Correct & Recovery (\%) & \shortstack{Row-floor\\activations (count)} & \shortstack{Met stopping\\test (count)} & \shortstack{Reached iteration\\budget (count)} \\
\midrule\endfirsthead
\multicolumn{6}{l}{\tablename~\thetable\ (continued)} \\
\toprule
Method & Correct & Recovery (\%) & \shortstack{Row-floor\\activations (count)} & \shortstack{Met stopping\\test (count)} & \shortstack{Reached iteration\\budget (count)} \\
\midrule\endhead
Huber IRLS & 188 & 97.9 & 0 & 192 & 0 \\
Cauchy IRLS & 190 & 99.0 & 0 & 192 & 0 \\
IGG~III IRLS & 170 & 88.5 & 180 & 192 & 0 \\
Tukey IRLS & 178 & 92.7 & 110 & 192 & 0 \\
\bottomrule
\end{longtable}

\begingroup
\setlength{\LTcapwidth}{\linewidth}
\begin{longtable}{llrrrlrrrr}
\caption{Internal-reference reconstruction within the original accepted arc of each target--reference pair covering the central epoch (the later central epoch for even $N$). Exclude the window and five minutes on either side. The two-sided rule requires two eligible sides agreeing on the integer. The one-side-allowed rule also accepts a single eligible side. Complete vectors require a reference for each target satellite. The 256 windows comprise 64 per date and length, with 52 complete vectors by the two-sided rule and 245 when one side is allowed. None of the available complete vectors changes relative to its original reference.}\label{tab:label-sensitivity}\\
\toprule
 & & & & & & \multicolumn{4}{c}{Satellite-pair status counts} \\
Date & $N$ & Windows & \shortstack{Complete,\\two-sided} & \shortstack{Complete, one\\side allowed} & \shortstack{Changed,\\two/one} & Both & One & Conflict & Unavailable \\
\midrule\endfirsthead
\multicolumn{10}{l}{\tablename~\thetable\ (continued)} \\
\toprule
 & & & & & & \multicolumn{4}{c}{Satellite-pair status counts} \\
Date & $N$ & Windows & \shortstack{Complete,\\two-sided} & \shortstack{Complete, one\\side allowed} & \shortstack{Changed,\\two/one} & Both & One & Conflict & Unavailable \\
\midrule\endhead
03-31 & 5 & 64 & 12 & 61 & 0/0 & 388 & 86 & 0 & 3 \\*
03-31 & 10 & 64 & 9 & 60 & 0/0 & 384 & 94 & 0 & 4 \\*
04-30 & 5 & 64 & 16 & 63 & 0/0 & 391 & 94 & 0 & 1 \\*
04-30 & 10 & 64 & 15 & 61 & 0/0 & 381 & 100 & 0 & 3 \\
\bottomrule
\end{longtable}
\endgroup

In a descriptive comparison at the nominal 10\% target, 12 of the 16 nominal Wilson intervals for measured-observation cells with zero observed wrong fixes contain the simulation point estimate for the same method and law. The four exceptions are five-epoch Gaussian-mixture cells. Each method--law pair reuses its simulation estimate across both dates and record lengths.
\begin{longtable}{llllrrrl}
\caption{Descriptive comparison of zero-wrong-fix cells under original-order Cholesky whitening at the nominal 10\% target. Each interval runs from zero to the nominal Wilson upper endpoint under independent Bernoulli sampling. The simulation point estimate is for the same method and error law at the 10\% calibration target in main-text Table~9.}\label{tab:zero-event-comparison}\\
\toprule
Evaluation & $N$ & Law & Method & Accepted & \shortstack{Nominal Wilson\\upper endpoint (\%)} & \shortstack{Simulation\\$q$ (\%)} & \shortstack{Simulation estimate\\in interval?} \\
\midrule\endfirsthead
\multicolumn{8}{l}{\tablename~\thetable\ (continued)} \\
\toprule
Evaluation & $N$ & Law & Method & Accepted & \shortstack{Nominal Wilson\\upper endpoint (\%)} & \shortstack{Simulation\\$q$ (\%)} & \shortstack{Simulation estimate\\in interval?} \\
\midrule\endhead
04-30 & 5 & Gaussian mixture & Huber IRLS & 66 & 5.50 & 11.87 & no \\
04-30 & 5 & Gaussian mixture & Cauchy IRLS & 55 & 6.53 & 9.12 & no \\
04-30 & 5 & Cauchy & Huber IRLS & 67 & 5.42 & 5.33 & yes \\
04-30 & 5 & Cauchy & Cauchy IRLS & 62 & 5.83 & 3.63 & yes \\
04-30 & 10 & Gaussian mixture & Huber IRLS & 28 & 12.06 & 11.87 & yes \\
04-30 & 10 & Gaussian mixture & Cauchy IRLS & 28 & 12.06 & 9.12 & yes \\
04-30 & 10 & Cauchy & Huber IRLS & 21 & 15.46 & 5.33 & yes \\
04-30 & 10 & Cauchy & Cauchy IRLS & 23 & 14.31 & 3.63 & yes \\
03-31 & 5 & Gaussian mixture & Huber IRLS & 42 & 8.38 & 11.87 & no \\
03-31 & 5 & Gaussian mixture & Cauchy IRLS & 51 & 7.00 & 9.12 & no \\
03-31 & 5 & Cauchy & Huber IRLS & 42 & 8.38 & 5.33 & yes \\
03-31 & 5 & Cauchy & Cauchy IRLS & 41 & 8.57 & 3.63 & yes \\
03-31 & 10 & Gaussian mixture & Huber IRLS & 22 & 14.87 & 11.87 & yes \\
03-31 & 10 & Gaussian mixture & Cauchy IRLS & 20 & 16.11 & 9.12 & yes \\
03-31 & 10 & Cauchy & Huber IRLS & 29 & 11.70 & 5.33 & yes \\
03-31 & 10 & Cauchy & Cauchy IRLS & 29 & 11.70 & 3.63 & yes \\
\bottomrule
\end{longtable}

\begingroup
\small
\setlength{\tabcolsep}{4pt}
\setlength{\LTcapwidth}{0.9\linewidth}
\begin{longtable}{@{}ll*{8}{r}@{}}
\caption{Five-epoch Gaussian-mixture evaluation under alternative whitening representations. A separate threshold is calibrated for each representation and applied unchanged on the other date. WLS is shown once per direction because its outcomes coincide across representations. $A=n_a/n$ and $q=n_w/n_a$, where $n$ counts planned evaluation windows, $n_a$ acceptances and $n_w$ wrong accepted vectors. Both rates are in percent. A displayed $q$ of 0.00 means zero observed wrong fixes. The study reference for $q$ is 0.5\%.}\label{tab:whitening-cross-date-acceptance}\\
\toprule
Estimator & Whitening & \multicolumn{2}{c}{Nominal 5\%} & \multicolumn{2}{c}{Nominal 10\%} & \multicolumn{2}{c}{Nominal 25\%} & \multicolumn{2}{c}{Nominal 50\%} \\
\cmidrule(lr){3-4}\cmidrule(lr){5-6}\cmidrule(lr){7-8}\cmidrule(lr){9-10}
 &  & $A$ (\%) & $q$ (\%) & $A$ (\%) & $q$ (\%) & $A$ (\%) & $q$ (\%) & $A$ (\%) & $q$ (\%) \\
\midrule
\endfirsthead
\multicolumn{10}{@{}l}{\tablename~\thetable\ (continued)} \\
\toprule
Estimator & Whitening & \multicolumn{2}{c}{Nominal 5\%} & \multicolumn{2}{c}{Nominal 10\%} & \multicolumn{2}{c}{Nominal 25\%} & \multicolumn{2}{c}{Nominal 50\%} \\
\cmidrule(lr){3-4}\cmidrule(lr){5-6}\cmidrule(lr){7-8}\cmidrule(lr){9-10}
 &  & $A$ (\%) & $q$ (\%) & $A$ (\%) & $q$ (\%) & $A$ (\%) & $q$ (\%) & $A$ (\%) & $q$ (\%) \\
\midrule
\endhead
\bottomrule
\endlastfoot
\multicolumn{10}{@{}l}{\textbf{31 March calibration, 30 April evaluation}, $n=532$} \\*
WLS & All representations & 3.38 & 16.67 & 8.08 & 51.16 & 23.31 & 71.77 & 50.56 & 82.53 \\*
Huber IRLS & Original Cholesky & 6.58 & 0.00 & 12.41 & 0.00 & 25.38 & 0.00 & 47.56 & 0.00 \\*
Huber IRLS & Reversed Cholesky & 4.89 & 0.00 & 9.40 & 0.00 & 27.07 & 2.78 & 48.50 & 10.47 \\*
Huber IRLS & Symmetric & 6.39 & 0.00 & 9.77 & 0.00 & 26.50 & 0.00 & 48.50 & 1.94 \\*
Cauchy IRLS & Original Cholesky & 6.02 & 0.00 & 10.34 & 0.00 & 25.00 & 0.00 & 44.74 & 0.42 \\*
Cauchy IRLS & Reversed Cholesky & 5.45 & 0.00 & 10.15 & 0.00 & 23.31 & 3.23 & 48.68 & 7.72 \\*
Cauchy IRLS & Symmetric & 6.58 & 0.00 & 10.53 & 0.00 & 24.81 & 0.00 & 47.93 & 1.57 \\
\addlinespace[4pt]
\multicolumn{10}{@{}l}{\textbf{30 April calibration, 31 March evaluation}, $n=535$} \\*
WLS & All representations & 7.29 & 46.15 & 11.40 & 50.82 & 25.98 & 71.94 & 48.79 & 80.08 \\*
Huber IRLS & Original Cholesky & 2.43 & 0.00 & 7.85 & 0.00 & 24.49 & 0.00 & 51.78 & 0.00 \\*
Huber IRLS & Reversed Cholesky & 5.23 & 0.00 & 11.21 & 0.00 & 23.36 & 0.00 & 51.40 & 10.55 \\*
Huber IRLS & Symmetric & 4.11 & 0.00 & 11.59 & 0.00 & 24.49 & 0.00 & 51.78 & 2.53 \\*
Cauchy IRLS & Original Cholesky & 4.67 & 0.00 & 9.53 & 0.00 & 25.05 & 0.00 & 55.14 & 0.34 \\*
Cauchy IRLS & Reversed Cholesky & 4.30 & 0.00 & 9.91 & 0.00 & 25.98 & 0.72 & 52.52 & 7.47 \\*
Cauchy IRLS & Symmetric & 2.62 & 0.00 & 9.53 & 0.00 & 25.23 & 0.00 & 51.21 & 1.09 \\
\end{longtable}
\endgroup

\section{Additional method definitions and retained controls}
\normalsize
The main text retains the full estimator definitions and reports finite-budget results. The following details preserve the numerical conventions and auxiliary explanations without treating them as new performance evidence.
\subsection{Fixed-weight geometry and metric scale}
Consider the equivalent observation $\bm y_W$ and weight $\bm W_{\rm eq}$ for fixed geometry in~\eqref{supp:fixed-geometry-equivalent}. With the weight held fixed, choose a factor $\bm C$ satisfying $\bm W_{\rm eq}=\bm C^{\mathsf T}\bm C$. Write $\tilde{\bm y}=\bm C\bm y_W$, $\tilde{\bm A}=\bm C\bm A$ and $\tilde{\bm B}=\bm C\bm B$. For full-column-rank $\tilde{\bm A}$, $\bm P_\perp=\bm I-\tilde{\bm A}(\tilde{\bm A}^{\mathsf T}\tilde{\bm A})^{-1}\tilde{\bm A}^{\mathsf T}$ projects onto the orthogonal complement of its column space. Then
\begin{equation}
\bm Q_z^{-1}=\tilde{\bm B}^{\mathsf T}\bm P_\perp\tilde{\bm B},\qquad
\hat{\bm z}=\bm Q_z\tilde{\bm B}^{\mathsf T}\bm P_\perp\tilde{\bm y},
\label{eq:float}
\end{equation}
and the ranking of integer candidates depends on the data only through $\hat{\bm z}$ under this fixed metric. For an arbitrary integer candidate $\bm z$ and a nonzero integer displacement $\bm a$, the boundary between $\bm z$ and $\bm z+\bm a$ is the hyperplane
\begin{equation}
\bm a^{\mathsf T}\bm Q_z^{-1}\Bigl(\hat{\bm z}-\bm z-\frac{\bm a}{2}\Bigr)=0
\label{eq:boundary}
\end{equation}
in the float space. Changing the weights can move both the float solution and the decision boundaries. The integer estimate changes if the new float solution is assigned to a different integer under the new metric.

Give the contaminated block a fixed fraction $\varepsilon$ of its nominal weight and retain the nominal weights on the clean block. If the clean block has full column rank, the float solution and its working covariance tend to the clean-block values as $\varepsilon\to0$. When all fixed bias is confined to the downweighted block and the retained errors have zero mean, the limiting float bias vanishes. The working $\bm Q_z$ increases in the positive-semidefinite order as $\varepsilon$ decreases. Sampling covariance depends on the error model and requires separate evaluation.

In one dimension, let the float error be $\hat z-z_0$. Rounding is correct when that error is in $[-1/2,1/2]$, with ties of probability zero under a continuous model. For a normal float error with mean $\beta$ and standard deviation $\sigma$ and standard normal distribution function $\Phi$, the correct-rounding probability is
\begin{equation}
\Phi\Bigl(\frac{1/2-\beta}{\sigma}\Bigr)
-\Phi\Bigl(\frac{-1/2-\beta}{\sigma}\Bigr).
\label{eq:round}
\end{equation}
For a fixed \(\sigma\), a smaller \(|\beta|\) raises this probability. When \(|\beta|>1/2\), the mean already lies outside the correct rounding interval, and a moderate increase in \(\sigma\) can raise the probability over a finite range. A moderately biased but precise float solution can therefore outscore an unbiased but imprecise one.

The fixed-geometry sample-mean WLS implementation can use the nominal epoch covariance for the mean even when its sampling covariance is the epoch covariance divided by $N$. Positive common rescaling of the weight preserves the optimizer and changes the absolute covariance and objective scales. It preserves the distance ratio when the $10^{-15}$ denominator guard is inactive.
\subsection{GMM updates and implementation fallbacks}

The GMM models each epoch residual $\bm r_j=\bm y_j-\bm A\bm x-\bm B\bm z$ in fixed geometry by a zero-mean Gaussian mixture with $K$ components. The implementation uses non-negative mixing weights $\pi_k$ with initial values $1/K$. When $\sum_\ell\pi_\ell>0$, the mixture density is $\sum_{k=1}^K (\pi_k/\sum_\ell\pi_\ell)\varphi_m(\bm r;\bm0,\bm\Sigma_k)$, with $\varphi_m$ the $m$-variate Gaussian density. A common positive scale of the weights cancels in the epoch responsibility:
\begin{equation}
\gamma_{jk}=
\frac{\pi_k\,\varphi_m(\bm r_j;\bm 0,\bm\Sigma_k)}
{\sum_{\ell=1}^{K}\pi_\ell\,\varphi_m(\bm r_j;\bm 0,\bm\Sigma_\ell)}.
\label{eq:gamma}
\end{equation}
Responsibilities use the current residuals and component parameters. The scatter and regularized covariances are then updated, and the same responsibilities are combined with the updated covariance inverses for the next fit. With $\nu_k=\sum_j\gamma_{jk}$, define $\mathcal K_{\rm act}=\{k:\nu_k\ge10^{-8}\}$. For each active component, update
\begin{equation}
\bm S_k=\frac1{\nu_k}\sum_{j=1}^{N}\gamma_{jk}\,\bm r_j\bm r_j^{\mathsf T},
\qquad
\pi_k=\nu_k/N.
\label{eq:gmm-m}
\end{equation}
Inactive components receive zero mixing weight and are omitted from the precision update. The covariance update adds the trace-scaled ridge
\begin{equation}
\bm\Sigma_k^{\rm new}=\bm S_k+10^{-4}\max\!\left(\frac{\operatorname{tr}(\bm S_k)}{m},10^{-12}\right)\bm I_m,
\label{eq:gmm-ridge}
\end{equation}
followed by symmetrization. The ridge gives a positive-definite update in exact arithmetic. A direction absent from the scatter has variance $10^{-4}d_k$, where $d_k=\max(\operatorname{tr}(\bm S_k)/m,10^{-12})$. At $N=5$ and $m=16$, at least eleven directions obtain variance in this way. With fixed geometry and positive-definite $\bm W_{\rm eq}$, the equivalent weight is the sum of the epoch weights, and the equivalent observation is the weighted observation sum premultiplied by $\bm W_{\rm eq}^{-1}$:
\begin{equation}
\bm W_{\rm eq}=\sum_{k\in\mathcal K_{\rm act}}\nu_k(\bm\Sigma_k^{\rm new})^{-1},
\qquad
\bm y_W=\bm W_{\rm eq}^{-1}
\sum_{j=1}^{N}\sum_{k\in\mathcal K_{\rm act}}\gamma_{jk}(\bm\Sigma_k^{\rm new})^{-1}\bm y_j.
\label{eq:gmm-w}
\end{equation}
For time-varying geometry, use residuals $\bm r_j=\bm y_j-\bm A_j\bm x-\bm B_j\bm z$ and epoch weights $\bm W_j=\sum_{k\in\mathcal K_{\rm act}}\gamma_{jk}(\bm\Sigma_k^{\rm new})^{-1}$ in the stacked problem. When the sum of weighted densities for an epoch is zero, the implementation assigns responsibility $1/K$ to each component.

The component means are fixed at zero. The initial covariance for component $k$ is $(2k-1)\bm I_m$ in squared working-observation units. The fitted component covariances act on entire epoch residual vectors. The Gaussian-mixture law used to generate some inputs acts on individual coordinates.

For a common invertible coordinate transform, transforming the residual and every component covariance together multiplies all component densities by the same Jacobian factor. This factor cancels in the responsibilities. In the fitted estimator, the coordinate-wise median start, covariance pooling across epochs and trace-scaled isotropic ridge affect the metric. Epoch-dependent transforms can give different transformed covariances for different epochs. The responsibility identity applies at fixed component parameters under a common transform.

On the 40-epoch random benchmark, correctness indicators differ across GMM component counts on 11.8\% of Gaussian-mixture trials and 28.8\% of Cauchy trials. In the measured-observation sparse-injection experiment, at least two of the three GMM configurations return different complete vectors on 13 of 4,680 matched window--law inputs. Tables S9 and S10 give the short-record diagnostic results.

In the fixed-geometry path, an all-zero equivalent precision is replaced by $N\bm I_m$ and the equivalent observation is the arithmetic mean. In the stacked path, an epoch with zero total precision uses $\bm I_m$. The fixed-geometry path substitutes $10^{-100}$ for a component density evaluation that raises an exception. The inversion routines can use a pseudoinverse if ordinary inversion fails. Fallback activations were not exhaustively recorded. Section~7.7 describes the finite-budget stopping rule.
\subsection{Reference screening and persistent-bias design}
We construct long-arc reference integers separately from the short-window estimates. Arcs are split at gaps exceeding 60\,s or when successive residuals differ by more than 0.5 cycles. The latter condition flags suspected cycle slips. The fixed screening rules require at least 30 epochs and 900\,s duration, MAD at most 0.1 cycle, and at least 95\% of values within 0.25 cycle of the median. The full and half-arc medians must select the same candidate integer and lie within 0.25 cycle of it. Stable common bias can affect these internal labels. Arc screening accepts 706 of 742 arcs on 31 March and 724 of 772 on 30 April. Tables~S12 and~S13 give the subsequent window counts and arc screening outcomes. Chronological non-overlapping selection within date and length gives 535/532 five-epoch and 247/246 ten-epoch windows on 31 March/30 April. All 2,470 and 2,460 epochs selected at $N=10$ also appear at $N=5$ on the respective dates. The injection conditions reuse these windows.

In each evaluation window, add positive biases of 1.0\,m to code and 0.050\,m to phase on the two lowest-elevation target satellites for $d$ consecutive epochs. The start is $(i+2k)\bmod(N-d+1)$. Here $i$ is the window index and $k$ the selected-satellite index. Both indices and the returned start index are zero-based. The durations are 1, 2, 3 and 4 for $N=5$, and 2, 3, 4 and 6 for $N=10$. The code-only and phase-only controls use $d=2$ for $N=5$ and $d=4$ for $N=10$. With $d=4$ and $N=5$, WLS recovers 67.5\%/68.6\% and IGG~III recovers 53.8\%/56.2\% on 31 March/30 April. Main-text Section 4.3 reports this comparison.

\subsection{Separate Wilson-reference selection}
Form distinct thresholds from the 0th through the 99th sample percentiles of valid calibration ratios. Among thresholds with at least 100 calibration acceptances and a two-sided 95\% Wilson upper endpoint for $q$ at or below 0.5\%, select the one with the largest calibration accepted count. Freeze that threshold for evaluation. No threshold qualifies under Gaussian-mixture or Cauchy noise for any of the eight configurations. Evaluation endpoints are nominal independent-Bernoulli summaries.
\begin{table}[htbp]
\caption{Evaluation after threshold selection by the 0.5\% conditional Wilson rule. Each entry gives evaluation $100A$ / the nominal Wilson upper endpoint for evaluation $100q$, in percent. A dash indicates no qualifying calibration threshold.}\label{tab:accept}
\centering
\footnotesize
\setlength{\tabcolsep}{4pt}
\renewcommand{\arraystretch}{0.97}
\begin{tabular}{lrr}
\toprule
Estimator & Gaussian & Student-\(t_3\) \\
\midrule
WLS & 55.5/0.369 & 19.3/0.592 \\
GMM, \(K=2\) & 5.9/0.934 & -- \\
Huber IRLS & 56.4/0.394 & 35.0/0.470 \\
Cauchy IRLS & 55.5/0.417 & 37.1/0.540 \\
IGG~III IRLS & 40.5/0.582 & 9.0/0.499 \\
Tukey IRLS & 49.4/0.424 & 22.7/0.418 \\
\bottomrule
\end{tabular}
\end{table}

\footnotesize
\FloatBarrier

\normalsize
\subsection{Fixed-geometry compression of the weighted stack}\label{supp:fixed-geometry-compression}
The common observation and integer-search model in the main text uses every epoch explicitly. The compression below applies only when all epochs have the same design matrices $\bm A_j=\bm A$ and $\bm B_j=\bm B$, with the current positive-definite working weights $\bm W_j$ fixed. Let $\bm c=\bm A\bm x+\bm B\bm z$ be the common $m$-vector predicted at each epoch. Define
\begin{equation}
\bm W_{\rm eq}=\sum_{j=1}^{N}\bm W_j,\qquad
\bm y_W=\bm W_{\rm eq}^{-1}\sum_{j=1}^{N}\bm W_j\bm y_j,\qquad
\bm\Sigma_W=\bm W_{\rm eq}^{-1}.
\label{supp:fixed-geometry-equivalent}
\end{equation}
The matrix $\bm\Sigma_W$ is the working covariance of the equivalent observation. With independent epoch errors and weights equal to their true inverse covariances, it is also the sampling covariance of $\bm y_W$. Expanding the squares gives
\begin{equation}
\sum_{j=1}^{N}(\bm y_j-\bm c)^{\mathsf T}\bm W_j(\bm y_j-\bm c)
=(\bm y_W-\bm c)^{\mathsf T}\bm W_{\rm eq}(\bm y_W-\bm c)+C,
\label{supp:fixed-geometry-completion}
\end{equation}
where $C=\sum_j\bm y_j^{\mathsf T}\bm W_j\bm y_j-\bm y_W^{\mathsf T}\bm W_{\rm eq}\bm y_W$ is independent of $\bm x$ and $\bm z$. The full stack and equivalent problem have the same set of minimizers at frozen weights. With time-varying geometry, retain the full stack.

Section~7.1 gives the archived WLS scale convention and its relation to the denominator guard.

\subsection{Exact mean-preserving information augmentation}\label{supp:mean-preserving-augmentation}
At frozen weights, $y_{ij}$ is component $i$ of the observation vector at epoch $j$, $\bm d_{ij}^{\mathsf T}$ is row $i$ of $[\bm A_j\ \bm B_j]$, and $\bm\theta=(\bm x^{\mathsf T},\bm z^{\mathsf T})^{\mathsf T}$. The scalar residual is $r_{ij}=y_{ij}-\bm d_{ij}^{\mathsf T}\bm\theta$. The selected IRLS function supplies $w_{ij}$ from $\eta_{ij}=|r_{ij}|/\hat\sigma$, where $\hat\sigma$ is the robust scale for the current iteration, with a numerical floor of $10^{-12}$ applied before summation. The bound $\tilde a_i$ applies to the total weight across epochs, with $\ell=1,\ldots,m$. Define
\begin{equation}
a_i=\sum_{j=1}^{N}w_{ij},\qquad b_i=\sum_{j=1}^{N}w_{ij}y_{ij},\qquad
\tilde a_i=\max\!\left(a_i,10^{-3}\max_\ell a_\ell,10^{-8}\right),\qquad
\delta_i=\tilde a_i-a_i.
\label{supp:row-information-floor}
\end{equation}
The row means use the original total $a_i$:
\begin{equation}
\bar y_i=\frac{b_i}{a_i},\qquad
\bar{\bm d}_i=\frac{1}{a_i}\sum_{j=1}^{N}w_{ij}\bm d_{ij}.
\label{supp:weighted-row-means}
\end{equation}
For $J(\bm\theta)=\sum_{i,j}w_{ij}(y_{ij}-\bm d_{ij}^{\mathsf T}\bm\theta)^2$, the common augmented objective is
\begin{equation}
J_{\rm mean}(\bm\theta)=J(\bm\theta)+
\sum_i\delta_i(\bar y_i-\bar{\bm d}_i^{\mathsf T}\bm\theta)^2.
\label{supp:mean-preserving-objective}
\end{equation}
For each $\delta_i>0$, add an equation with response $\bar y_i$, design row $\bar{\bm d}_i^{\mathsf T}$ and weight $\delta_i$. In the stacked solver, multiply both the response and the design row by $\sqrt{\delta_i}$ and use unit weight for this equation. The original epoch rows remain in the fit.

For fixed geometry, $\bm d_{ij}=\bm d_i$ and
\begin{equation}
\sum_jw_{ij}(y_{ij}-\bm d_i^{\mathsf T}\bm\theta)^2
=a_i(\bar y_i-\bm d_i^{\mathsf T}\bm\theta)^2+C_i,
\label{supp:fixed-row-completion}
\end{equation}
Here $C_i$ is independent of $\bm\theta$. At fixed geometry the augmentation changes the coefficient from $a_i$ to $\tilde a_i$ at the same centre $\bar y_i=b_i/a_i$. The equivalent working variance is $1/\tilde a_i$. The historical augmentation rule computes the centre by dividing $b_i$ by $\tilde a_i$, changing the centre. Checks on 60 ordinary and 30 forced-floor cases agree between the compressed and stacked forms. Tables~S29--S30 report the historical sensitivities.

\subsection{Finite-budget stopping and retained weight diagnostics}\label{supp:finite-budget-diagnostics}
The starting fit uses the componentwise epoch medians and, for time-varying geometry, the epoch means of the design matrices. The outer iteration stops when $\|\bm x^{(t+1)}-\bm x^{(t)}\|_2<10^{-6}$ and the integer vectors are identical, or after 20 updates. The final diagnostic search uses the weights from the last update, including its augmentation when applicable. The reported finite-budget results include fits reaching the budget.

Among records with a valid stopping flag, the largest within-cell fractions of fits reaching the budget without satisfying the stopping test are 5.4\%, 6.3\%, 16.7\% and 12.3\% in the core grid, and 7.5\%, 6.2\%, 26.0\% and 18.0\% in the duration grid, for Huber, Cauchy, IGG~III and Tukey IRLS, respectively.

Across the four injection-law--length groups in the 48-case audit, IGG~III has a group mean of 15.5\% to 19.0\% for the within-window percentage of final weights at or below $10^{-11}$. Huber and Cauchy have zero such weights in these groups. Setting the augmentation response to zero changes the complete vector on five IGG~III and three Tukey cases out of 48 per rule. The Huber and Cauchy comparisons return the same vectors. Tables~S11 and~S24 report the weight and implementation diagnostics, respectively.

\footnotesize
\section{Historical configuration sensitivities}
The historical Huber, Cauchy and Tukey settings use the unscaled MAD. The historical centre rule divides the weighted observation sum by the floored total $\tilde a_i$. IGG~III uses standardized MAD. Tables~S29--S30 report these sensitivities separately from the comparisons using the principal protocol in main-text Section 3.1.
\begin{longtable}{llrrrrrrrr}
\caption{Full-vector recovery (\%) in the historical random-benchmark sensitivity under legacy weighting settings. The base uses $N=40$, $n_z=6$ and $\sigma=2$. Each setting changes the stated parameter and uses 1,000 paired problems per law.}\\
\toprule
 & & & \multicolumn{3}{c}{GMM components $K$} & & & & \\
Setting & Noise & WLS & 2 & 3 & 4 & \shortstack{Huber\\IRLS} & \shortstack{Cauchy\\IRLS} & \shortstack{IGG~III\\IRLS} & \shortstack{Tukey\\IRLS} \\
\midrule
\endfirsthead
\multicolumn{10}{l}{\tablename~\thetable\ (continued)} \\
\toprule
 & & & \multicolumn{3}{c}{GMM components $K$} & & & & \\
Setting & Noise & WLS & 2 & 3 & 4 & \shortstack{Huber\\IRLS} & \shortstack{Cauchy\\IRLS} & \shortstack{IGG~III\\IRLS} & \shortstack{Tukey\\IRLS} \\
\midrule
\endhead
$N=20$ & Cauchy & 0.2 & 58.3 & 67.0 & 71.4 & 84.1 & 84.2 & 86.2 & 85.3 \\
$N=20$ & Mixture & 0.7 & 34.1 & 54.8 & 56.1 & 59.2 & 59.4 & 72.2 & 66.9 \\
$N=80$ & Cauchy & 0.0 & 40.5 & 48.5 & 56.1 & 99.4 & 99.5 & 99.6 & 99.5 \\
$N=80$ & Mixture & 16.2 & 95.8 & 98.7 & 98.7 & 99.0 & 98.9 & 99.4 & 99.2 \\
base & Cauchy & 0.2 & 49.0 & 61.5 & 67.2 & 98.0 & 98.0 & 98.2 & 98.1 \\
base & Mixture & 4.1 & 75.6 & 87.4 & 87.4 & 88.7 & 88.5 & 94.3 & 91.5 \\
$n_z=4$ & Cauchy & 1.5 & 52.0 & 62.5 & 66.7 & 96.3 & 96.4 & 96.7 & 96.4 \\
$n_z=4$ & Mixture & 9.7 & 80.1 & 86.5 & 86.8 & 87.9 & 87.8 & 93.2 & 90.2 \\
$n_z=8$ & Cauchy & 0.0 & 51.4 & 63.1 & 69.1 & 98.6 & 98.6 & 99.0 & 98.9 \\
$n_z=8$ & Mixture & 1.8 & 73.1 & 89.9 & 90.6 & 91.5 & 91.4 & 96.4 & 93.4 \\
$\sigma=1$ & Cauchy & 2.9 & 91.0 & 93.0 & 95.7 & 100.0 & 100.0 & 100.0 & 100.0 \\
$\sigma=1$ & Mixture & 50.6 & 100.0 & 100.0 & 100.0 & 100.0 & 100.0 & 100.0 & 100.0 \\
$\sigma=4$ & Cauchy & 0.0 & 2.4 & 7.2 & 11.1 & 55.5 & 57.3 & 60.7 & 60.2 \\
$\sigma=4$ & Mixture & 0.0 & 0.6 & 10.4 & 20.6 & 24.4 & 25.8 & 36.5 & 32.6 \\
\bottomrule
\end{longtable}

\par\noindent\makebox[\linewidth][c]{%
\begin{minipage}{.90\linewidth}
\centering
\setlength{\abovecaptionskip}{4pt}
\setlength{\belowcaptionskip}{5pt}
\captionof{table}{Full-vector recovery (\%) by probability of drawing the high-variance Gaussian component under historical weighting settings. Each cell uses 1,000 problems. At 0\% and 20\%, the GMM inputs differ from the WLS and IRLS inputs. At the other probabilities all eight configurations share inputs. Huber, Cauchy and Tukey IRLS use the historical unscaled MAD.}\label{tab:historical-ratio}
\begin{tabular}{lcccccc}
\toprule
Estimator & 0\% & 5\% & 10\% & 20\% & 30\% & 50\% \\
\midrule
WLS & 98.4 & 46.1 & 17.6 & 4.1 & 1.4 & 0.4 \\
GMM, \(K=2\) & 95.5 & 93.5 & 91.4 & 72.0 & 66.6 & 39.9 \\
GMM, \(K=3\) & 95.5 & 93.8 & 92.5 & 85.4 & 69.2 & 23.0 \\
GMM, \(K=4\) & 95.5 & 93.7 & 92.6 & 85.4 & 75.0 & 21.2 \\
Huber IRLS & 96.0 & 94.0 & 93.8 & 88.7 & 78.1 & 40.0 \\
 Cauchy IRLS & 95.5 & 93.8 & 93.4 & 88.5 & 78.7 & 46.7 \\
IGG~III IRLS & 96.3 & 94.8 & 95.4 & 94.3 & 91.1 & 53.5 \\
Tukey IRLS & 95.5 & 93.8 & 93.9 & 91.5 & 85.8 & 55.1 \\
\bottomrule
\end{tabular}
\end{minipage}}\par

\FloatBarrier
\section*{Fixed-weight Gaussian bound}
This auxiliary diagnostic uses a separate Gaussian grid of 160 cells, crossing two geometries, one or two biased satellites, four durations of one to four out of five epochs, five bias amplitudes from $1\sigma$ to $5\sigma$, and two elevation classes. Each cell has 1,000 trials. After nominal whitening, WLS uses identity weight. Let $\bm A_{\rm st}$ and $\bm B_{\rm st}$ be the stacked real and integer designs, $\bm b$ the stacked fixed bias known in simulation, and $\bm P_\perp$ the orthogonal projector onto the complement of the column space of $\bm A_{\rm st}$. Write $\bm G=\bm B_{\rm st}^{\mathsf T}\bm P_\perp\bm B_{\rm st}$ and let $\bm z_0$ be the simulated true vector. Conditional on the full-rank design and fixed bias, $\hat{\bm z}-\bm z_0$ is Gaussian with mean $\bm G^{-1}\bm B_{\rm st}^{\mathsf T}\bm P_\perp\bm b$ and covariance $\bm G^{-1}$. For a nonzero integer displacement $\bm a$, define
\begin{equation}
h_{\bm a}=\frac12\sqrt{\bm a^{\mathsf T}\bm G\bm a},\qquad
d_{\bm a}=\frac{\bm a^{\mathsf T}\bm B_{\rm st}^{\mathsf T}\bm P_\perp\bm b}{\sqrt{\bm a^{\mathsf T}\bm G\bm a}}.
\end{equation}
Here $h_{\bm a}$ is the half-width of the pairwise decision interval in standard deviation units, and $d_{\bm a}$ is the standardized mean shift in the same direction. Let $\Phi$ be the standard normal distribution function. Conditional on the design and bias, the probability that $\bm z_0$ beats both $\bm z_0+\bm a$ and $\bm z_0-\bm a$ is $\Phi(h_{\bm a}-d_{\bm a})-\Phi(-h_{\bm a}-d_{\bm a})$. Let $p_{\min}$ be the minimum of these probabilities over the $3^{n_z}-1$ nonzero vectors with coordinates in $\{-1,0,1\}$, with $n_z=8$ here. A correct fix requires beating both candidates in every tested pair, so $p_{\min}$ bounds conditional WLS success from above.

The trial-specific bound is conditional on design and bias. Across the 160 cells, the rank correlation between cell-mean $p_{\min}$ and empirical WLS failure is $-0.99$. The cell-mean bound exceeds empirical WLS success by 0.9\,pp under favourable geometry and 9.7\,pp under weak geometry. Its scope is the fixed-weight Gaussian model.